\documentclass{article}
\usepackage{enumerate}
\usepackage{graphicx}
\usepackage{latexsym,amsmath,amsfonts,amscd, amsthm}
\usepackage{epsfig}
\usepackage{changebar}
\usepackage{color}
\usepackage{bm}
\usepackage{tikz}
\usepackage{multirow}
\usetikzlibrary{arrows,backgrounds,snakes}

\newtheoremstyle{plainNoItalics}{}{}{\normalfont}{}{\bfseries}{.}{ }{}

\theoremstyle{plain}
\newtheorem{thm}{Theorem}[section]

\theoremstyle{plainNoItalics}

\newtheorem{defn}[thm]{Definition}
\newtheorem{rem}[thm]{Remark}
\newtheorem{prop}[thm]{Proposition}
\newtheorem{exa}[thm]{Example}

\renewcommand{\theequation}{\thesection.\arabic{equation}}
\renewcommand{\thefigure}{\thesection.\arabic{figure}}
\renewcommand{\thetable}{\thesection.\arabic{table}}

\newcommand{\beq}{\begin{equation}}
\newcommand{\eeq}{\end{equation}}
\newcommand{\beqa}{\begin{eqnarray}}
\newcommand{\eeqa}{\end{eqnarray}}
\newcommand{\bit}{\begin{itemize}}
\newcommand{\eit}{\end{itemize}}
\newcommand{\bedef}{\begin{defn}}
\newcommand{\edefn}{\end{defn}}
\newcommand{\bpro}{\begin{prop}}
\newcommand{\epro}{\end{prop}}

\newcommand{\eps}{\varepsilon}

\newcommand{\bq}{{\bf q}}
\newcommand{\bu}{{\bf u}}

\newcommand{\bx}{{\bf x}}

\newcommand{\Dt}{{\Delta t}}
\newcommand{\R}{{\mathbb{R}}}
\newcommand{\Hf}{{\mathcal{H}}}

\newcommand{\ds}{\displaystyle}

\begin{document}

%\baselineskip=1.8pc

%\vspace*{.10in}

%%=============  title  =========================
\begin{center}
{
\bf High Order Semi-implicit WENO Schemes for All Mach Full Euler System of Gas Dynamics
}
\end{center}

\vspace{.2in}

%\author{Sebastiano Boscarino
%\thanks{Department of Mathematics and Computer Science, University of Catania, Catania, 95127, E-mail: {\tt boscarino@dmi.unict.it}. S. Boscarino and G.Russo are members of the INdAM Research group GNCS.
%They have been partially supported by ITN-ETN Horizon 2020 Project \emph{ModCompShock} (Modeling and Computation on Shocks and Interfaces) Project Reference 642768, by the Department of Mathematics and Computer Science, University of Catania \emph{Piano triennale della Ricerca 2016--2018}, and by INdAM-GNCS 2018 project {\em Numerical methods for multi-scale control problems and applications}.}, 
%Jing-Mei Qiu\thanks{Department of Mathematical Sciences, University of Delaware, Newark, DE, 19716, U.S.A.
% {\tt jingqiu@udel.edu}. 
% Research is partially supported by NSF grant NSF-DMS-1522777 and 1818924, Air Force Office of Scientific Computing FA9550-18-1-0257.},  
%Giovanni Russo \thanks{Department of Mathematics and Computer Science, University of Catania, Catania, 95125, E-mail: {\tt russo@dmi.unict.it.}}, 
%Tao Xiong \thanks{School of Mathematical Sciences, Fujian Provincial Key Laboratory of Mathematical Modeling and High-Performance Scientific Computing, Xiamen University, Xiamen, Fujian, P.R. China, 361005. Email: {\tt txiong@xmu.edu.cn.} Research is partially supported by NSFC grant 11601455, NSAF grant U1630247, NSF grant of Fujian Province No. 2016J05022 and the Fundamental Research Funds for the Central Universities No. 20720160009.
%}
%}

\centerline{ 
Sebastiano Boscarino \footnote{Department of Mathematics and Computer Science, University of Catania, Catania, 95127. E-mail: {\tt boscarino@dmi.unict.it}. S. Boscarino and G. Russo are members of the INdAM Research group GNCS. They would like to thank the Italian Ministry of Instruction, University and Research (MIUR) to support this research with funds coming from PRIN Project 2017 (No. 2017KKJP4X, entitled ``Innovative numerical methods for evolutionary partial differential equations and applications''). G. Russo has also been supported by ITN-ETN Horizon 2020 Project ModCompShock, Modeling and Computation on Shocks and Interfaces, Project Reference 642768 and S. Boscarino has been supported by the University of Catania (“Piano della Ricerca 2016/2018, Linea di intervento 2”).}, 
Jing-Mei Qiu\footnote{Department of Mathematics, University of Delaware,
Newark, DE 19716. E-mail: {\tt jingqiu@udel.edu}. Research is supported by NSF grant NSF-DMS-1818924, Air Force Office of Scientific Research FA9550-18-1-0257. },
Giovanni Russo \footnote{Department of Mathematics and Computer Science, University of Catania, Catania, 95125. E-mail: {\tt russo@dmi.unict.it.}}, 
Tao Xiong \footnote{Corresponding author. School of Mathematical Sciences, Fujian Provincial Key Laboratory of Mathematical Modeling and High-Performance Scientific Computing, Xiamen University, Xiamen, Fujian, P.R. China, 361005. Email: {\tt txiong@xmu.edu.cn.} Research is supported by NSFC No. 11971025, NSF of Fujian Province No. 2019J06002, the Strategic Priority Research Program of Chinese Academy of Sciences Grant No. XDA25010401, the Science Challenge Project No. TZ2016002.}
%\footnote{ S. Boscarino and G.Russo are members of the INdAM Research group GNCS.}
}

\bigskip

\noindent
{\bf Abstract.}
In this paper,  we propose {a new} high order semi-implicit scheme for the all Mach full Euler equations of gas dynamics. Material waves are treated explicitly, while acoustic waves are treated implicitly, thus avoiding severe CFL restrictions for low Mach flows. High order accuracy in time is obtained by semi-implicit temporal integrator based on the IMEX Runge-Kutta (IMEX-RK) framework. High order in space is achieved by finite difference WENO schemes with characteristic-wise reconstructions adapted to the semi-implicit IMEX-RK time discretization. {Type A IMEX schemes are constructed to handle not well-prepared initial conditions. Besides, these schemes are proven to be asymptotic preserving and asymptotically accurate as the Mach number vanishes for well-prepared initial conditions.}  {Divergence-free property of the time-discrete schemes is proved}. The proposed scheme can also well capture discontinuous solutions in the compressible regime, especially for two dimensional Riemann problems. Numerical tests in one and two space dimensions will illustrate the effectiveness of the proposed schemes.

\bigskip


\bigskip


%\newpage

\section{Introduction}
\label{sec:1}
\setcounter{equation}{0}
\setcounter{figure}{0}
\setcounter{table}{0}

Computational fluid dynamics (CFD) has been a very active research field in the past decades.
Numerical methods developed in this area generally can be divided into two categories, which are classified by the dimensionless Mach number. For moderate to high Mach number compressible effects have to be taken into account, while for low Mach number the flow can be considered incompressible or weakly compressible. For compressible flows, most numerical solutions are obtained by Godunov type shock capturing schemes for compressible Euler equations, which have the structure of a hyperbolic system of conservation laws \cite{leveque2002finite, toro2009riemann, GoRa2014, shu2009high, cockburn2006advanced}, while for the incompressible flows, preserving incompressibility and resolving vortex dynamics are among the main purposes \cite{chorin1968numerical, Temam1984, guermond2006overview}. 

There are, however, circumstances in which flows with a wide range of Mach number appear, making it desirable to develop numerical methods which can be applied for fluid flows at any speed, as already shown in the pioneering work of Harlow and Amdsen \cite{Harlow1968,Harlow1971}. However, due to the different physical mechanisms
and mathematical characteristics for the governing equations at different speeds, the development of efficient and effective numerical methods to capture flows with different compressibility is challenging \cite{Klein1995,Denner2020,park-munz-2005multiple-all-Mach}, and a lot of progress has been made only recently. For hyperbolic systems, waves propagate at finite speeds. Numerical methods have to resolve all the space and time scales that characterize these waves. Most shock capturing schemes devoted to such systems are obtained by explicit time discretization, and the time step has to satisfy a stability restriction, known as the CFL condition: it is limited by the size of the spatial mesh divided by the fastest wave speed. For compressible flows with Mach number greater than, say,  $0.1$, such a restriction is not a problem: indeed, if one is interested in resolving all the waves, accuracy and stability restrictions on space and time discretization are of similar nature. However, for low Mach flows, acoustic waves usually carry a negligible amount of energy. If one is not interested in resolving them, then the system becomes {\em stiff}: stability limitations on the time step are much stricter than the restrictions imposed by accuracy \cite{Turkel1987,Turkel1993}. In such cases, one may resort to implicit time discretization to avoid the acoustic CFL restriction. However, 
shock capturing schemes are highly non-linear, and a {\em naive\/} implicit version of them risks to be very inefficient. Furthermore, numerical viscosities for Godunov-type schemes are inversely proportional to the Mach number, introducing excessive numerical dissipations on the slow waves \cite{dellacherie2010analysis}. Preconditioning techniques are adopted to cure the large numerical diffusion as discussed in \cite{Turkel1987,viozat1997implicit,Ropke2015Mach}, but such techniques are effectively applicable only if Mach numbers are not too small.

On the other hand, as the Mach number vanishes, the flow converges to the incompressible limit. For the full Euler equations, at the incompressible limit, the density remains constant along the fluid particle trajectories, and the pressure waves propagate with infinite speed. Energy conservation equation reduces to the incompressibilty condition $\nabla\cdot \bu=0$ on the velocity field \cite{chorin1993}, so that the pressure and the density are decoupled. The pressure turns out to act as a Lagrange multiplier to enforce incompressibility of the flow \cite{Klein1995}. A rigorous proof for the compressible flow converging to the incompressible one as the Mach number goes to zero is given in \cite{klainerman1981singular}. 
%Numerically it is also important to preserve the incompressibility $\nabla\cdot \bu=0$ as in well-developed incompressible or weakly compressible fluid solvers \cite{colella1999projection, Gresho1990, guermond2006overview}. Among them, 
An effective approach to deal with low Mach flows is given by pressure-based algorithms, such as, for example, the one by Casulli and Greenspan \cite{casulli1984pressure}, in which a semi-implicit treatment of the pressure is incorporated in a scheme for compressible flow. The authors use an upwind discretization on the material wave, and an implicit equation for the pressure, which is solved by a SOR-type method. Several authors have subsequently worked on the development of semi-implicit methods \cite{Klein1995,munz2003extension-low-Mach} based on low-Mach asymptotics \cite{klainerman1981singular}. However, many of such schemes are specifically designed to deal with low Mach flows. When the fluid flow is compressible at large speed, shock discontinuities may form and propagate. In these cases, it is necessary to resort to conservative schemes (density-based schemes) which correctly capture possible shocks. 

Recently several papers have been written along these lines, see for example \cite{degond2011all, haack2012all, tang2012, dimarco2017study, BoRuSc18, boscarino2019high} for isentropic Euler and Navier-Stokes equations, or \cite{cordier2012asymptotic,tavelli2017pressure,Denner2018,Jonas2020novel,Denner2020, boscheri2021high} for full Euler and Navier-Stokes equations. However, most finite volume and finite difference schemes for full Euler equations are second order accurate in space and time, while existing high order schemes developed for  all Mach flows in isentropic Euler equations are not robust enough to be directly extended to the full Euler equations. 

The aim of the present paper is to design a {new} high order finite difference shock capturing scheme for the full compressible Euler equations. Finite difference weighted essentially non-oscillatory (WENO) schemes are used for spatial discretization, while high order semi-implicit IMEX Runge-Kutta (SI-IMEX-RK) methods are adopted for time discretization. {New IMEX schemes are suitably designed for stability and accuracy, with time stepping size independent of the Mach number $\varepsilon$.} 
% meanwhile it is a consistent and high order incompressible solver in the limit of vanishing $\varepsilon$. 
A key feature of the scheme is the implicit treatment of acoustic waves, while material waves are treated explicitly by WENO reconstructions of numerical fluxes. In particular, a suitable local Lax-Friedrich flux with characteristic-wise WENO reconstruction has been adopted for the explicit convective terms, while component-wise WENO reconstructions with zero numerical diffusion is used for implicit acoustic terms. {The method is able to capture shocks and discontinuities in an essentially non-oscillatory fashion in the compressible regime.}
%acoustic terms are treated implicitly with zero numerical diffusion and , thus achieving at the same time good accuracy, low dissipation \cite{boscarino2019high} and essentially non-oscillatory capture of solutions. 
In order to avoid the nonlinearity from the equation of state (EOS), a semi-implicit treatment similar to the one adopted in \cite{BoRuSc18} is used, leading to a linearized elliptic equation for the pressure, as described in Section \ref{HOrder}. {Another essential ingredient of the scheme design is to split the pressure into a thermodynamic pressure and a hydrodynamic one, using a similar idea adopted in \cite{cordier2012asymptotic} but with a fixed splitting parameter $\alpha$.} The thermodynamic pressure is used for the characteristic reconstructions, while the hydrodynamic pressure is obtained by solving an elliptic system. {We show that the resulting scheme is asymptotic preserving (AP) and asymptotically accurate \cite{jin1999efficient,jin2010asymptotic}, i.e., it is a consistent and high order discretization of the compressible Euler equations and, in the limit as $\varepsilon \to 0$, with $\Delta x$ and $\Delta t$ fixed, it becomes a consistent and high order discretization of the incompressible Euler equations.}

The rest of the paper is organised as follows. We recall the low Mach limit for the full compressible Euler equations in Section 2. We start Section 3 by introducing a first order semi-implicit scheme in time, then we describe the extension to high order time discretization using IMEX methods in Section 3.2, {in particular we design IMEX methods called of type A \cite{boscarino2008error}, which will be robust enough to solve the elliptic equation for the pressure}. We close the section with a description of high order spatial discretization obtained by characteristic-wise and component-wise WENO strategies.
The asymptotic preserving (AP) and asymptotic accuracy (AA) properties of the scheme are given in Section 4. Numerical tests are performed in Section 5, with the conclusion drew in the last section.

\section{Low Mach limit for the full Euler equations}
\label{Sec:2}

We consider the compressible Euler equations for an ideal gas in the non-dimensional form \cite{noelle2014weakly, dellacherie2010analysis, cordier2012asymptotic}:
\begin{equation}\label{Start1}
\left\{
\begin{array}{ll}
\rho_t + \nabla \cdot(\rho \mathbf{ u} ) \,=\, 0, \\ [3mm]
(\rho \mathbf{ u} )_t + \nabla \cdot(\rho \mathbf{ u}  \otimes \mathbf{ u} ) + \frac{1}{{\varepsilon^2}} {\nabla p} \,=\, 0, \\ [3mm]
E_t + \nabla\cdot[(E + p)\mathbf{ u} ] \,=\, 0,
\end{array}\right.
\end{equation}
with the EOS for a polytropic gas satisfying
\beq
\label{eq:EOS}
E = \frac{p}{\gamma -1} + \frac{\varepsilon^2}{2}\rho |\textbf{u}|^2,
\eeq
where $\gamma  > 1$ being the ratio of specific heats. 
The parameter $\varepsilon$ represents a global Mach number characterizing the non-dimensionalization. 
System (\ref{Start1}) is \emph{hyperbolic} and the eigenvalues along the direction {$\textbf{n}$} are: {$\lambda_1 = \textbf{u}\cdot \textbf{n} - c_s/\varepsilon$}, {$\lambda_2 = \textbf{u}\cdot \textbf{n}$}, {$\lambda_3 = \textbf{u}\cdot \textbf{n} + c_s/\varepsilon$} with $\displaystyle c_s = \sqrt{\gamma p/\rho}$.

When the reference Mach number is of order one, namely $\eps=\mathcal{O}(1)$, modern shock capturing methods are able to compute the formation and evolution of shocks and other complex structures with high resolutions at a reasonable cost. On the other hand, when the flows are slow compared to the speed of sound, i.e. $\eps\ll1$, we are near the incompressible regime. In such a situation, pressure waves become very fast compared to material waves. Standard explicit shock-capturing methods require a CFL time restriction dictated by the sound speed $c_s/\varepsilon$ to integrate the system. This leads to the stiffness in time, see e.g. \cite{haack2012all,degond2011all,cordier2012asymptotic}, where the time discretization is constrained by a stability condition given by
 $
   \Delta t < \Delta x/\lambda_{ \max} = \mathcal{O}(\varepsilon \Delta x),
$
here $\Delta t$ is the time step size,  $\Delta x$ is the mesh size and $\lambda_{\max} = \max_{\Omega} (|\textbf{u}|+c_s/\varepsilon)$ on the computational domain $\Omega$.
This restriction results in an increasingly large computational time for low Mach fluid flows. Moreover, excessive numerical viscosity (scales as $\varepsilon^{-1}$) in standard upwind schemes, leads to highly inaccurate solutions \cite{Turkel1987, Turkel1993}. Thus, it is of challenge and great importance to design numerical schemes, not only for shock-capturing, but also with consideration on stability and consistency in the incompressible limit, i.e. with \emph{asymptotic preserving} (AP) property.
In fact, in the low Mach limit, one is not interested in resolving the pressure waves; instead the fluid pressure serves as a Lagrangian multiplier in preserving the incompressibility of the velocity field. For the theoretical analysis of convergence from compressible flow to incompressible equations, such as $\eps\to0$, we refer to Klainerman and Majda \cite{klainerman1981singular, klainerman1982compressible} for a rigorous study in this low Mach limit. 

Here we recall the classical formal derivation of the incompressible Euler equations from the rescaled compressible Euler equations for an ideal gas (\ref{Start1}) with the EOS \eqref{eq:EOS}. 
We consider an asymptotic expansion \emph{ansatz} for the following two main variables:
\begin{equation}\label{ExEpBis}
\begin{array}{l}
p(\textbf{x},t) = p_0(\textbf{x},t) + 
\varepsilon^2 \, p_2(\textbf{x},t)+\cdots,\\ [3mm]
u(\textbf{x},t) = \bu_0(\textbf{x},t) + \varepsilon \, \bu_1(\textbf{x},t)+ \cdots,
\end{array}
\end{equation}
and insert them into the full Euler equations (\ref{Start1}). First, for the leading order $\mathcal{O}(\varepsilon^{-2})$, we formally find: $\nabla p_0(\bx,t) = 0$, i.e. pressure $p_0$ is constant in space, up to fluctuations of $\mathcal{O}(\varepsilon^2)$ and by \eqref{eq:EOS}, assuming for $E$ an expansion as $p$ in (\ref{ExEpBis}), we get $E_0 = p_0/(\gamma-1)$.
Then we formally find to $\mathcal{O}(1)$, 
\begin{equation}\label{eq: incomp}
\left\{
\begin{array}{l}
\displaystyle \partial_t \rho_0 + \nabla \cdot (\rho_0\bu_0) \,=\, 0, \\ [3mm]
\displaystyle \partial_t (\rho_0\bu_0) + \nabla \cdot \left( \rho_0 \bu_0 \otimes \bu_0 \right) + \nabla p_2 \,=\, 0,\\ [3mm]
\displaystyle \partial_t E_0 +  \nabla \cdot \left[ (E_0 + p_0) \bu_0\right] \,=\, 0.
\end{array}
\right.
\end{equation} 
From the energy equation in \eqref{eq: incomp}, with the EOS \eqref{eq:EOS} and denoting $d/dt = \partial/\partial t + \bu_0 \cdot \nabla$, we have
\beq
\label{eq: energy0}
\displaystyle \nabla \cdot  \bu_0 = -\frac{1}{ p_0 \gamma} \frac{d p_0}{dt}.
\eeq
 Now integrating the equation (\ref{eq: energy0}) over a spatial domain $\Omega$ with no-slip or periodic boundary conditions, we get $\int_{\Omega} \nabla \cdot \bu_0 \,d{\bf x} = \int_{\partial \Omega} \bu_0 \cdot {\bf n}\,ds= 0$, where ${\bf n}$ is the unit outward normal vector along $\partial \Omega$, and this implies $p_0$ is constant both in space and time, i.e. $p_0 = $ Const. 
Then, back into \eqref{eq: energy0}, one finds the divergence constraint $\nabla \cdot \bu_0 = 0$ in the zero Mach number limit. In summary, we have:
 \begin{equation}\label{incomp2_bis}
\left\{
\begin{array}{l}
p_0 = \textrm{Const.}, \quad \nabla \cdot {\bf \bu}_0 = 0,\\ [3mm]
\displaystyle \partial_t \rho_0 +\nabla \cdot ( \rho_0\bu_0)  = 0,\\ [3mm]
\displaystyle \partial_t (\rho_0\bu_0) + \nabla \cdot \left( \rho_0 \bu_0 \otimes \bu_0 \right) + \nabla p_2 = 0.
\end{array}
\right.
\end{equation}

As a direct consequence of $\nabla \cdot \bu_0 = 0$, we obtain from the mass-continuity
equation $d\rho/dt = 0$, namely the material derivative of the density is zero. This means that the density is constant along the particle trajectories. In particular, if the initial density is constant in space, the density of the fluid is constant in space and time.
Note that 
\begin{equation}
	\label{eq:p2}
p_2 = \lim_{\eps \to 0} \frac{1}{\eps^2} (p - p_0)
\end{equation}
is implicitly defined by the constraint $\nabla \cdot \textbf{u}_0 = 0$, which satisfies the following elliptic equation: 
\begin{equation}\label{equ_p_2}
\displaystyle -\nabla \cdot \left( \frac{1}{\rho_0} \nabla p_{2}\right) = \nabla \cdot (( \bu_0 \cdot \nabla) \bu_0).
\end{equation}

Finally, we assume the initial condition is {\em well-prepared} \cite{klainerman1981singular, klainerman1982compressible, dellacherie2010analysis, Klein1995},  that is, the initial condition for (\ref{ExEpBis}) is compatible with the equations at various orders of $\varepsilon$: 
\begin{equation}\label{ExEpBis0}
\left\{\begin{array}{l}
p(t = 0, \textbf{x}) = p_0 + \varepsilon^2 \,p_2(0,\textbf{x})+\cdots\\[3mm]
\bu(t = 0, \textbf{x}) = {\bu}_0(\textbf{x}) + \mathcal{O}(\varepsilon),
\end{array}\right.
\end{equation}
with $p_0  = \textrm{Const.}$ and $\nabla \cdot{\bu}_0 = 0$ and we impose $\rho(0,\bx) = \rho_0(\bx)$, with $\rho_0(\bx)$ being a strictly positive function. Note that well-prepared initial conditions are required if we want that the solution to the $\varepsilon$-dependent problem smoothly converges to the solution of the limiting incompressible problem. Furthermore, well-prepared initial condition is an important requirement 
to design AP schemes. It is crucial to preserve the constant state for leading order terms of $p$ and $E$, as well as that the divergence free constraint on the leading order term of $\bu$. For an arbitrary initial condition, an initial layer will appear, which requires a numerical resolution at the $\varepsilon$-scale.

\section{Numerical scheme}
\label{Sec:scheme}

In this section, we aim to construct and analyze a class of high order finite difference schemes with the AP property for unsteady compressible flows, when the Mach number $\varepsilon$ spans several orders of magnitude. The features of our scheme are the following: we design a semi-implicit IMEX (SI-IMEX) time discretization strategy, so that the scheme is stable with a time stepping constraint independent of the Mach number $\varepsilon$, the AP property is preserved in the zero Mach number limit and the scheme can be implemented in a semi-implicit manner \cite{boscarino2016high, boscarino2016linearly, boscarino2015linearly, boscarino2014high} to enable effective and efficient numerical implementations.  Our scheme preserves the incompressible velocity field in the zero Mach number limit by involving an elliptic solver for the hydrostatic pressure. 
In this section, the strategy of numerical discretizations is different from the traditional {\em method-of-lines} approach, since  we first perform time discretization to ensure AP property, after which we apply a suitable space discretization.
In particular we adopt high order WENO strategies with characteristic reconstructions tailored to IMEX-type methods in time. The final scheme can successfully capture shocks in the compressible regime, and efficiently solves the equations in the low Mach regime, with CFL condition depending only on fluid velocity.

\subsection{Semi-implicit treatment}
\label{sec: si}

We introduce our proposed strategy of implicit and explicit time discretizations, which is similar in spirit to the first order scheme in \cite{cordier2012asymptotic} with a slight modification. We emphasize the special treatment to avoid solving nonlinear equations required by a fully implicit scheme.
We rewrite (\ref{Start1}) as
\begin{equation}\label{Eeq2}
\frac{dU}{dt}= -\nabla \cdot \mathcal{F}_{E} - \nabla \cdot \mathcal{F}_{SI}
\end{equation}
where $U = (\rho, \rho{\bf u}, E)^T$ and 
\begin{equation}
  \label{Hfu}
\mathcal{F}_E \doteq \left( 
   \begin{array}{c}
   \displaystyle   \bq_E \\[3mm]
  \displaystyle  \left(\frac{\bq_E \otimes \bq_E}{\rho_E}\right) +{\alpha} \,p_E \mathbb{I}\\[3mm]
\displaystyle  0 
\end{array}
\right), 
\qquad
\mathcal{F}_{SI}\doteq
\left( 
   \begin{array}{c}
   \displaystyle  0 \\[3mm]
 \displaystyle   \frac{1-{\alpha}\varepsilon^2}{\varepsilon^2}  p_I \mathbb{I}\\[3mm]
\displaystyle  \frac{E_E + p_E}{\rho_I} \bq_I 
\end{array}
\right).
\end{equation}
Subscripts $E$ and $SI$ of $\mathcal{F}$ indicate the explicit and semi-implicit treatment of the first and the second term respectively. Several remarks are in order:
\begin{enumerate}
\item 
The parameter $\alpha$ determines the splitting between the explicit and implicit contribution of the pressure:
the former ensures that the explicit part is still hyperbolic, with a much smaller sound speed than the physical one, when $\varepsilon\ll 1$, while the latter will ensure much milder stability restrictions. As $\varepsilon$ increases, the explicit contribution becomes more and more relevant. This form of splitting is similar to the one in \cite{cordier2012asymptotic}, but differs in two main aspects: one is the implicit treatment of $\rho$ in the energy equation in (\ref{Hfu}), which gives a better asymptotic preserving property as will be elaborated in Section~\ref{sec: AP}; the other is the choice of the parameter $\alpha$ in splitting the pressure. In \cite{cordier2012asymptotic}, $\alpha$ is chosen depending on the Mach number $\varepsilon$. In our proposed scheme, however, $\alpha$ is chosen to be equal to 1 for all $\varepsilon < 1 $, and $\alpha = 1/\varepsilon^2$ for $\varepsilon \ge 1$.
\item 
We use the following EOS for (\ref{eq:EOS}) to avoid the nonlinearity in the semi-implicit scheme:
\begin{equation}\label{E_E}
\displaystyle E_E = \frac{1}{\gamma-1} p_E + \eps^2 \frac{|\bq_E|^2}{2\rho_E}, \quad
\displaystyle E_I = \frac{1}{\gamma-1} p_I + \eps^2 \frac{|\bq_E|^2}{2\rho_E}. 
\end{equation}
The subscripts $E$ and $I$ of $\bq$, $\rho$, $E$ and $p$ indicate the explicit and implicit treatments of the corresponding variables, respectively. 
\item 
{For the implicit term $p_I$ in \eqref{Hfu}, it is convenient to introduce a pressure perturbation $p_{I,2}$ \cite{park-munz-2005multiple-all-Mach}}, corresponding to the hydrodynamic pressure in the incompressible limit, defined as 
\beq
p_{I, 2} \doteq \frac{p_I - \bar{p}_E}{{\eps^2}},
\label{eq: pI2_a}
\eeq 
where $\bar{p}_E$ denotes the spatial average of the pressure $p_E$ computed from the EOS \eqref{E_E}.
In this way the term $p_{I, 2}$ will remain finite even as $\varepsilon\to 0$.
Then we have
\beq
\label{eq: p2}
\frac{1}{\eps^2}\nabla p_{I} = \nabla p_{I, 2}.
\eeq
\end{enumerate}
As an example, we present the scheme as well as the flow chart to update the numerical solution $U^{n+1} = (\rho^{n+1}, \bq^{n+1}, E^{n+1})^T$ for the first order semi-implicit scheme solving system (\ref{Start1}). We focus on the time discretization while keeping the space continuous, whose discretizations will be discussed in detail in Section \ref{ssec: spatial}:
 \begin{subequations}\label{eq:scheme_1bis_bis}
 \begin{align}
 \label{eq: rho_1}
&\frac{\rho^{n+1}-\rho^n}{\Dt} + \nabla \cdot \bq^{n} = 0,\\ \vspace{6mm}
\label{eq: q_1}
&\frac{\bq^{n+1}-\bq^n}{\Dt}  +\nabla \cdot \left(\frac{\bq^{n}  \otimes \bq^{n}}{\rho^{n}}+p^{n} \,{\mathbb I} \right) + \frac{1 - {\eps^2}}{\eps^2}\,\nabla p^{n+1} = 0,\\ \vspace{6mm}
\label{eq: E1}
&\frac{E^{n+1}-E^n}{\Dt} + \nabla\cdot \left(\frac{E^{n} + p^{n}}{\rho^{n+1}}\bq^{n+1} \right) =0,
\end{align}
 \end{subequations}
 with  
\begin{equation}
  \label{HfuNum}
\mathcal{F}_E (U^n)\doteq \left( 
   \begin{array}{c}
   \displaystyle   \bq^n \\[3mm]
  \displaystyle  \left(\frac{\bq^n \otimes \bq^n}{\rho^n}\right) +  p^n \, \mathbb{I}\\[3mm]
\displaystyle  0 
\end{array}
\right), 
\quad
\mathcal{F}_{SI}(U^n, U^{n+1})\doteq
\left( 
\begin{array}{c}
   \displaystyle  0 \\[3mm]
   \displaystyle  \frac{1-\varepsilon^2}{\varepsilon^2}  \,p^{n+1} \,\mathbb{I}\\[3mm]
   \displaystyle  \frac{E^n + p^n}{\rho^{n+1}} \bq^{n+1} 
\end{array}
\right).
\end{equation}
%Note that the pressure is split into an explicit pressure $p_E=p^n$ and an implicit one $p_I=p^{n+1}$ in (\ref{HfuNum}). 
The flow chart based on the semi-implicit scheme is the following:
%Below is the procedure to update solutions. Note that 
\begin{enumerate}
\item update $\rho^{n+1}$ from \eqref{eq: rho_1}.  
\item we rewrite \eqref{eq: q_1} as
\beq
\label{eq: scheme_1bis_bb}
\bq^{n+1}   = {\bq^*} - {\Dt} \frac{1 - {\eps^2}}{\eps^2}\nabla p^{n+1},
\eeq
with ${\bq^*} = {\bq^n} - {\Dt}\nabla \cdot \left(\frac{\bq^{n}  \otimes \bq^{n}}{\rho^{n}}+p^{n} {\mathbb I}\right) $. We substitute $\bq^{n+1}$ into \eqref{eq: E1} to get
\begin{equation} \label{eq: E1_bis}
E^{n+1} = E^* + {\Dt}^2\frac{1 - {\eps^2}}{\eps^2}\nabla\cdot\left(H^n \nabla p^{n+1}\right),
\end{equation}
where $H^n  = (E^{n} + p^{n})/\rho^{n+1}$ and $E^* = E^n - {\Dt}\nabla\cdot (H^n{\bq^*})$.
\item Now we replace $E^{n+1}$ by $p^{n+1}/(\gamma-1) + \eps^2 {(\bq^n)^2}/{(2\rho^n)}$ in \eqref{eq: E1_bis} using \eqref{E_E}, 
and with the introduction of \eqref{eq: pI2_a} i.e.
$p^{n+1}_{2} \doteq ({p^{n+1} - \bar{p}^n})/{{\eps^2}}$ where $p_I = p^{n+1}$ and $\bar{p}^{n}=\bar{p}_E$,  we rewrite \eqref{eq: E1_bis} as  
\beq
\label{eq: scheme_1bis_bbb}
\frac{\eps^2}{\gamma -1}p^{n+1}_2 - {\Dt}^2(1 - \eps^2) \nabla \cdot \left(H^n \nabla p_2^{n+1}\right) = E^{**},
\eeq
where $E^{**} = E^* -  \bar{p}^n/({\gamma -1}) -\eps^2 {(\bq^n)^2}/({2\rho^n})$ is explicitly computed. 
We obtain an elliptic equation (\ref{eq: scheme_1bis_bbb}) for $p^{n+1}_{2}$. 
\item Finally, we update $\bq^{n+1}$ from \eqref{eq: scheme_1bis_bb}, and then $E^{n+1}$ from \eqref{eq: E1}. 
\end{enumerate}
Note that if $\varepsilon\ge1$ the implicit pressure contribution in \eqref{eq:  q_1} vanishes, so the momentum $\bq^{n+1}$ is evaluated explicitly. With updated $\rho^{n+1}$ and $\bq^{n+1}$, $E^{n+1}$ can also be updated in an explicit way from \eqref{eq: E1}.

%, and the whole scheme becomes explicit.
%The above mentioned treatments will be generalized to high order temporal discretization in the next subsection. 
%
%  
\subsection{High order semi-implicit (SI) temporal discretization using IMEX}
\label{HOrder}

\subsubsection{High order SI-IMEX R-K scheme for all Mach number full Euler system}

We generalize the first order scheme (\ref{eq:scheme_1bis_bis}) to high order in the framework of IMEX R-K methods. In order to do so, 
we follow the idea introduced in \cite{boscarino2016high}. We consider an autonomous system of the form
\[
U' = \Hf(U,U), \quad U(t_0) = U_0,
\]
%we re-write system (\ref{Start1}) as an autonomous system following \cite{boscarino2016high}
%\[
%U' = \Hf(U,U), \quad U(t_0) = U_0,
%\]
where $\mathcal{H}: \R^n \times \R^n \to \R^n $ is a sufficiently regular mapping. We assume an explicit treatment of the first argument of $\Hf$ (using subscript "E"), and an implicit treatment to the second argument (using subscript "I" ), i.e.
\label{subsection: h-time}
\begin{equation}\label{PartitionedSyst}
\left\{
\begin{array}{l}
U'_E =  \mathcal{H}(U_E,U_I),\\[3mm]
U'_I =  \mathcal{H}(U_E,U_I),
\end{array}
\right.
\end{equation}   
 with initial conditions 
\begin{equation}\label{InitialCondU}
 U_E(t_0) = U_0, \quad U_I(t_0) = U_0.
\end{equation} 
Then system (\ref{Eeq2}) with (\ref{Hfu}) can be rewritten in the form (\ref{PartitionedSyst}) where $U \doteq (U_E, U_I)^T$ with $U_E = (\rho_E, \bq_E, E_E)^T$, and  $U_I = (\rho_I, \bq_I, E_I)^T$, and $\mathcal{H}(U_E,U_I) = -\nabla \cdot \mathcal{F}_{E} - \nabla \cdot \mathcal{F}_{SI}$.

System (\ref{PartitionedSyst}) is a particular case of the \emph{partitioned system}, \cite{hairer1993solving2}.
One can apply an IMEX R-K scheme to (\ref{PartitionedSyst}), using the corresponding pair of Butcher  $tableau$ \cite{butcher2016},
\begin{equation}\label{DBT}
\begin{array}{c|c}
\tilde{c} & \tilde{A}\\
\hline
\vspace{-0.25cm}
\\
 & \tilde{b^T} \end{array}, \ \ \ \ \  
\begin{array}{c|c}
{c} & {A}\\
\hline
\vspace{-0.25cm}
\\
 & {b^T} \end{array},
\end{equation}
where $\tilde{A} = (\tilde{a}_{ij})$ is an $s \times s$ matrix for an explicit scheme, with $\tilde{a}_{ij}=0$ for $j \geq i$ and $A = ({a}_{ij})$ is an  $s \times s$ matrix for an implicit scheme. For the implicit part of the methods, we use a diagonally implicit scheme, i.e. $a_{ij}=0$, for $j > i$, in order to guarantee simplicity and efficiency in solving the algebraic equations corresponding to the implicit part of the discretization. The vectors $\tilde{c}=(\tilde{c}_1,...,\tilde{c}_s)^T$, $\tilde{b}=(\tilde{b}_1,...,\tilde{b}_s)^T$, and $c=(c_1,...,c_s)^T$, $b=(b_1,...,b_s)^T$ complete the characterization of the scheme. The coefficients $\tilde{c}$ and $c$ are given by the usual relation
\begin{eqnarray}\label{eq:candc}
\tilde{c}_i = \sum_{j=1}^{i-1} \tilde a_{ij}, \ \ \ c_i = \sum_{j=1}^{i} a_{ij}.
\end{eqnarray}
{From now on, it is useful to characterize different IMEX schemes we will consider in the sequel accordingly to the structure of the DIRK method. Following \cite{boscarino2008error} we call: an IMEX-RK method of type A if the matrix $A \in \mathbb{R}^{s \times s}$  is invertible, and we call an IMEX-RK method of type CK if the matrix $A$ can be written as
\[
A = \left(
\begin{array}{cc}
0 & 0\\
a & \hat{A}
\end{array}
\right),
\] 
with $a = (a_{21},...,a_{s1})^T \in \mathbb{R}^{(s-1)}$ and the submatrix $\hat{A} \in \mathbb{R}^{(s-1) \times (s-1)}$ is invertible, or equivalently $a_{ii} \neq 0$, $i = 2,...,s$. In the special case $a = 0$, $b_1 = 0$ the scheme is said to be of type ARS and the DIRK method is reducible to a method using $s-1$ stages. 
}
Later, for the consideration of AP property, we consider {\em stiffly accurate} (SA) implicit schemes, i.e., the implicit part of the Butcher table satisfies  
the condition $b^T = e^T_s A$, with $e_s = (0,...,0,1)$ and $c_s = 1$. We will see that SA guarantees that the numerical solution is identical to the last internal stage value of the scheme.
%\end{defn}

Now an SI-IMEX  R-K scheme applied to (\ref{PartitionedSyst}) reads
\begin{subequations}
\label{IMEXschem}
	\begin{equation}\label{IMEXschem1}
	U^{(i)}_E = U^n_E + \Delta t \sum_{j= 1}^{i-1} \tilde{a}_{ij}\mathcal{H}(U^{(j)}_E, U^{(j)}_{I}), \quad
	U^{(i)}_I = U^n_I + \Delta t \sum_{j= 1}^{i} {a}_{ij}\mathcal{H}(U^{(j)}_E, U^{(j)}_{I}),
	\end{equation}
	\begin{equation}\label{LIbis}
	U^{n+1}_E = U^n_E + \Delta t \sum_{i= 1}^{s} \tilde{b}_{i}\mathcal{H}(U^{(i)}_E, U^{(i)}_{I}),\quad
	U^{n+1}_I = U^n_I + \Delta t \sum_{i= 1}^{s} {b}_{i}\mathcal{H}(U^{(i)}_E, U^{(i)}_{I}).
	\end{equation} 
\end{subequations}
	 Apparently, by writing in the form (\ref{PartitionedSyst}), we increased the computational cost since we double the number of variables.
	 However, the cost is mainly due to the number of function evaluations, which essentially depends on the semi-implicit scheme, because in (\ref{IMEXschem}) the identical term 
	 $\mathcal{H}(U^j_E, U^j_{I})$ appears in both explicit and implicit parts. Furthermore,  
	  if $\tilde{b}_i = b_i$ for all $i$, then $U^{n}_E = U^n_I$ for all $n> 0$ (provided $U^{0}_E = U^0_I$) and therefore the duplication of variables is not necessary. Such a property is of particular relevance for designing time discretization schemes which are asymptotic preserving.

	 We rewrite the scheme \eqref{IMEXschem} in the following new form for the convenience of further discussion
%	 \begin{subequations}\label{IMEXscheme}
	 \begin{equation}\label{IMEXscheme1}
	U^{(i)}_E = U^n + \Delta t \sum_{j= 1}^{i-1} \tilde{a}_{ij}\mathcal{H}(U^{(j)}_E, U^{(j)}_{I}), \quad 	{U}_{I}^{(i)}   = \tilde{U}_I^{(i)}  +\Delta t\, {a}_{ii}\mathcal{H}(U^{(i)} _E, U^{(i)}_{I}), 
	\end{equation}
%	\begin{equation}\label{IMEXscheme2}
where $\tilde{U}_I^{(i)} = U^n + \Delta t \sum_{j= 1}^{i-1} {a}_{ij}\mathcal{H}(U^{(j)}_E, U^{(j)}_{I})$. Denoting $U^n=U^n_E=U^n_I$ for all $n>0$ (due to $\tilde{b}_i = b_i$ for all $i$), finally we have 
%	\end{equation}
	\begin{equation}\label{LI}
	U^{n+1} = U^n + \Delta t \sum_{i= 1}^{s} {b}_{i}\mathcal{H}(U^{(i)} _E, U^{(i)}_{I}).
	\end{equation}
%	\end{subequations}
Similar to the first order scheme, we split the pressure into explicit ($p_E$) and implicit ($p_I$) terms. In order to avoid the nonlinearity in the semi-implicit step we use the following equations of state \eqref{E_E}, i.e., for the generic stages $i=1,\ldots,s$
\begin{subequations}\label{E_I}
\beq\label{E_I1}
\displaystyle E^{(i)}_E  = \frac{1}{\gamma-1} p^{(i)}_E + \eps^2 \frac{|\bq^{(i)}_E|^2}{2\rho^{(i)}_E}, 
\eeq
\beq\label{E_I2}
\displaystyle E^{(i)}_I = \frac{1}{\gamma-1} p^{(i)}_I + \eps^2 \frac{|\bq^{(i)}_E|^2}{2\rho^{(i)}_E}.
\eeq
\end{subequations} 

\begin{rem}
The first order scheme \eqref{eq:scheme_1bis_bis} is the same as applying the following Butcher table to \eqref{PartitionedSyst} 
\begin{equation}\label{Afirst}
\mbox{\bf Type A I:} \quad 
\begin{array}{c|c}
              0 & 0 \\
               \hline
              &1 
\end{array} \qquad
\begin{array}{c|c}
                1 & 1     \\
              \hline 
               &  1
\end{array}. 
\end{equation}
Formally applying the above tableau (\ref{Afirst})  to the partitioned system (\ref{PartitionedSyst}), it reads 
\begin{equation*}
\begin{array}{ll}
U^{(1)}_E &= U^n,\\[2mm]
U^{(1)}_I &= U^n + \Delta t \mathcal{H}(U^{(1)}_{E}, U^{(1)}_{I}),\\[2mm]
U^{n+1} &= U^n + \Delta t \mathcal{H}(U^{(1)}_{E}, U^{(1)}_{I}), 
\end{array}
\end{equation*}
which is the first order scheme \eqref{eq:scheme_1bis_bis}, considering $U^{(1)}_E = U^n$ and $U^{(1)}_I = U^{n+1}$. 
\end{rem}

\begin{rem}
{The authors in \cite{boscheri2021high} proposed a different semi-implicit discretization of the EOS (\ref{E_I2}) where the kinetic energy in the total energy definition splits into an explicit and an implicit contribution, namely
\begin{equation}\label{E_I2_bis}
E^{(i)}_I  = \frac{1}{\gamma-1} p^{(i)}_I + \eps^2 \frac{|(\bq^{(i)}_E)^T \bq^{(i)}_I |}{2\rho^{(i)}_I}. 
\end{equation}
We performed numerical tests to compare such semi-implicit treatment of the kinetic energy (\ref{E_I2_bis}) with (\ref{E_I2}). It is found numerically that such semi-implicit treatment (\ref{E_I2_bis}) may produce slightly better error in some test cases. Overall comparable performances are observed. In order to save space, we decide not to report these results in this paper.}
\end{rem}

\subsubsection{Construction of an IMEX R-K scheme}
Next, we construct an IMEX R-K Butcher tableau for semi-implicit discretization of all-Mach full Euler system. The construction is based on the following considerations for accuracy and for handling non-well prepared initial data.  
\begin{enumerate}
\item 
We require that the implicit part of the IMEX R-K scheme is SA. 
With such an assumption, one can derive the AP and AA properties of the scheme as discussed in Section~\ref{sec: AP}. Note that if the implicit part of the scheme is $A$-stable, SA is a sufficient condition to make it $L$-stable, see \cite{hairer1993solving2}. 
\item The invertibility of the implicit matrix $A$ of the SI-IMEX R-K scheme is important to handle non-well prepared initial conditions and for proper initialization of the hydrodynamic pressure. {In particular, $a_{11} \neq 0$ for IMEX R-K schemes of type A. This is critical to solve the pressure wave equation (\ref{eq: scheme_1bis_bbb}).  If $a_{11} = 0$ (e.g. in an IMEX scheme of  type ARS), then we can not have the pressure wave equation, hence no proper value of $p_{I, 2}$ (e.g. see eq.~\eqref{eq34:LIbis}) at the first IMEX stage.}
%In fact, when the initial condition is not well-prepared, there is usually an initial layer for which an invertible matrix $A$ ensures robustness of an IMEX RK scheme. 
For further related discussions and analysis, we refer the reader to the following papers \cite{pareschi2005implicit, boscarino2013implicit, boscarino2017asymptotic}. {Thus, we construct IMEX R-K of type $A$ for the time discretization.}
\item To synchronize $U_E$ and $U_I$, the weights for the final stage of the double tableaus should be the same, i.e. $\tilde{b}_i = b_i$, $i=1\cdots s$. That is, we keep only one set of numerical solution in the process of updating \cite{boscarino2016high}.
Alternatively, we can select a different vector of weights for the $U_E(t)$, say $\tilde{b}_i \neq b_i$, which will provide a lower/higher order approximation of the solution for $U_E(t)$; this can be used to implement a procedure of automatic time step control \cite{hairer1993solving1}. 
\item
In order to simplify order conditions, given the equations of state used in the semi-implicit step \eqref{E_I}, 
%to keep the correct evaluations of the EOS  \eqref{E_I} at the same time level, 
we impose that 
%\beq
%\label{eq: c}
$c_i =\tilde{c}_i, \, i = 2 \cdots s$.
%\eeq 
Note that if matrix $A$ is invertible we have 
$c_1 \neq \tilde{c}_1=0$ \cite{pareschi2005implicit}. 
\end{enumerate}
Based on the above considerations, we design a high-order IMEX scheme with matrix $A$ invertible, which is SA and satisfies $c_i  = \tilde{c}_i$ for $i = 2,\cdots,s$. We impose the conditions required for a third order scheme \cite{pareschi2005implicit} and obtain the following double Butcher tableau
\begin{align} \label{IMEX1_(4,4,3)}
 &\textrm{\bf Explicit :} \nonumber \\ \vspace{8mm}
 &\begin{array}{c|cccc}
                    0 & 0 & 0 & 0 & 0 \\
                    \gamma & \gamma & 0 & 0 & 0\\
                   0.717933260754 & 0.435866521508& 0.282066739245 & 0 & 0\\
                    1 &   -0.733534082748750 & 2.150527381100 &  -0.416993298352& 0\\
               \hline 
                  0 & 0 & 1.208496649176& -0.644363170684 & \gamma
\end{array},\nonumber \\ \vspace{8mm}
&\textrm{\bf Implicit :}  \\ \vspace{8mm}
&\begin{array}{c|cccc}
               \gamma&  \gamma & 0  & 0 & 0\\
              \gamma & 0&  \gamma & 0 & 0\\
               0.717933260754 &0 & 0.282066739245 & \gamma & 0\\
               1 &0 & 1.208496649176& -0.644363170684 & \gamma\\
               \hline 
                &0 & 1.208496649176& -0.644363170684 & \gamma
\end{array},\nonumber
\end{align}
with $\gamma = 0.435866521508$. We call it SI-IMEX(4,4,3), the triplet $(s, \sigma, p)$ characterizing the number of stages of the implicit scheme ($s = 4$), the number of stages of the explicit scheme ($\sigma=4$ ) and the order of the scheme ($p = 3$). % \cite{ascher1997implicit}.
{A corresponding 2nd order scheme, which we call SI-IMEX(3,3,2) is given as
	\begin{align} \label{IMEX_(3,3,2)}
	&\textrm{\bf Explicit :} \hspace{3.4cm}\textrm{\bf Implicit :}  \nonumber \\ \vspace{8mm}
	&\begin{array}{c|ccc}
	0 & 0 & 0 & 0  \\
	\gamma & \gamma & 0 & 0 \\
	1 & \delta & 1-\delta & 0 \\
	\hline 
	0 & 0 & 1-\gamma & \gamma
	\end{array}, \hspace{2cm}
	\begin{array}{c|ccc}
	\gamma&  \gamma & 0  & 0 \\
	\gamma & 0&  \gamma & 0 \\
	1 &0 & 1-\gamma&  \gamma\\
	\hline 
	0& 0& 1-\gamma & \gamma
	\end{array},
	\end{align}
	where $\gamma=1-\sqrt{2}/2$ and $\delta=-2\sqrt{2}/3$.
}

%%%%%%%%%%%%%%%%%
\subsection{High order spatial discretization for all-Mach fluid flows}
\label{ssec: spatial}
Below we describe our spatial discretization strategies that incorporate WENO mechanism to capture shocks in the compressible regime and produce a high order incompressible solver for the flow in the zero Mach limit.  
%In particular, the WENO spatial discretization will be performed for $\nabla \cdot \mathcal{F}_E$ in \eqref{Hfu} along local characteristic directions, and $\nabla \cdot \mathcal{F}_{SI}$ in \eqref{Hfu} will be handled in the spirit of an incompressible Euler solver. 
One major difficulty, hence the new ingredient in the scheme design, in extending the high order AP scheme for the isentropic Euler system \cite{boscarino2019high} to the full Euler system is about the pressure. In the isentropic case, the pressure is an explicit function of $\rho$, whereas in the full Euler case the pressure comes from the EOS involving all conserved variables. As such, we split the pressure into the part involving $p_E$ obtained from the EOS \eqref{E_I1}, and the part $p_{I, 2}$ obtained by formulating an elliptic equation from a semi-implicit solver of the system. In this process, the spatial discretization in the scheme formulation becomes critical for its robustness. There are several ingredients in our spatial discretization that differentiate our approach from the existing ones in the literature. Firstly, we carefully apply the WENO procedure for spatial reconstructions of fluxes for the full Euler system which is a nontrivial generalization from the isentropic one \cite{boscarino2019high}. In particular, we propose to apply the fifth order characteristic-wise WENO procedure to $\mathcal{F}_E$ (in \eqref{Hfu}) of the full Euler system so that in the compressible regime oscillations could be best controlled; and to apply a component-wise WENO procedure for the flux functions of $\mathcal{F}_{SI}$ (in \eqref{Hfu}). We found that such a choice is optimal in compromising the need from both compressible and incompressible fluid solvers. The characteristic-wise WENO turns out to be important for compressible Euler system, as shown in the numerical section. Secondly, following the strategy in our previous work \cite{boscarino2019high}, we apply a compact high order spatial discretization to second order derivatives terms in \eqref{eq33:LIbis}. This is discussed in Section~\ref{ssec: flowchart}. 
%, to be elaborated in the 

\subsubsection{WENO spatial discretization}
\label{ssec: weno}
Without loss of generality, we describe our algorithms in a 2D setting. %In a 2D setting, 
Let 
$$U = (\rho, \rho u, \rho v, E)^T,$$
$$\nabla \cdot \mathcal{F}_E = \partial_x \mathcal{F}^x_E +\partial_y \mathcal{F}^y_E,
\quad
\nabla \cdot \mathcal{F}_{SI} = \partial_x \mathcal{F}^x_{SI} +\partial_y \mathcal{F}^y_{SI}
$$
with flux functions in $x$ and $y$ directions 
$$
\mathcal{F}^x_E = (\rho u, \rho u^2+p, \rho u v,  0)^T, \quad 
\mathcal{F}^y_E = (\rho v, \rho u v, \rho v^2+p, 0)^T,$$ 
$$
\mathcal{F}^x_{SI} = \left(0, \frac{1-\eps^2}{\eps^2}p, 0, \frac{E+p} {\rho} (\rho u) \right)^T, \quad 
\mathcal{F}^y_{SI} = \left(0, 0, \frac{1-\eps^2}{\eps^2}p, \frac{E+p} {\rho} (\rho v)\right)^T.$$ 
We consider a rectangular domain discretized by a uniform Cartesian grid, with meshsize $\Delta x = \Delta y$, and grid points $(x_i, y_j)$, $i=1, \cdots N_x$, $j=1\cdots N_y$, located at cell centers. We use the subscript $\cdot_{i,j}$ to denote the solution point values at $(x_i, y_j)$, and $\cdot_{i\pm1/2, j}$ and $\cdot_{i, j\pm1/2}$ for the reconstructed numerical fluxes in approximating $x$- and $y$- derivatives respectively. 

Below we first present two types of WENO spatial discretizations, i.e, a characteristic-wise WENO and a component-wise WENO, for explicit and semi-implicit parts as specified in \eqref{Hfu}, respectively. 
\begin{enumerate}
\item {\em Characteristic-wise WENO $\nabla_{CW}$ with Lax-Friedrichs splitting for discretizing $\nabla \cdot \mathcal{F}_E$.} 
For compressible hyperbolic systems with shocks, WENO reconstructions in the component-wise fashion may still lead to oscillations and a local characteristic decomposition will be needed \cite{shu1998essentially}. 
We let $\mathcal{F}^{x}_{E, i, j}$ be the flux function at each grid point $(x_i, y_j)$. 
To approximate the $\partial_x \mathcal{F}^{x}_{E}|_{x_i, y_j}$, we fix a $j$-index and perform a global Lax-Friedrichs splitting for the flux terms, i.e. 
\[
\mathcal{F}^{x}_{E, i, j} = \mathcal{F}^{x, +}_{E, i, j}+\mathcal{F}^{x, -}_{E, i, j},\; \forall i
\]
with 
\beq
\label{eq: lxf}
\mathcal{F}^{x, \pm}_{E, i, j} =\frac12 (\mathcal{F}^{x}_{E, i, j} \pm \Lambda \,U_{E,i,j}), \quad \Lambda = \max_{u,v}\{|u|+|v|+\min\{1, \frac1\eps\}\,c_s\}
\eeq
%\textcolor{magenta}{Are we sure it is global Lax-Friedrichs? It should be local... Also, the method becomes not practical if $\varepsilon>1$, so we may assume $\varepsilon<1$ and specify that in case $\varepsilon\ge 1$ one should use just an explicit solver}
with $c_s = \sqrt{\gamma p/\rho}$ being the sound speed multipled by $\varepsilon$, and the $\max$ is taken over appropriate range of $u$, $v$. {Notice that such a choice of $\Lambda$ is the same as in our previous work \cite{boscarino2019high}.} We then project $\{\mathcal{F}^{x, \pm}_{E, i, j}\}_{i=1}^{N_x}$ to local characteristics directions {for the full Euler system in the compressible regime (corresponding to $\varepsilon=1$)}, perform WENO reconstruction of fluxes there and then project the reconstructed fluxes back, to obtain the flux terms as
$$\hat{\mathcal{F}}^{x}_{E, i+\frac12, j}=\hat{\mathcal{F}}^{x, +}_{E, i+\frac12, j}+\hat{\mathcal{F}}^{x, -}_{E, i+\frac12, j}.$$
Then 
\[
\partial_x \mathcal{F}^{x}_{E}|_{x_i, y_j} \approx \frac{1}{\Delta x} (\hat{\mathcal{F}}^{x}_{E, i+\frac12, j}-\hat{\mathcal{F}}^{x}_{E, i-\frac12, j}).
\]
%\textcolor{magenta}{We may consider putting the description of the CW WENO procedure in an Appendx, if space allows it} 
Similar procedure could be performed to approximate spatial derivatives in the $y$-direction. 
We introduce the notation of $\nabla_{CW}$ for the characteristic-wise WENO spatial discretization of $\mathcal{F}_E$. 
\item {\em A component-wise WENO $\nabla_{W}$ with Lax-Friedrichs flux splitting for discretizing $\nabla \cdot \mathcal{F}_{SI}$}. We apply a component-wise WENO reconstruction of fluxes for discretizing $\nabla \cdot \mathcal{F}_{SI}$ denoted as $\nabla_W \cdot \mathcal{F}_{SI}$. In particular, a Lax-Friedrichs flux splitting, followed by a component-wise WENO reconstruction is performed to obtain the numerical fluxes. 
\end{enumerate}
%{In addition, following the compact discretization for second order operator proposed in \cite{boscarino2019high}, }
%\SB{Note that in the evaluation of the momentum component in $U^{(i)}$ (\ref{IMEXscheme2}), the derivative term $\nabla p^{(i)}_I$ in $\mathcal{F}_{SI}$ will be discretized by a central difference approximation denoted as $\nabla_C$ and for the energy equation we split the operator $\nabla$ as $\nabla_{W}$ for the explicit quantities and $\nabla_{C}$ for the evaluating part ???? .}

%%%%%%%%%%%%%%%%%

\subsubsection{Flowchart and compact discretization of second order derivatives}
\label{ssec: flowchart}
With the introduction of characteristic-wise and component-wise WENO $\nabla_{CW}$ and $\nabla_W$, we summarize the flow chart of the SI-IMEX R-K time discretization (\ref{IMEXscheme1})-(\ref{LI}), coupled with WENO spatial discretization, for solving all-Mach full Euler system. In this process, a compact high order discretization for second order spatial derivatives in \eqref{eq33:LIbis} is applied.  
 \begin{enumerate}
\item Start from $U^n$ at time $t^n$, we first compute $U_E^{(i)}$  from (\ref{IMEXscheme1})
   \begin{subequations}
   	\begin{equation}\label{rhoE}
       \displaystyle \rho^{(i)}_E  = \rho^n - \Delta t\, \sum_{j = 1}^{i-1} \tilde{a}_{ij} \nabla_{CW} \cdot \bq_E^{(j)},  \\%[2mm]
               \end{equation}
       \begin{equation}\label{QE}
       \displaystyle {\bq}^{(i)}_E = \bq^n  - {\Delta t}\,\sum_{j = 1}^{i-1}\tilde{a}_{ij}\left(\nabla_{CW} \cdot \left(\frac{\bq^{(j)}_E \otimes \bq^{(j)}_E}{\rho^{(j)}_E} + p^{(j)}_E \,{\mathbb I}\right) +  (1-\varepsilon^2) \frac{\nabla_W p_{I}^{(j)}}{\varepsilon^2}\right), \\%[2mm]
               \end{equation}
       \begin{equation}\label{EE}
       \displaystyle  {E}^{(i)}_E = E^n - \Delta t \,\sum_{j = 1}^{i-1}\tilde{a}_{ij} \nabla_W \cdot \left( \bar{H}^{(j)}{\bq}^{(j)}_I \right).
               \end{equation}
      \end{subequations}
\item We compute $\tilde{U}^{(i)}$ in (\ref{IMEXscheme1}).%(\ref{LIbis2})
   \begin{subequations}
       \begin{equation}\label{rhoI}
       \displaystyle \tilde{\rho}^{(i)}  = \rho^n - \Delta t\, \sum_{j = 1}^{i-1} {a}_{ij} \nabla_{CW} \cdot \bq_E^{(j)},%\\[2mm]
       \end{equation}
        \begin{equation}\label{QI}
       \displaystyle \tilde{\bq}^{(i)}= \bq^n  - {\Delta t}\,\sum_{j = 1}^{i-1}{a}_{ij}\left(\nabla_{CW} \cdot \left(\frac{\bq^{(j)}_E \otimes \bq^{(j)}_E}{\rho^{(j)}_E} + p^{(j)}_E \,{\mathbb I}\right) + (1-\varepsilon^2)\frac{\nabla_W p_{I}^{(j)}}{\varepsilon^2}\right), %\\[2mm]
        \end{equation}
        \begin{equation}\label{EexI}
       \displaystyle  \tilde{E}^{(i)} = E^n - \Delta t \,\sum_{j = 1}^{i-1}a_{ij} \nabla_W \cdot \left( \bar{H}^{(j)}{\bq}^{(j)}_I \right),
        \end{equation}
      \end{subequations}
      with $\ds \bar{H}^{(j)} =  (E^{(j)}_E+ p^{(j)}_E)/\rho^{(j)}_I$. 
\item Solve $U_I^{(i)}$ from (\ref{IMEXscheme1}). 
\begin{enumerate}
\item In components, $U_I^{(i)}$ satisfies
	\begin{subequations}\label{eq: LIbis}
	\begin{equation} \label{eq1:LIbis}
      \displaystyle {\rho}^{(i)}_I  = \displaystyle \tilde{\rho}^{(i)} -  \Delta t\, a_{ii} \nabla_{CW} \cdot \bq^{(i)}_E, 
      \end{equation}
      \begin{equation}\label{eq2:LIbis}
     \displaystyle  {\bq}^{(i)}_I = \displaystyle \tilde{\tilde{\bq}}^{(i)} - \Delta t \, a_{ii}\left(  
      \frac{1-\varepsilon^2}{\varepsilon^2} \nabla p_{I}^{(i)}\right), 
      \end{equation}
             \begin{equation}\label{eq3:LIbis}
     \displaystyle E_I^{(i)} = \tilde{E}^{(i)}  -\Delta t a_{ii} \nabla \cdot (\bar{H}_{i}\bq_I^{(i)}),
     %  (1-\varepsilon^2)\Delta t^2 \, a_{ii}^2\nabla \left( \bar{H}^{(i)}  \cdot \left(      
%      \frac{\nabla_C p_{I}^{(i)}}{\varepsilon^2} \right) \right), 
           \end{equation}
     \end{subequations}
    where
      \begin{equation}
    \tilde{\tilde{\bq}}^{(i)} = \displaystyle \tilde{\bq}^{(i)} - \Delta t \, a_{ii}\left(\nabla_{CW} \cdot \left(
         \frac{\bq^{(i)}_E \otimes \bq^{(i)}_E}{\rho^{(i)}_E} + p^{(i)}_E \,{\mathbb I}\right) \right). 
	\end{equation}	
%         and with
        %\tilde{\tilde{{E}}}^{(i)} = \displaystyle \tilde{E}^{(i)} - \Delta t \,a_{ii} \left(\nabla_W \cdot \left( \bar{H}^{(i)}\tilde{\tilde{{\bq}}}^{(i)}_I \right)\right).
      %\end{equation}
    %  The Lax-Friedrich parameter $\Lambda$ for the discrete operator $\nabla_W$ in \eqref{Ett} for the momentum equation  is taken to be $0$, and is taken to be $\max\{|u|+|v|+\min\{1, \frac1\eps\} c_s\}$ for the energy equation \eqref{EexI} and \eqref{eq3:LIbis}. 
\item 
To solve the system \eqref{eq: LIbis}, we substitute ${\bq}^{(i)}_I$ of \eqref{eq2:LIbis} into \eqref{eq3:LIbis} and obtain
             \begin{equation}\label{eq33:LIbis}
     \displaystyle E_I^{(i)} = E^{**}_I  +  (1-\varepsilon^2)\Delta t^2 \, a_{ii}^2\nabla \left( \bar{H}^{(i)}  \left(      
      \frac{\nabla p_{I}^{(i)}}{\varepsilon^2} \right) \right), 
           \end{equation}
           with $E_I^{(i)}$ following the EOS \eqref{E_I2} and
%    \begin{equation}\label{Ett}
$E_I^{**}  = \tilde{E}^{(i)} - \Delta t\, a_{ii} \nabla_W \cdot \left(\bar{H}^{(i)}\tilde{\tilde{\bq}}^{(i)} \right)$.
%	\end{equation}			
\eqref{eq33:LIbis} is an implicit equation about $p^{(i)}_{I}$. By the introduction of pressure perturbation $p^{(i)}_{I, 2}$ in \eqref{eq: pI2_a}, we solve $p^{(i)}_{I, 2}$ from 
\begin{equation}\label{eq34:LIbis}
\displaystyle \frac{\eps^2}{\gamma-1}\, p_{I,2}^{(i)} = E^{***}_I  +  (1-\varepsilon^2)\Delta t^2 \, a_{ii}^2\nabla \left( \bar{H}^{(i)}  \left(      
\nabla p_{I,2}^{(i)} \right) \right), 
\end{equation}
where $E^{***}_I=E^{**}_I-\bar{p}_E/(\gamma-1)-\eps^2|\bq^{(i)}_E|^2/(2\rho^{(i)}_E)$.
Notice that in the process of substitution to obtain \eqref{eq33:LIbis} or \eqref{eq34:LIbis}, the gradient and the divergence operators are kept to be continuous, obtaining the second order operator $\nabla \left( \bar{H}^{(i)} \nabla p_{I,2}^{(i)} \right)$. The second order spatial derivative is then discretized by a compact discretization as proposed in \cite{boscarino2019high}. 
\item 
With $p^{(i)}_{I,2}$ solved from \eqref{eq34:LIbis}, we update $\bq^{(i)}_I$ from 
      \begin{equation}\label{eq22:LIbis}
     \displaystyle  {\bq}^{(i)}_I = \displaystyle \tilde{\tilde{\bq}}^{(i)} - \Delta t \, a_{ii}\left(  
      1-\varepsilon^2\right) \nabla_W p_{I,2}^{(i)}, 
      \end{equation}
and successively update $E^{(i)}_I$ from
             \begin{equation}
     \displaystyle E_I^{(i)} = \tilde{E}^{(i)}  -\Delta t a_{ii} \nabla_{W} \cdot (\bar{H}_{i}\bq_I^{(i)}). 
           \end{equation} 
\end{enumerate}
\item Finally, update the numerical solution $U^{n+1} = U^{(s)}_I$ with the assumption on the SA property of IMEX R-K schemes.
\end{enumerate}
This completes the description of the high order SI-IMEX R-K time discretization to the all-Mach full Euler equations.

\section{ Asymptotic preserving (AP) and Asymptotically Accurate (AA) property}
\label{sec: AP} 
\subsection{AP property}
In this section we prove the AP property of scheme (\ref{eq:scheme_1bis_bis}). 
%It is straightforward to extend this AP analysis to the Type A (II) scheme in (\ref{Afirst}).
 In particular, we  %consider the scheme \eqref{eq:scheme_1bis_bis} and 
prove that its limiting scheme is consistent with the continuous limit model (\ref{incomp2_bis}) at $\varepsilon = 0$.
We focus on the AP analysis on time discretizations, while keeping the space continuous. We assume that the data at time $t^n$ are well-prepared in the sense of (\ref{ExEpBis}), i.e. $p^n({\bf x}) := p(t^n,{\bf x})$ and $\bu^n({\bf x}) := \bu(t^n,{\bf x})$ admit the decomposition: 
\begin{equation}\label{InitWP}
p^n({\bf x}) = p^n_0 + \varepsilon^2 \,p_2^n({\bf x}), \quad \bu^n({\bf x}) = \bu_0^n({\bf x}) + \mathcal{O}(\varepsilon), 
\end{equation}
where %we assume%d that the boundary conditions such that 
{$p^n_0 = (\gamma -1) E^{n}_0$} is a constant and $\nabla \cdot \bu^n_0({\bf x}) = 0$. 
%Here $p_2$ is a hydrostatic pressure, while $p_0$ is a thermodynamic pressure with $\nabla p_0 = 0$ under the assumption that $p_0$ is independent of space.

Then we consider an expansion in powers of $\varepsilon$ of the form (\ref{ExEpBis}) and we plugin it into the semi-discrete scheme \eqref{eq:scheme_1bis_bis}. 
Equating to zero the $\mathcal{O}(\varepsilon^{-2})$ term we have 
\begin{equation}\label{eq:p0n+1}
\nabla p_0^{n+1} = 0.
\end{equation} 
Equating to zero the  $\mathcal{O}(\varepsilon^{0})$ terms, we have
\begin{subequations}\label{eq: incomp_discr}
\beq\label{eq: incomp_discr1}
\displaystyle \frac{ \rho_0^{n+1}-\rho_0^{n}}{\Delta t} + \nabla \cdot (\rho^{n}_0\bu_0^{n}) = 0,
\eeq
\begin{equation}\label{eq: incomp_discr2}
\displaystyle\frac{\rho^{n+1}_0\bu_0^{n+1}-\rho_0^{n}\bu_0^{n}}{\Delta t} + \nabla \cdot(\rho^{n}_0\bu_0^n \otimes \bu_0^n) + \nabla p_2^{n+1} = 0,
\end{equation}
\begin{equation}\label{eq: incomp_discr3}
E_0^{n+1} = p_0^{n+1}/(\gamma - 1), 
\end{equation}
\begin{equation}\label{eq: incomp_discr4}
\displaystyle  \frac{ E_0^{n+1}-E_0^{n}}{\Delta t}+  \nabla \cdot \left( \bar{H}^n_0 (\rho_0^{n+1}\bu_0^{n+1})\right)= 0,
\end{equation}
\end{subequations} 
with 
\begin{equation}\label{H0}
\bar{H}^n_0 = \frac{p_0^{n} + E_0^n}{\rho_0^{n+1}}=\frac{\gamma}{\gamma-1}\frac{p^n_0}{\rho^{n+1}_0}.
\end{equation}
By (\ref{eq: incomp_discr3}) and (\ref{H0}), (\ref{eq: incomp_discr4}) is equivalent to 
\begin{equation}\label{limitP}
\displaystyle  \frac{ p_0^{n+1}-p_0^{n}}{\Delta t}+  \gamma \, p^n_0\,\nabla \cdot \bu_0^{n+1}= 0.
\end{equation}
%where by (\ref{H0})
%\begin{equation*}%\label{HObis}
%\nabla \cdot \left(\bar{H}^n_0 \rho_0^{n+1}\bu_0^{n+1}\right)= \nabla \cdot [\left(p_0^{n} + E_0^n\right)\bu_0^{n+1}] =\frac{\gamma}{\gamma - 1} p^n_0\nabla \cdot \left(\bu^{n+1}_0\right)
%\end{equation*}
%have been used. 
Note that (\ref{limitP}) is the time discretization of the limiting equation for the pressure (\ref{eq: energy0}). 

From (\ref{eq:p0n+1}), it follows that $p^{n+1}_0$ is independent of space. Now, as in the continuous case, integrating (\ref{limitP}) over a spatial domain $\Omega$ and assuming some boundary conditions (for example, no-slip or periodic), we get  $p^{n+1}_0 = p^n_0$, i.e. $p^{n+1}_0$ is also independent of time. Using this in the equation (\ref{limitP}), we obtain the divergence-free condition for the velocity $\nabla \cdot \bu_0^{n+1} = 0$. 

{This yields the following theorem: 
\begin{thm}
 The time-discrete scheme (\ref{eq: incomp_discr}) is asymptotic preserving in the sense that at the leading order asymptotic expansion the equations (\ref{eq: incomp_discr1}), (\ref{eq: incomp_discr2}) with  $p^{n+1}_0$ constant and $\nabla \cdot \bu_0^{n+1} = 0$, are a consistent approximation of the incompressible Euler equations (\ref{incomp2_bis}).
\end{thm}
}
{Another noticeable feature of the time-discrete scheme (\ref{eq: incomp_discr}) is that we can obtain a time-discrete version for the elliptic equation (\ref{equ_p_2}). We get it by applying the divergence operator to the  momentum equation (\ref{eq: incomp_discr2}), after some algebraic manipulations, and making use of the density equation (\ref{eq: incomp_discr1}). The resulting relation is
\begin{equation*}%\label{eqElEq}
%\displaystyle  \frac{\gamma}{\gamma -1} p_n^0 \Delta t \left(
\nabla \cdot \left( \frac{\rho^{n}_0}{\rho_0^{n+1}}(\bu_0^n \cdot \nabla) \bu_0^n\right) = -\nabla \cdot \left( \frac{1}{\rho_0^{n+1} }\nabla p_2^{n+1}\right),
%\right)= 0.
\end{equation*}
which is a consistent discretization of Eq. (\ref{equ_p_2}) because $\rho_0^n/\rho_0^{n+1}  = 1 + \mathcal{O}(\Delta t)$.
}
%Then system (\ref{eq: incomp_discr}) is a consistent discretization of the continuous one (\ref{eq: incomp}), i.e. the scheme (\ref{eq:scheme_1bis_bis}) is AP. % asymptotic preserving.

Note that once it is established that $\nabla\cdot \bu_0^{n+1} = 0$, from the first equation of (\ref{eq: incomp}) it follows that if the initial density $\rho_0^0$ is constant, then it remains constant for any later time.

%The system for the unknowns $\bu_0^{n+1}$ and $p_2$ then becomes
%\begin{align}\label{eq:monolitic}
%\frac{\bu_0^{n+1}-\bu_0^n}{\Delta t} + \nabla\cdot (\rho_0^n \bu_0^n \otimes \bu^n_0 ) + \frac{1}{\rho_0}\nabla p_2^{n+1} & =  0\\
%\nabla\cdot \bu_0^{n+1} & =  0
%\end{align}
%which, by suitable space discretization, become the classical first order monolithic scheme for the incompressible Euler equations. 

\begin{rem}
{Note that in \cite{cordier2012asymptotic} the authors proposed a scheme similar to (\ref{eq:scheme_1bis_bis}), but with a slight difference in the explicit treatment of the density in the energy equation, i.e.
\begin{equation*}%\label{rho_n}
\frac{E^{n+1}-E^n}{\Dt} + \nabla\cdot\left(\frac{E^n + p^n}{\rho^{n}}\bq^{n+1} \right)  =0.
\end{equation*}
Their scheme is AP, i.e.  it is consistent with the limiting equations (\ref{eq: incomp}) in the incompressible regime. However, as $\varepsilon \to 0$, they prove that the divergence-free condition on the velocity is explicitly satisfied up to the order of the approximation with the CFL condition independent of the Mach-number $\varepsilon$, i.e. $\nabla \cdot \bu_0^{n+1} =\mathcal{O}(\Delta t$) (see the Proposition in Section 4.2 in \cite{cordier2012asymptotic}). In our case, instead, scheme (\ref{eq:scheme_1bis_bis}) has the correct discrete divergence-free condition for the leading order velocity $\bu^{n+1}_0$, i.e. $\nabla \cdot \bu_0^{n+1} = 0$. }
\end{rem}

%%%%%%%%%%%%%%%%%%%%%%%%%%%%%%%
\subsection{AA property}
The AP property guarantees only the consistency of the scheme, but in general the AP property does not guarantee the high order accuracy of IMEX schemes in the limit for $\varepsilon \ll 1$, i.e.\ as $\varepsilon \to 0$ the order of accuracy may degrade. In what follows, we first formally state the definition of the AA property, and then we recognize that the {SA} condition is crucial to guarantee the AA property
of our SI-IMEX R-K scheme. Similarly as the AP property, here we focus on the AA analysis on time discretizations, while keeping space continuous.

\begin{defn}
A semi-implicit temporal discretization for the Euler system (\ref{Start1}) is said to be asymptotically accurate (AA), if  it maintains its order of temporal accuracy for the limiting system (\ref{eq: incomp}) when $\varepsilon \to 0$.
\end{defn}

\begin{prop}
\label{prop: AA}
Consider an SI-IMEX  R-K scheme (\ref{IMEXscheme1})-(\ref{LI}) of order $r$ applied to system (\ref{Start1}) in a bounded domain $\Omega \subset \mathbb{R}$ with zero Neumann condition. Assume that the implicit part of the IMEX R-K scheme is SA and that the initial conditions $(\rho^0({\bf x}), \rho^0({\bf x})\bu^0({\bf x}), p^0)^T$ are well-prepared in the form of
(\ref{ExEpBis0}). 
%with {$p^0 $} being constant and $\nabla \cdot {\bu}_0({\bf x}) = 0$. 
Let us denote by $(\rho^1({\bf x};\varepsilon)$, $\rho^1({\bf x};\varepsilon)\bu^1({\bf x};\varepsilon)$, $p^1({\bf x};\varepsilon))^T$ the numerical solution after one time step. Then we have:
\begin{equation}\label{Prop1}
\lim_{\varepsilon \to 0} p^1({\bf x} ;\varepsilon)=p_*, \quad \lim_{\varepsilon \to 0} \nabla \cdot \bu^1({\bf x};\varepsilon)=0,
\end{equation}
with $p_*$ a constant.

Furthermore, let ${\bf V}_{inc}({\bf x},t) = (\rho_{inc}({\bf x},t), {\rho_{inc}({\bf x},t)\bu_{inc}({\bf x},t)}, p_{inc}({\bf x},t))^T$ be the exact solution of the incompressible Euler equations (\ref{incomp2_bis}) with {the same initial data.} 
Then one has the following one-step error estimate
\begin{equation}\label{Prop1_2}
\lim_{\varepsilon\to 0} {\bf V}^1({\bf x};\varepsilon)={\bf V}_{inc}({\bf x}, \Delta t)+\mathcal{O}(\Delta t^{r + 1}),
\end{equation}
i.e., the scheme is AA.
\end{prop}

{\bf Proof.} We consider the first step from $t^0=0$ to $t^1=\Delta t$ for the SI-IMEX R-K scheme  (\ref{IMEXscheme1})-(\ref{LI}) of order $r$ applied to system (\ref{Start1}) 
 with well-prepared initial data (\ref{InitWP}):
\begin{equation*}%\label{InitWP2}
\rho^0({\bf x}) = \rho_{inc}^0 + \varepsilon^2\, \rho_2^0({\bf x}), \quad
p^0 ({\bf x})= p_* + \varepsilon^2\, p_2^0({\bf x}), \quad \bu^0({\bf x}) = \bu_{inc}^0 + \mathcal{O}(\varepsilon), 
\end{equation*}
where $\rho^0_{inc} := \rho_{inc}({\bf x},0)$, $\bu^0_{inc}:= \bu_{inc}({\bf x},0)$ and by well prepared assumption we have: $p_{inc}({\bf x},0) := p_*$ constant independent of time and space, and $ \nabla \cdot \bu^0_{inc} =0$. 

Now we consider a formal $\varepsilon$-expansion of the quantities $U_I^{(i)} = ( \rho^{(i)}_I,   \bq^{(i)}_I , E^{(i)}_I)^T$, and $U_E^{(i)} =  ( \rho^{(i)}_E,  \bq^{(i)}_E , E^{(i)}_E)^T$ with  $\bq^{(i)}_I = \rho^{(i)}_I\bu^{(i)}_I$ and $\bq^{(i)}_E=\rho^{(i)}_E\bu^{(i)}_E$, as an example for the density and pressure:
\begin{equation}
\label{QuantExp}
\begin{array}{lll}
&\rho^{(i)}_I = \rho^{(i)}_{0,I} + \varepsilon^2\, \rho^{(i)}_{2,I} + ..., \quad   
&\rho^{(i)}_E = \rho^{(i)}_{0,E} + \varepsilon^2\, \rho^{(i)}_{2,E} + ..., \\ [3mm]
&p^{(i)}_I = p^{(i)}_{0,I} + \varepsilon^2\, p^{(i)}_{2,I} + ..., \quad   &p^{(i)}_E = p^{(i)}_{0,E} + \varepsilon^2\, p^{(i)}_{2,E} + ... .
\end{array}
\end{equation}
%and
%\begin{subequations}\label{QuantExp2}
%\beq\label{eq:1ex2}
%\rho^{(i)}_I = \rho^{(i)}_{0,I} + \varepsilon^2\, \rho^{(i)}_{2,I} + ..., 
%\eeq
%\beq\label{eq:2ex2}
%\bq^{(i)}_I  = \bq^{(i)}_{0,I} + \varepsilon^2 \, \bq^{(i)}_{2,I}   + ...,\\
%\eeq
%\beq\label{eq:3ex}
%E^{(i)}_I = E^{(i)}_{0,I} + \varepsilon^2\, E^{(i)}_{2,I} + ..., \\
%\eeq
%\beq\label{eq:4ex2}
%p^{(i)}_I = p^{(i)}_{0,I} + \varepsilon^2\, p^{(i)}_{2,I} + ..., \\
%\eeq
%\end{subequations}
In order to prove the theorem, we use the mathematical induction. 
\begin{itemize}
\item 
\emph{Asymptotic accuracy for the internal stages $i = 1, \cdots, s$}. 

Case $i = 1$ leads to the same AP analysis for the scheme 
(\ref{eq:scheme_1bis_bis}) with $\Delta t$ replaced by $a_{11} \Delta t$. 
To prove the result for $i> 1$ onwards, we make use of the induction hypothesis, assuming the property holds 
for $j\le i-1$ and prove that it holds for $j=i$.
%e 
%step, if the property holds for $i-1$, then we prove 
%that it holds for the next natural $i$ with $i = 2, \cdots, s$. 
Then for $j = 1, \cdots ,i-1$ we have 
\begin{equation}\label{cond0}
\ds p_{0,E}^{(j)} = p_*, \quad  E_{0,E}^{(j)}  = \frac{p_*}{\gamma -1},\quad  \nabla \cdot \bu_{0,I}^{(j)}=0. %\quad  \SB{\nabla \cdot \bu_{0,E}^{(j)}=0}.
\end{equation}
Now we insert the expansions (\ref{QuantExp}) into the explicit step in (\ref{IMEXscheme1}), and we get for the energy equation: 
\begin{equation}\label{Eet}
E^{(i)}_{0,E} = E_{inc}^0 -\Delta t\sum_{j = 1}^{i-1} \tilde{a}_{ij}\nabla \cdot \left( \bar{H}^{(j)}_{0} \bq_{0,I}^{(j)} \right),% \Delta t \,\gamma\, p_* \sum_{j = 1}^{i-1} \tilde{a}_{ij} \nabla \cdot \bu^{j}_{0,I} = p_*,
\end{equation}
with $E_{inc}^0 = p_{*}/(\gamma -1)$ and for $j = 1, \cdots ,i-1$,
\begin{equation}\label{H0_int}
\ds \bar{H}_0^{(j)} =  \frac{E^{(j)}_{0,E}+ p^{(j)}_{0,E}}{\rho^{(j)}_{0,I}}=\frac{\gamma}{\gamma-1}\frac{p_{*}}{\rho^{(j)}_{0,I}}.
\end{equation}
Now by (\ref{cond0}), (\ref{H0_int}) and $E^{(i)}_{0,E} = p^{(i)}_{0,E}/(\gamma -1)$, from (\ref{Eet}) we obtain 
\begin{equation}\label{P0E}
p^{(i)}_{0,E} = p_* -\Delta t \,\gamma\, p_* \sum_{j = 1}^{i-1} \tilde{a}_{ij} \nabla \cdot \bu^{(j)}_{0,I} = p_*,
\end{equation}
% \begin{equation}\label{Hh_0}
% \nabla \cdot \left( \hat{H}^{(j)}_{0} \bq_{0,I}^{(j)} \right)=  \frac{\gamma\, p_*}{(\gamma - 1)} \nabla \cdot \bu_{0,I}^{(j)}
%  \end{equation} 
%  for $j = 1, \cdots, i-1$, we get for the energy equation:
%and by (\ref{cond0})
that is, $p^{(i)}_{0,E} = p_*$ is constant. Then $E^{(i)}_{0,E} = p_*/(\gamma - 1)$ is also constant for the stage $i$. 

From (\ref{IMEXscheme1}), to $\mathcal{O}(1)$ we obtain for the density and momentum equations  
\begin{equation}\label{rho_E}
\rho_{0,E}^{(i)} =  \rho^0_{inc} 
\displaystyle  - \Delta t \sum_{j = 1}^{i-1} \tilde{a}_{ij} \nabla \cdot \bq^{(j)}_{0,E}, %\rho^{(j)}_{0,E}
\end{equation}
and
\begin{equation}\label{rhoq_E}
\bq_{0,E}^{(i)}= \bq_{inc}^0 - \Delta t \sum_{j = 1}^{i-1} \tilde{a}_{ij} \left(\nabla \cdot \left({\rho^{(j)}_{0,E}}\, \bu^{(j)}_{0,E} \otimes \bu^{(j)}_{0,E}\right) + \nabla p_{2,I}^{(j)}\right),
\end{equation}
with $\bq_{inc}^0 =(\rho \bu)_{inc}^0$ and $\nabla p^{(j)}_{0,E}=0$ for $j=1, \cdots ,i-1$.

%From (\ref{LI}), we have for the $i$-th stage
%\begin{equation}\label{iStage}
%     \left( \begin{array}{rl} 
%     & \rho^{(i)}_I\\[5mm]
%   & \bq^{(i)}_I \\[5mm]
%   & E^{(i)}_I 
%   \end{array}\right)  = \left( \begin{array}{rl} 
%     & \rho^{(i)}_*\\[5mm]
%   & \bq^{(i)}_* \\[5mm]
%   & E^{(i)}_* 
%   \end{array}\right) 
%    - \Delta t\, a_{ii}\nabla \cdot \left( 
%   \begin{array}{rl}
%  & \displaystyle  \bq^{i}_E \\[5mm]
%  &\displaystyle  \left(\frac{\bq^{(i)}_{E} \otimes \bq^{(i)}_E}{\rho^{(i)}_E} + p^{(i)}_EI\right) + (1-\varepsilon^2) \nabla p^{(i)}_{I,2}	\\[5mm]%p^{(2)}_I\\[5mm]
%&\displaystyle  \hat{H}^{(i)} \bq^{(i)}_I 
%\end{array}
%\right).
%\end{equation}
%and inserting again (\ref{QuantExp}) into it, for $\mathcal{O}(\varepsilon^{-2})$ term, it follows: $\nabla p_{0,I}^{(i)}= 0$, i.e., $p_{0,I}^{(i)}$ is independent of space.  

Similarly, inserting expansions (\ref{QuantExp}) into (\ref{QI}), up to $\mathcal{O}(1)$ we get for the intermediate explicit step $\tilde{U}_{0}^{(i)}$ in (\ref{IMEXscheme1}): 
\begin{equation}\label{q_star}
    \tilde{\bq}_{0}^{(i)} =  \bq^0_{inc} 
  \displaystyle  - \Delta t \sum_{j = 1}^{i-1} a_{ij} \left(\nabla \cdot \left({\rho^{(j)}_{0,E}} \bu^{(j)}_{0,E} \otimes \bu^{(j)}_{0,E}\right) + \nabla p_{2,I}^{(j)}\right),
  \end{equation}
  where from (\ref{eq: pI2_a}) and (\ref{cond0}), it follows $\nabla \bar{p}^{(j)}_{E} = 0$ for $j = 1,\cdots, i-1$. Furthermore, from (\ref{rhoI}) we have
  \begin{equation}\label{rho_star}
    \tilde{\rho}_{0}^{(i)} =  \rho^0_{inc} 
  \displaystyle  - \Delta t \sum_{j = 1}^{i-1} a_{ij}  \nabla \cdot \bq^{(j)}_{0,E},
    \end{equation}  
     and using (\ref{H0_int}),  we get 
    \begin{equation}\label{E_star}
    \tilde{E}_{0}^{(i)} =  E^0_{inc} 
  \displaystyle  - \Delta t \sum_{j = 1}^{i-1} a_{ij}  \nabla \cdot ( \bar{H}^{(j)}_{0} \bq_{0,I}^{(j)}) = \frac{p_*}{\gamma -1} - \frac{\gamma\, p_*}{\gamma - 1}\Delta t \sum_{j = 1}^{i-1} a_{ij}  \nabla \cdot  \bu_{0,I}^{(j)}.
    \end{equation} 
    Thus, by (\ref{cond0}) we get from (\ref{E_star})%Similarly as observed for the equation (\ref{Eet}) and by (\ref{cond0}) we get:
\begin{equation} \label{Ei0}
\tilde{E}_{0}^{(i)} = \frac{p_*}{\gamma -1}.
\end{equation}
Now from (\ref{LIbis}) and (\ref{Ei0}), it follows for the energy equation 
\begin{equation*}%\label{E_star1}
  E_{0,I}^{(i)} =  \tilde{E}^{(i)}_{0} 
  \displaystyle  - \Delta t\,  a_{ii}  \nabla \cdot ( \bar{H}^{(i)}_{0} \bq_{0,I}^{(i)}) =
  \frac{p_*}{\gamma -1}
  \displaystyle  - \Delta t \, a_{ii}  \frac{\gamma\, p_*}{\gamma - 1} \nabla \cdot  \bu_{0,I}^{(i)}.
\end{equation*} 
Considering the EOS \eqref{E_I2} to zeroth order in $\varepsilon$, we get $E_{0,I}^{(i)}=p_{0,I}^{(i)}/(\gamma-1)$, and we obtain for the pressure
  \begin{equation*}%\label{E_star2}
    p_{0,I}^{(i)}  =   p_*  + \Delta t \, \gamma \, p_* \,a_{ii}  \nabla \cdot \bu_{0,I}^{(i)}.
    \end{equation*}
 Integrating it over spatial bounded domain $\Omega$, and assuming some boundary conditions (for example, no-slip or periodic,) we first obtain $p_{0,I}^{(i)}  = p_*$ and by this we get $\nabla \cdot \bu_{0,I}^{(i)} = 0$ at the stage $i$. %Note that $p_{0,I}^{(i)}  = p_*$ is also consistent with the leading order when we take $p^{(i)}_I=\bar{p}^{(i)}_E+\eps^2\,p^{(i)}_{I,2}$ since $\bar{p}^{(i)}_{0,E}=p_*$.
 
Finally, from \eqref{LIbis}, considering (\ref{q_star}) and (\ref{rho_star}), we get for the density and momentum,
\begin{subequations} 
\begin{equation}\label{rho_complete}
    \rho_{0,I}^{(i)} =  \rho^0_{inc} 
  \displaystyle  - \Delta t \sum_{j = 1}^{i} a_{ij} \nabla \cdot \bq^{(j)}_{0,E},
\end{equation} 
\begin{equation}\label{rhoq_0}
\bq_{0,I}^{(i)}= \bq_{inc}^0 - \Delta t \sum_{j = 1}^{i}a_{ij} \left(\nabla \cdot \left({\rho^{(j)}_{0,E}}\, \bu^{(j)}_{0,E} \otimes \bu^{(j)}_{0,E}\right) + \nabla p_{2,I}^{(j)}\right),
\end{equation}
\end{subequations}
where from (\ref{eq: pI2_a}) and  $p^{(i)}_{0,E} = p_*$ it follows in the equation of the momentum $\nabla \bar{p}^{(i)}_{E} = 0$ for $i$.
 
Then equations (\ref{rho_E}), (\ref{rhoq_E}), \ref{rho_complete}), (\ref{rhoq_0}), with constant limiting pressure, i.e. $p_{0,E}^{(i)} = p_{0,I}^{(i)} = p_*$ and, the divergence free leading order velocity, i.e., $\nabla \cdot \bu_{0,I}^{(i)} = 0$, provide the discretization of system (\ref{incomp2_bis}) for the internal stage $i$ of SI-IMEX R-K scheme. This shows that in the limit $\varepsilon\to0$, the scheme becomes the same SI-IMEX RK time-discrete scheme for the incompressible Euler equations (\ref{incomp2_bis}). 

%
%Since this is valid for all $i=1,\ldots,s$, and the scheme is stiffly accurate, Then the above proof is valid for any $i$-stage, of the SI-IMEX-RK scheme for $i = 1, ...., s$.

\item \emph{Asymtotic accuracy for the numerical solution}. 

Assuming that SI-IMEX R-K scheme \eqref{DBT} is SA, then the numerical solution coincides with the last internal stage $s$, 
%i.e. 
%\[\left( \rho_0^1, \bq_{0}^{1},  E_0^1\right)^T =  \left( \rho^{(s)}_{0,I}, \bq_{0,I}^{(s)}, E^{(s)}_{0,I} \right)^T,
%\] 
 and then by setting $i = s$, we get
\begin{equation}\label{lim2}
p^1_{0} = p^{(s)}_{0,I}  = p_*, \quad \nabla \cdot \bu_{0}^1=\nabla \cdot \bu_{0,I}^{(s)} = 0,
\end{equation}
i.e., we have (\ref{Prop1}).

Now if we denote by ${\bf V}_{inc}({\bf x},t) = (\rho_{inc}({\bf x},t), \rho_{inc}({\bf x},t)\bu_{inc}({\bf x},t), p_{inc}({\bf x},t))^T$ the exact solutions of (\ref{incomp2_bis}), with initial data ${\bf V}_{inc}({\bf x},0) = (\rho^0({\bf x}), \rho^0({\bf x})\bu^0({\bf x}), p^0)^T$, {from equations (\ref{rho_E}), (\ref{rhoq_E}), (\ref{rhoq_0}) and (\ref{rho_complete}) with (\ref{lim2}), %$ p^1_{0} =p^{(s)}_{0,I} = p_*$ constant, and $\nabla \cdot \bu_{0}^1= \nabla \cdot \bu_{0}^{(s)} = 0$, 
one gets in the limit case $\varepsilon = 0$,  a SI-IMEX R-K scheme of order $r$ for the numerical solutions of equations (\ref{incomp2_bis}), that is, the SI-IMEX R-K scheme(\ref{IMEXscheme1})-(\ref{LI}) of order $r$ is asymptotically accurate (AA), and the conclusion (\ref{Prop1_2}) is obtained.}
\end{itemize}

\section{Numerical Tests}
\label{sec_numer}
\setcounter{equation}{0}
\setcounter{figure}{0}
\setcounter{table}{0}

Through an extensive set of 1D and 2D numerical tests, we will show that our scheme is uniformly stable, effective and can capture the correct asymptotic limit. We use the 3rd order SI-IMEX(4,4,3) scheme in time \eqref{IMEX1_(4,4,3)}, 4th order compact central difference discretization for 2nd order derivatives in the elliptic equation \eqref{eq33:LIbis}, and 5th order finite difference WENO reconstruction in space for both $\nabla_W$ and $\nabla_{CW}$ in Section \ref{ssec: spatial}. Overall the scheme is 4th order in space and 3rd order in time, denoted as ``S4T3". {On the other hand, an ``S2T2" scheme refers to using a 2nd order TVB reconstruction with the parameter $M=1$ instead of a 5th order WENO reconstruction in the S4T3 scheme, and a 2nd order SI-IMEX(3,3,2) scheme in time \eqref{IMEX_(3,3,2)}.} Reference solutions are computed by a 5th order finite difference WENO scheme with 3rd order explicit R-K time discretization \cite{shu1998essentially}, denoted as ``WENO5RK3".
For simplicity, we all take $\gamma=1.4$ with an ideal EOS. The time step is 
$
\Delta t = \text{CFL}\,\Delta x/\Lambda, 
$
where 
$\Lambda = \max_x \left(|u|+\min(1/\eps,1)\,c_s\right)$ in 1D, and $ \Lambda = \max_{x,y} \left(|u|+|v|+\min(1/\eps,1)\,c_s\right) $ in 2D
respectively. $c_s=\sqrt{\gamma{p}/{\rho}}$ is the scaled sound speed. We take $\text{CFL}=0.25$ for all tests. 
%{Due to the limited space, we present below only a selected set of examples, and organize the rest as a supplementary material.}

\begin{exa} \label{Ex1D}
   (Two colliding acoustic pulses \cite{Klein1995, noelle2014weakly}.) This problem is defined on the domain $-L\le x \le L=2/\eps$ with periodic boundary condition and initial data
	\begin{equation}
	\label{1dacoustic}
	\left\{
	\begin{array}{ll}
	\rho(x,0) = \rho_0 + \frac12 \varepsilon \rho_1 \left(1-\cos\left(2\pi x/L\right)\right), \quad & \rho_0 = 0.955, \quad \rho_1 = 2.0;\\ [1mm] 
	 u(x,0) = \frac12 u_0\, \text{sign}(x) \left(1-\cos\left(2\pi x/L\right)\right), \quad & u_0=2\sqrt{\gamma};\\ [1mm]
	 p(x,0) = p_0 + \frac12 \varepsilon p_1 \left(1-\cos\left(2\pi x/L\right)\right), \quad & p_0 = 1.0, \quad p_1 = 2 \gamma.
	\end{array}
	\right.
	\end{equation}	
	We first test the order of accuracy for our scheme. Due to $u(x,0)$ in \eqref{1dacoustic} is not smooth enough, in order to observe more than 2nd order accuracy, we modify it as
	\beq
	u(x,0) = u_0\sin\left(2\pi x/L\right) \left(1-\cos\left(2\pi x/L\right)\right)/2.
	\label{uacoustic}
	\eeq 
	Others are the same as in \eqref{1dacoustic}. We take $\eps=10/11$ to avoid order reduction from the high order IMEX time discretization. We take mesh sizes with $N_k = 2^k\cdot N_{0}$, where $k = 0,1,2,3$ and $N_{0} = 40$. Reference solutions are computed with $N=2560$. The errors are computed by comparing the numerical solutions for the pressure $p_E$ to its reference solution. In Table \ref{tab12}, we show the errors and orders at time $t=0.1$. The order of convergence is between $4$ and $5$, due to the dominance of spatial errors. 
	%Because the final time is relatively short, spatial errors are dominant at the resolution we tested. The observed order of accuracy for ``S4T3" is therefore between $4$ and $5$. 
	
	\begin{table}[!h]
	\caption{Example \ref{Ex1D}. Convergence test for the two colliding acoustic pulses problem with initial condition \eqref{1dacoustic}, and the velocity is replaced by \eqref{uacoustic}. $t=0.1$. $\eps=10/11$. }
    \begin{center}
%	\begin{tabular}{|c|c|c|c|c|c|}\hline
%		$N$ & $L^1$ error  & order   &  $L^\infty$ error &    order \\ \hline
%    40 &     1.62E-02 &       --&     3.95E-02 &       -- \\ \hline
%    80 &     9.97E-04 &     4.02&     3.77E-03 &     3.39 \\ \hline
%   160 &     3.54E-05 &     4.82&     1.76E-04 &     4.42 \\ \hline
%   320 &     1.34E-06 &     4.72&     8.73E-06 &     4.34 \\ \hline
%	\end{tabular}
	\begin{tabular}{|c|c|c|c|c|c|}\hline
	$N$          &  40       & 80        & 160      & 320      \\ \hline
	$L^1$ error  &  1.62E-02 & 9.97E-04  & 3.54E-05 & 1.34E-06 \\ \hline
	order        &  --       & 4.02      & 4.82     & 4.72     \\ \hline
\end{tabular}
	\label{tab12}
    \end{center}
    \end{table}

    Then we take $\eps=1/11$ with initial condition \eqref{1dacoustic} and compute the solution up to $T=1.63$ by both ``S2T2'' and ``S4T3'' schemes. A reference solution is computed with $N=2200$. In Fig. \ref{pulse}, we compare the reference solution with the numerical ones obtained with $N=22$ grid points at time $t=1.63$. We can see that the results of the ``S4T3'' scheme on the very coarse mesh match the reference solutions better than the ``S2T2'' scheme. 
 
    \begin{figure}[!h]
	\centering
	{\includegraphics[width=0.28\textwidth]{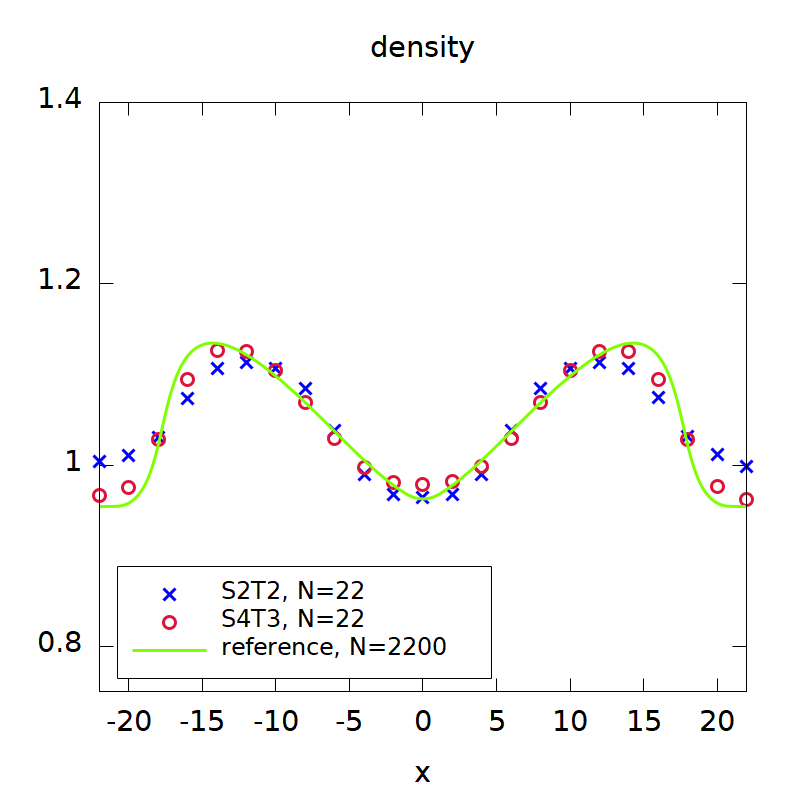}}	
	{\includegraphics[width=0.28\textwidth]{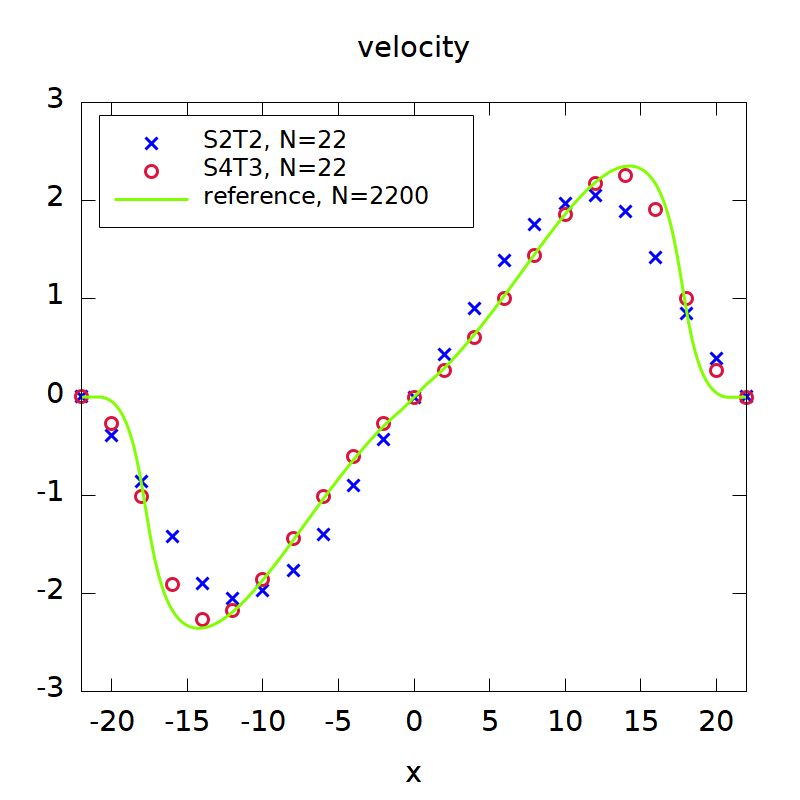}}
	{\includegraphics[width=0.28\textwidth]{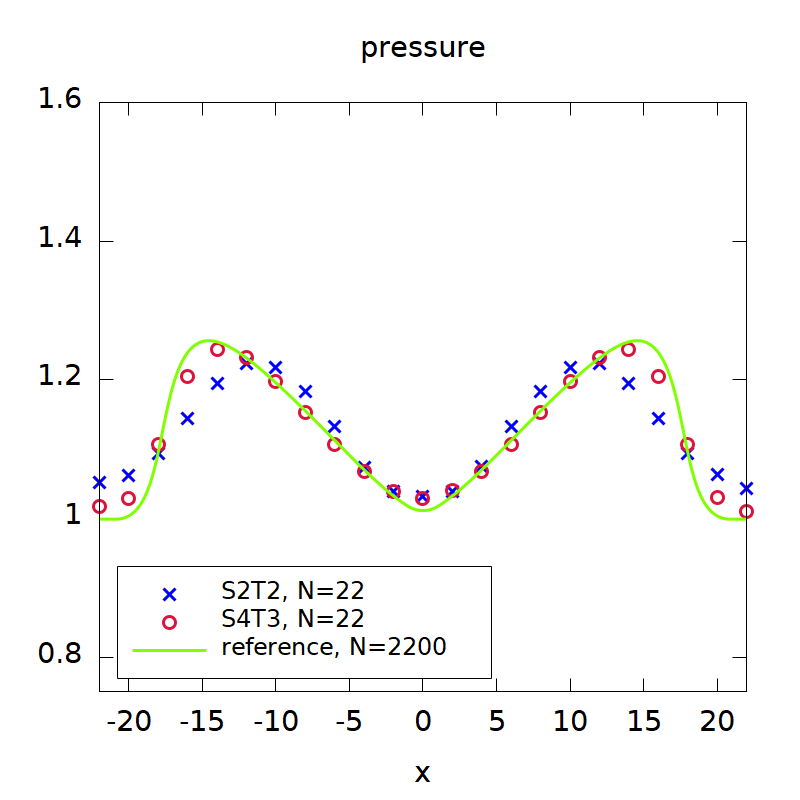}}
	\caption{Example \ref{Ex1D}. Two colliding acoustic pulses with initial condition \eqref{1dacoustic}. $t=1.63$.}
	\label{pulse}
    \end{figure} 

\end{exa}

\begin{exa}
	\label{1DRiemann}
	(1D shock tube problem.) In this example, we consider two 1D shock tube problems in the compressible regime when the Mach number is of $\mathcal{O}(1)$. We take the initial data, one is the Sod problem, where
	\begin{equation}
	\label{1Dsod}
	(\rho, u, p) = (1,0,1), \, \text{ if } x < 0.5; \quad (\rho, u, p) = (0.125,0,0.1), \, \text{ otherwise};
	\end{equation}
	another is the Lax problem, where
	\begin{equation}
	\label{1Dlax}
	(\rho, u, p) = (0.445,0.698,3.528), \, \text{ if } x < 0.5; \quad (\rho, u, p) = (0.5,0,0.571), \, \text{ otherwise},
	\end{equation}
	both with $\eps=1$ on the domain $x \in [0, 1]$. Reflective boundary conditions are considered and we take $N=50$. We compute the numerical solution by both ``S2T2'' and ``S4T3'' schemes. The results are shown in Fig. \ref{shock1d}, and compared to the exact solutions. For both the Sod and Lax shock tube problems, the results match the exact solutions very well, which show that our scheme in the moderate Mach regime ($\eps$ of order 1) can capture strong discontinuities without any observable numerical oscillations. 
	
	\begin{figure}[!h]
	\centering
	{\includegraphics[width=0.28\textwidth]{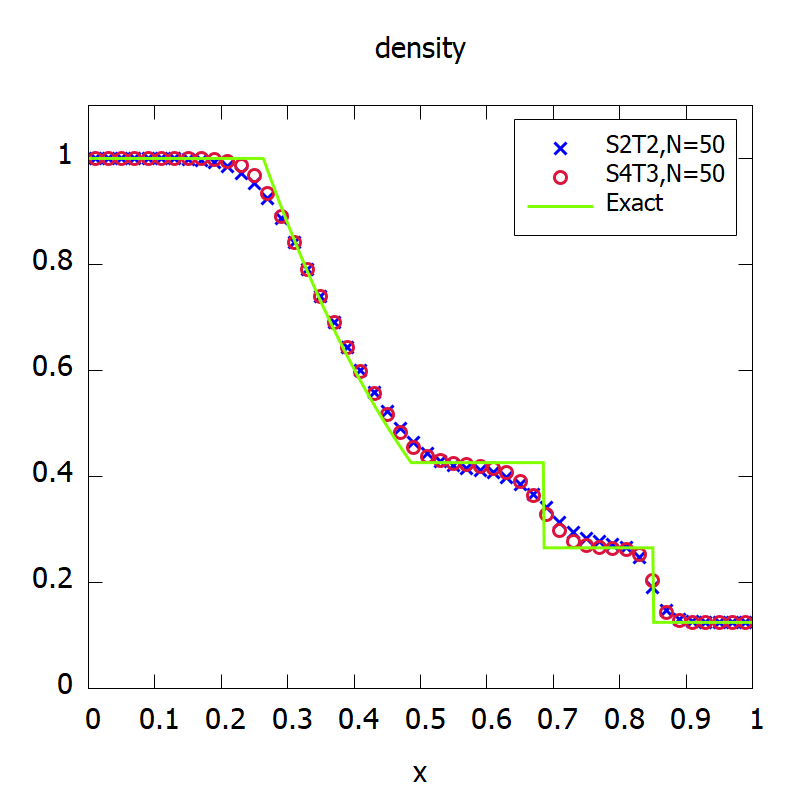}}
	{\includegraphics[width=0.28\textwidth]{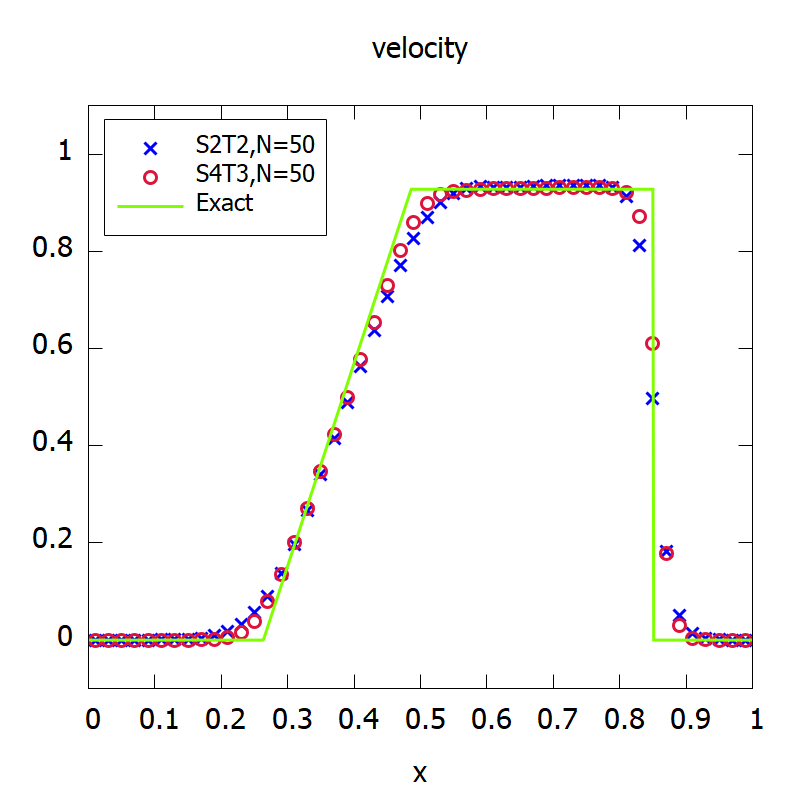}}
	{\includegraphics[width=0.28\textwidth]{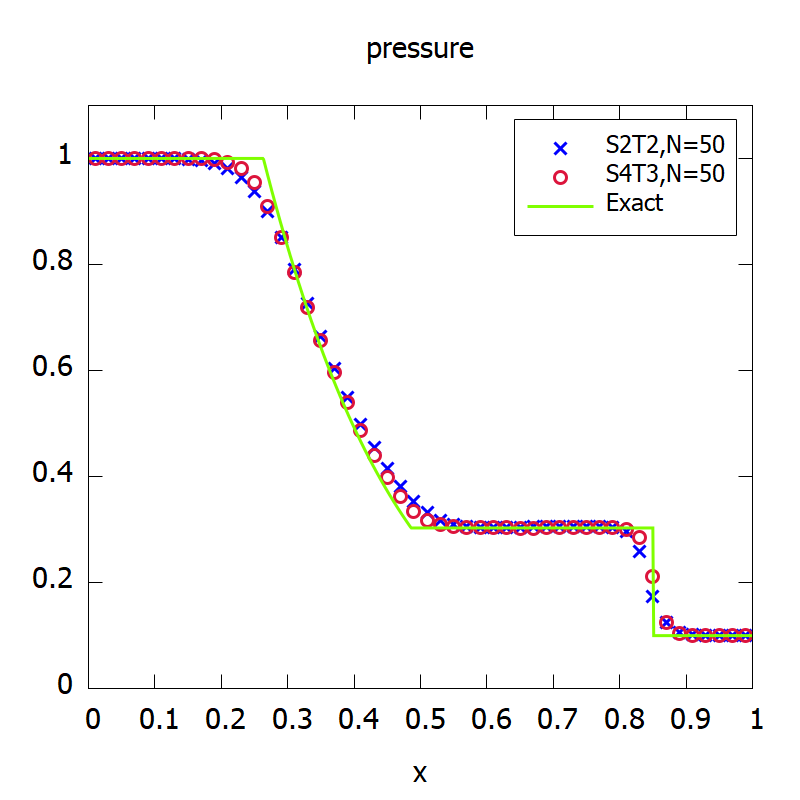}}
	{\includegraphics[width=0.28\textwidth]{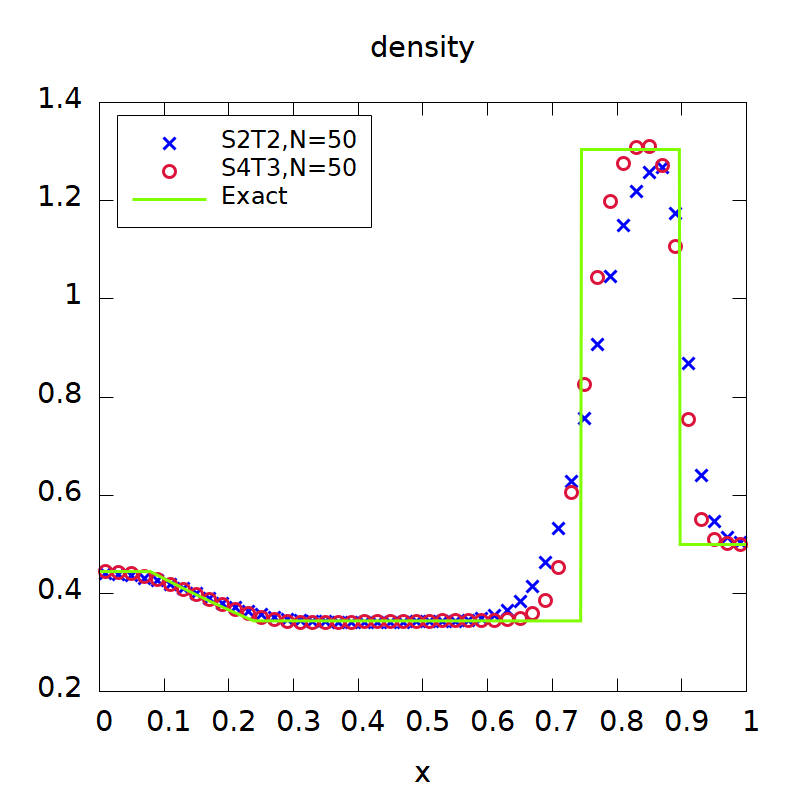}}	
	{\includegraphics[width=0.28\textwidth]{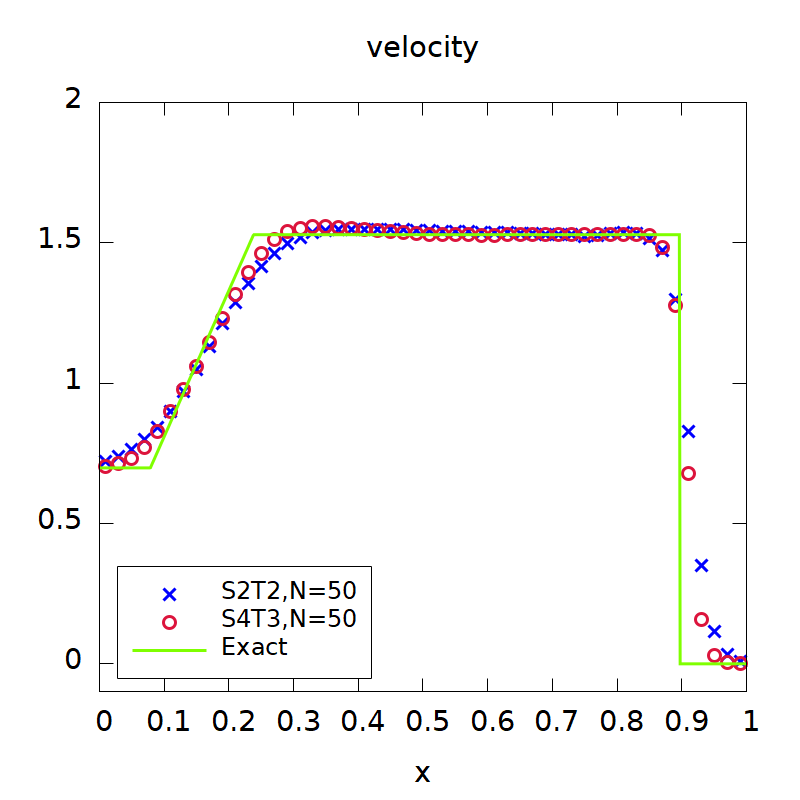}}	
	{\includegraphics[width=0.28\textwidth]{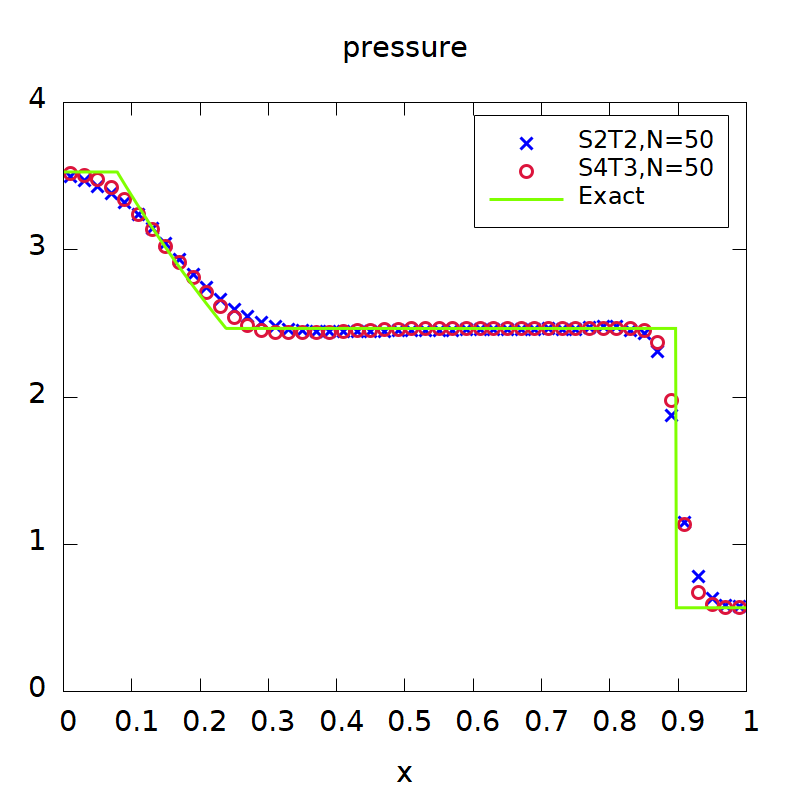}}	
	\caption{Example \ref{1DRiemann}. 1D shock tube problem for $\eps=1$. Top: the Sod problem at $t=0.2$. Bottom: the Lax problem at $t=0.16$. Mesh grid: $N=50$. The solid lines are the exact solutions.}
	\label{shock1d}
    \end{figure} 

\end{exa}

\begin{exa}
\label{Ex2D} (2D Convergence test.) We take the initial data
\begin{equation}
\left\{
\begin{array}{rcl}
\rho(x,y,0) &=& 1 + \varepsilon^2 \,\sin^2(2\pi(x+y)), \\ [1mm]
\rho(x,y,0)\,u(x,y,0) &=& \sin(2\pi(x-y)) + \varepsilon^2 \, \sin(2\pi(x+y)), \\ [1mm]
\rho(x,y,0)\,v(x,y,0) &=& \sin(2\pi(x-y)) + \varepsilon^2 \, \cos(2\pi(x+y)),
\end{array}
\right.
\label{acc2d}
\end{equation}
on the domain $\Omega=[0,1]^2$ with periodic boundary conditions. 
Initially we take $p=\rho^\gamma$. We choose $N=N_x=N_y$ and refine the mesh size by $N_k = 2^k\cdot N_{0}$, for $k = 0,1,2$, with $N_{0} = 32$. The numerical errors are computed by comparing the numerical solutions of momentum $q_2$ to the reference solution, which is computed by the ``S4T3'' scheme on the mesh $N_x=N_y=512$. In Table \ref{tab2}, we show the errors and orders at time $t=0.02$. Around 4th order for $\eps=1$ and 5th order for $\eps=10^{-6}$ are observed. For the intermediate value of $\eps=10^{-2}$, order reduction is observed. {In Fig. \ref{2dacc}, we display the orders versus $\eps$ for $\eps\in[10^{-6},1]$, where the order is computed by comparing the errors on the mesh $64\times64$ and $128\times128$. Order reduction for intermediate $\eps$'s is observed. 
%We can see that the orders gradually reduce to zero, and then get back to $5$th order for smaller $\eps$'s, which is a typical behavior of high order IMEX schemes for these multi-scale problems \cite{boscarino2013implicit}. 
}

%We   
%\begin{table}[!h]	
%\caption{Example \ref{Ex2D}. Convergence test for 2D full Euler equations with initial condition 
%		\eqref{acc2d} at $t=0.02$. Errors are computed by comparison with a reference solution. }
%	\begin{center}
%		\begin{tabular}{*{9}{|c|cc|cc|cc|cc|cc|cc|cc|}} 
%			\hline
%			\multicolumn{1}{|c|}{} & \multicolumn{2}{|c|}{ $\varepsilon=1$ } &\multicolumn{2}{|c|}{ $\varepsilon=10^{-1}$} &\multicolumn{2}{|c|}{ $\varepsilon=10^{-2}$}   &\multicolumn{2}{|c|}{ $\varepsilon=10^{-3}$} & \multicolumn{2}{|c|}{ $\varepsilon=10^{-6}$} \\ \hline
%			$N$ & $L^1$ error & order & $L^1$ error & order & $L^1$ error & order & $L^1$ error & order & $L^1$ error & order \\ \hline
%			%16  & 3.54e-02  &  --     &  2.54e-03   &  --   &  1.00e-03 &  --  \\
%            %32  & 4.64e-03  &  2.92   &  2.45e-03   &  --   &  7.26e-05 & 3.79 \\
%            32  & 4.64e-03  &  --     & 1.79E-04 & --   &  2.45e-03   &  --   & 4.76E-05 & --   & 7.26e-05 & --   \\
%            64  & 2.82e-05  &  4.05   & 1.82E-05 & 3.30 &  2.68e-03   &  --   & 1.75E-06 & 4.77 & 1.79e-06 & 5.34 \\
%            128 & 1.34e-06  &  4.39   & 2.11E-06 & 3.11 &  1.43e-03   &  0.91 & 7.64E-07 & 1.19 & 4.81e-08 & 5.22 \\
%			\hline
%		\end{tabular}
%		\label{tab2}
%	\end{center}
%\end{table}

\begin{table}[!h]	
	\caption{Example \ref{Ex2D}. Convergence test for 2D full Euler equations with initial condition 
		\eqref{acc2d} at $t=0.02$. Errors are computed by comparison with a reference solution. }
	\begin{center}
		\begin{tabular}{*{9}{|c|cc|cc|cc|cc|cc|}} 
			\hline
			\multicolumn{1}{|c|}{} & \multicolumn{2}{|c|}{ $\varepsilon=1$ }  &\multicolumn{2}{|c|}{ $\varepsilon=10^{-2}$}   & \multicolumn{2}{|c|}{ $\varepsilon=10^{-6}$} \\ \hline
			$N$ & $L^1$ error & order & $L^1$ error & order & $L^1$ error & order \\ \hline
			%16  & 3.54e-02  &  --     &  2.54e-03   &  --   &  1.00e-03 &  --  \\
			%32  & 4.64e-03  &  2.92   &  2.45e-03   &  --   &  7.26e-05 & 3.79 \\
			32  & 4.64e-03  &  --     &  2.45e-03   &  --   &  7.26e-05 & --   \\
			64  & 2.82e-05  &  4.05   &  2.68e-03   &  --   &  1.79e-06 & 5.34 \\
			128 & 1.34e-06  &  4.39   &  1.43e-03   &  0.91 &  4.81e-08 & 5.22 \\
			\hline
		\end{tabular}
		\label{tab2}
	\end{center}
\end{table}

    \begin{figure}[!h]
	\centering 
	{\includegraphics[width=0.55\textwidth]{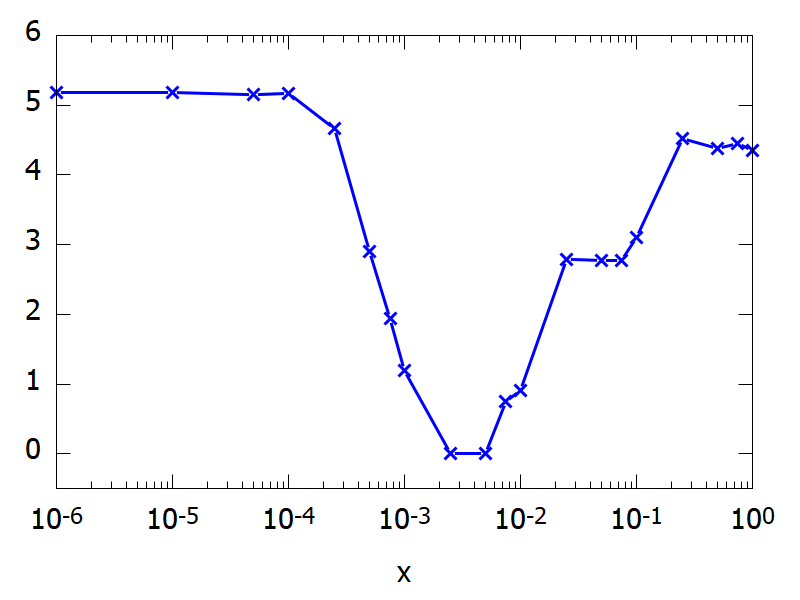}}
	\caption{{Example \ref{Ex2D}, 2D convergence test for the ``S4T3'' scheme. Convergence orders versus $\eps$'s for $\eps\in[10^{-6},1]$. As the space error dominates, the scheme is observed to have 5th order convergence for small and around 4th order for large values of $\eps$. Order degradation is observed for intermediate values of the Mach number.} }
	\label{2dacc}
\end{figure} 

\end{exa}

\begin{exa}
\label{2Dshock} (2D Riemann problem.) This is a 2D Riemann problem in the high Mach regime. We take $\eps=1$ and $\Omega=[-1,1]^2$. The initial data are defined in four quadrants, which are modified from Configuration 3 in \cite{lax1998} with four shocks
\begin{equation}
\label{shock2d1}
(\rho,u,v,p) = 
\left\{
\begin{array}{lll}
(1.5,0,0,1.5),             & x \ge 0.8,\, & y \ge 0.8; \\ [1mm]
(0.5323,1.206,0,0.3) ),    & x < 0.8,\,   & y \ge 0.8; \\ [1mm]
(0.138,1.206,1.206,0.029), & x < 0.8,\,   & y < 0.8;   \\ [1mm]
(0.5323,0,1.206,0.3),      & x \ge 0.8,\, & y < 0.8. 
\end{array}
\right. 
\end{equation}
The second one is the Configuration 5 in \cite{lax1998} with four contact discontinuities, where the initial data are
\begin{equation}
	\label{shock2d2}
	(\rho,u,v,p) = 
	\left\{
	\begin{array}{lll}
	(1,-0.75,-0.5,1),     & x \ge 0.5,\, & y \ge 0.5; \\ [1mm]
	(2,-0.75, 0.5,1),     & x < 0.5,  \, & y \ge 0.5; \\ [1mm]
	(1, 0.75, 0.5,1),     & x < 0.5,  \, & y < 0.5;   \\ [1mm]
	(3, 0.75,-0.5,1),     & x \ge 0.5,\, & y < 0.5.
	\end{array}
	\right.  
\end{equation}	
We take mesh size $N_x\times N_y=400\times 400$. For this example, {four different approaches are compared to the reference solution. The first one is our main approach, which is splitting by taking $\alpha=1$ in \eqref{Hfu}, using both spatial discretizations $\nabla_{CW}$ and $\nabla_{W}$ described in Section \ref{ssec: spatial} (``A1"); another one uses the same splitting, but replaces $\nabla_{CW}$ all by $\nabla_{W}$ (``A2"); the third one is no splitting, namely, we take $\alpha=0$ in \eqref{Hfu} and only $\nabla_{W}$ is used (``A3"); the last one is the ``S2T2'' scheme using the same approach as ``A1''.
In Fig. \ref{2dshock1}, we show the surface plots of the density at $T=0.8$ for the initial data \eqref{shock2d1}. We can see that the solution of A1 is the closest to the reference solution and no obvious oscillations are observed. For the other two approaches A2 and A3, the solutions do not perform well, numerical oscillations can be clearly seen in the middle region. The result of the 2nd order ``S2T2'' scheme is close to ``A1'', but it has poorer resolutions. We also show the cutting plots along two different lines. A1 approach is observed to perform better than the A2 and A3 approaches, and ``S2T2'' is in between but a little closer to ``A1'', which show the importance of characteristic reconstructions. In Fig. \ref{2dshock2}, we present same results for the 2D Riemann initial data \eqref{shock2d2} at $T=0.23$. Similar observations can be made.}
%{We present performance of these schemes on another Riemann initial data in the supplementary material.}
%Similar results hold for other 2D Riemann problems, here we omit them to save space.

%\begin{figure}[!h]
%	\centering
%	{\includegraphics[width=0.28\textwidth]{./pic/conf33.png}} 
%	{\includegraphics[width=0.28\textwidth]{./pic/conf3.png}} 
%	{\includegraphics[width=0.28\textwidth]{./pic/_conf3cut2.png}}\\
%	{\includegraphics[width=0.28\textwidth]{./pic/conf31.png}} 
%	{\includegraphics[width=0.28\textwidth]{./pic/conf32.png}} 
%	{\includegraphics[width=0.28\textwidth]{./pic/_conf3cut3.png}}	
%	\caption{Example \ref{2Dshock}. Surface and cut plots of the density for 2D Riemann problem \eqref{shock2d1} at $T=0.8$. Mesh grid: $400\times 400$. Top left: reference solution; top middle: the approach A1, splitting with characteristic-wise reconstruction $\nabla_{CW}$; bottom left: the approach A2, splitting with only component-wise reconstruction $\nabla_{W}$; bottom middle: the approach A3, no splitting. Top right and bottom right are 1D cuts of solutions.}
%	\label{2dshock1}
%\end{figure} 

\begin{figure}[!h]
	\centering
	{\includegraphics[width=0.28\textwidth]{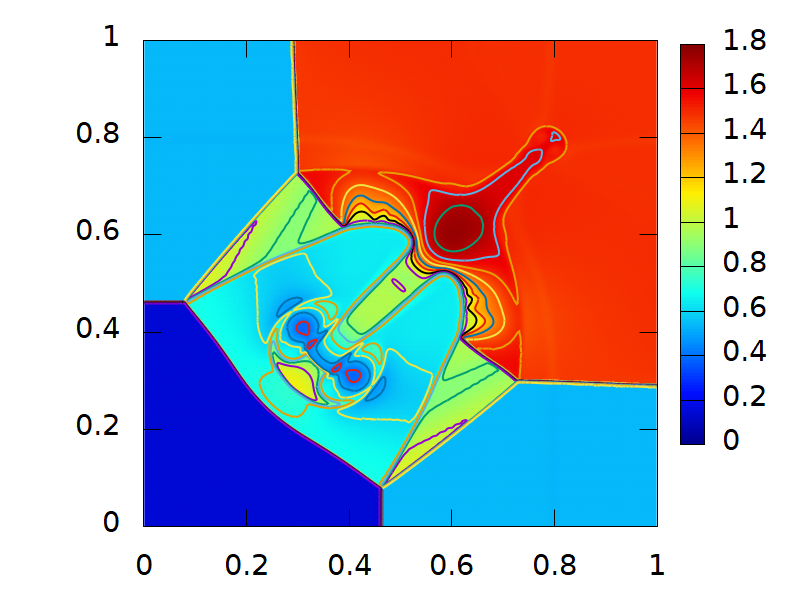}} 
	{\includegraphics[width=0.28\textwidth]{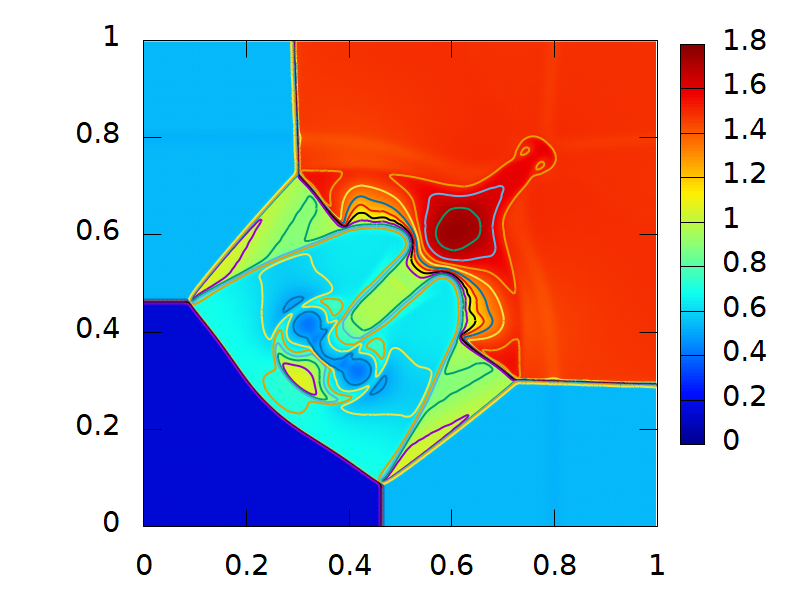}} \\
	{\includegraphics[width=0.28\textwidth]{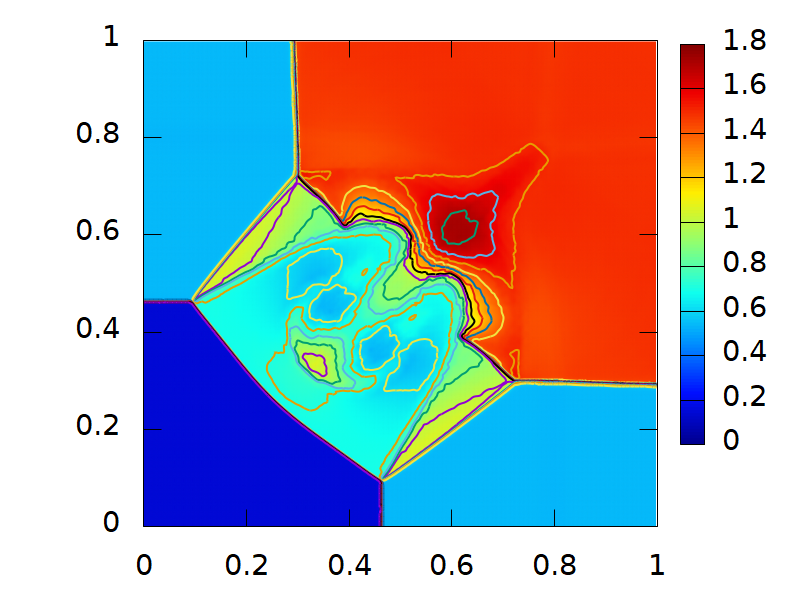}} 
	{\includegraphics[width=0.28\textwidth]{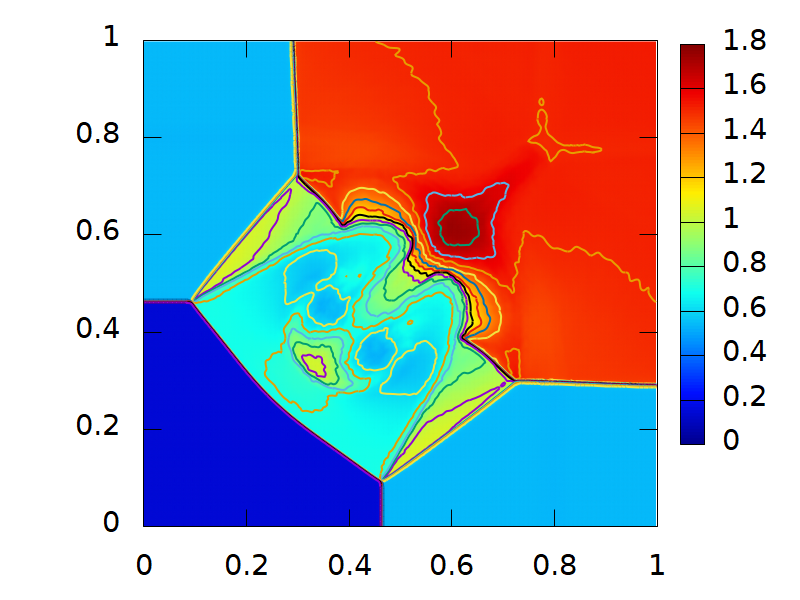}} 
	{\includegraphics[width=0.28\textwidth]{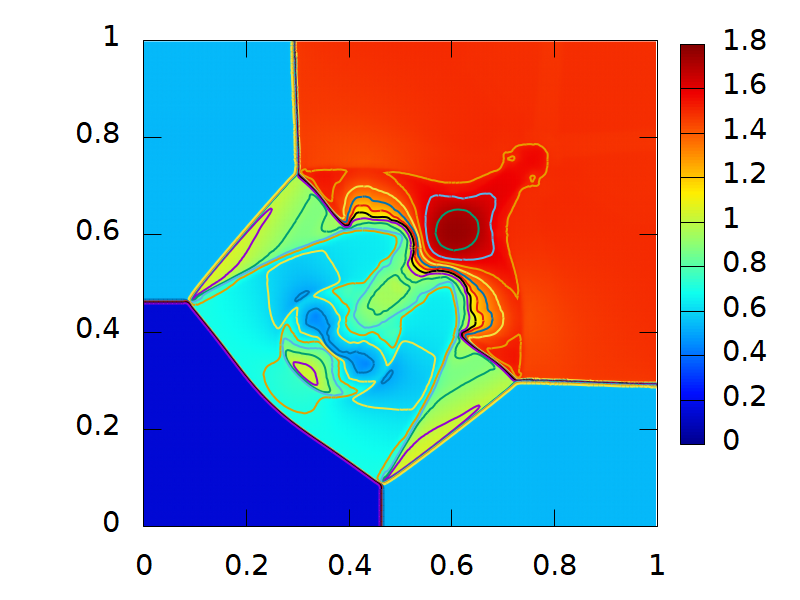}} \\
	{\includegraphics[width=0.28\textwidth]{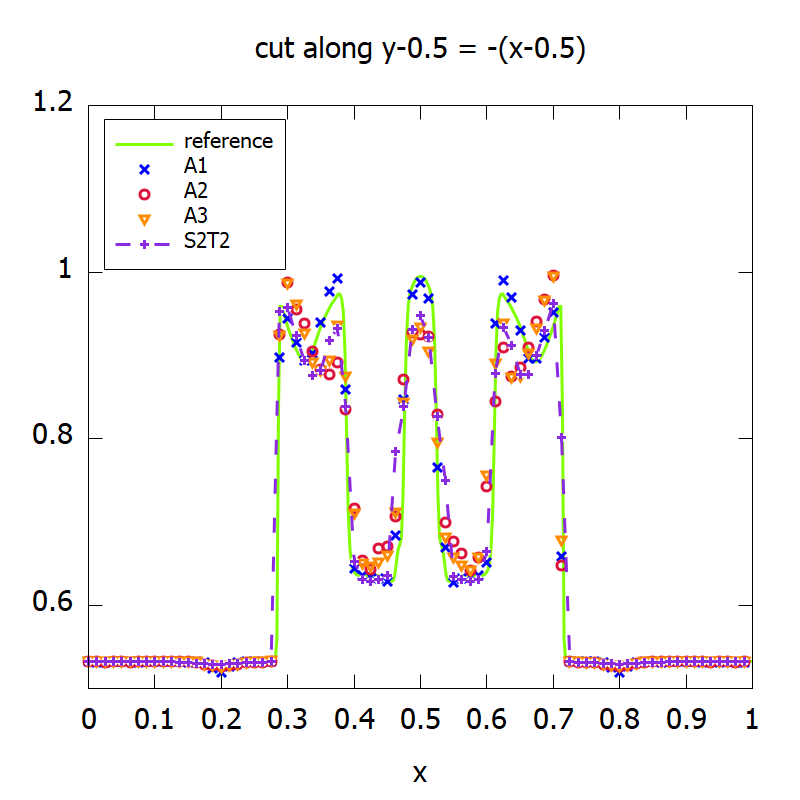}}
	{\includegraphics[width=0.28\textwidth]{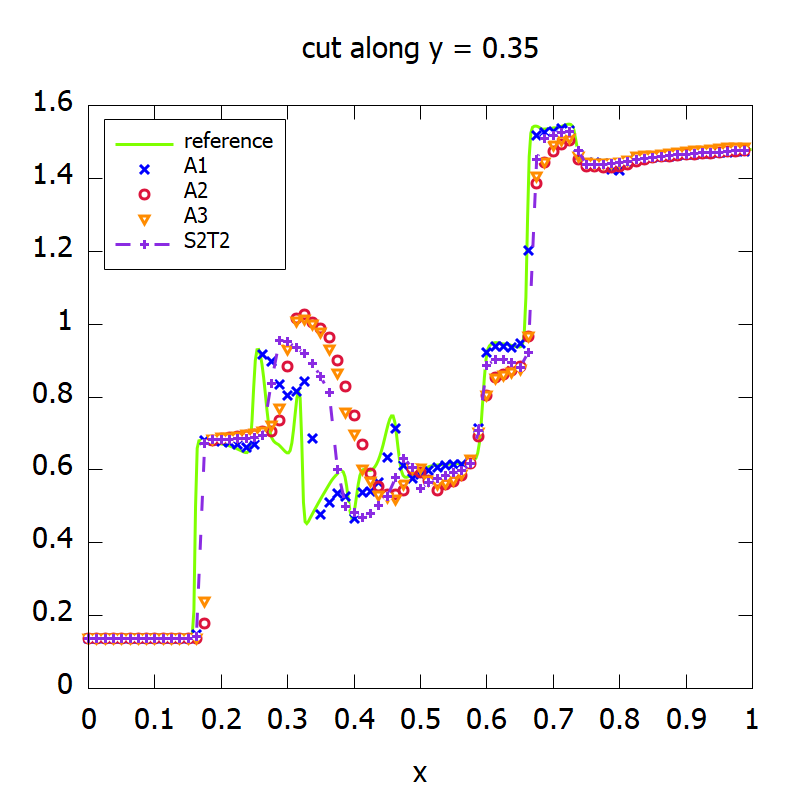}}	
	\caption{Example \ref{2Dshock}. Surface and cut plots of the density for 2D Riemann problem \eqref{shock2d1} at $T=0.8$. Mesh grid: $400\times 400$. {Top left: reference solution; top right: the approach A1, splitting with characteristic-wise reconstruction $\nabla_{CW}$; middle left: the approach A2, splitting with only component-wise reconstruction $\nabla_{W}$; middle middle: the approach A3, no splitting; middle right: S2T2. Bottom are 1D cuts of these solutions.}}
	\label{2dshock1}
\end{figure}

%\begin{figure}[!h]
%	\centering
%	{\includegraphics[width=0.28\textwidth]{./pic/conf53.png}} 
%	{\includegraphics[width=0.28\textwidth]{./pic/conf5.png}} 
%	{\includegraphics[width=0.28\textwidth]{./pic/_conf5cut2.png}}
%	{\includegraphics[width=0.28\textwidth]{./pic/conf51.png}} 
%	{\includegraphics[width=0.28\textwidth]{./pic/conf52.png}} 
%	{\includegraphics[width=0.28\textwidth]{./pic/_conf5cut3.png}}
%	\caption{Example \ref{2Dshock}. Surface plots of the density for 2D Riemann problem \eqref{shock2d2} at $T=0.23$. Mesh grid: $400\times 400$. Top left: reference solution; top middle: the approach A1, splitting with characteristic-wise reconstruction $\nabla_{CW}$; bottom left: the approach A2, splitting with only component-wise reconstruction $\nabla_{W}$; bottom middle: the approach A3, no splitting. Top right and bottom right are 1D cuts of solutions.}
%	\label{2dshock2}
%\end{figure}

\begin{figure}[!h]
	\centering
	{\includegraphics[width=0.28\textwidth]{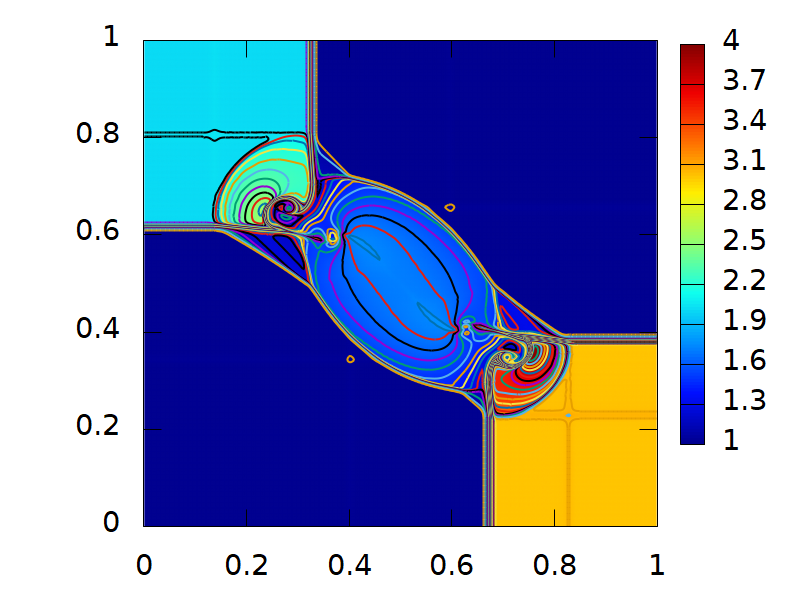}} 
	{\includegraphics[width=0.28\textwidth]{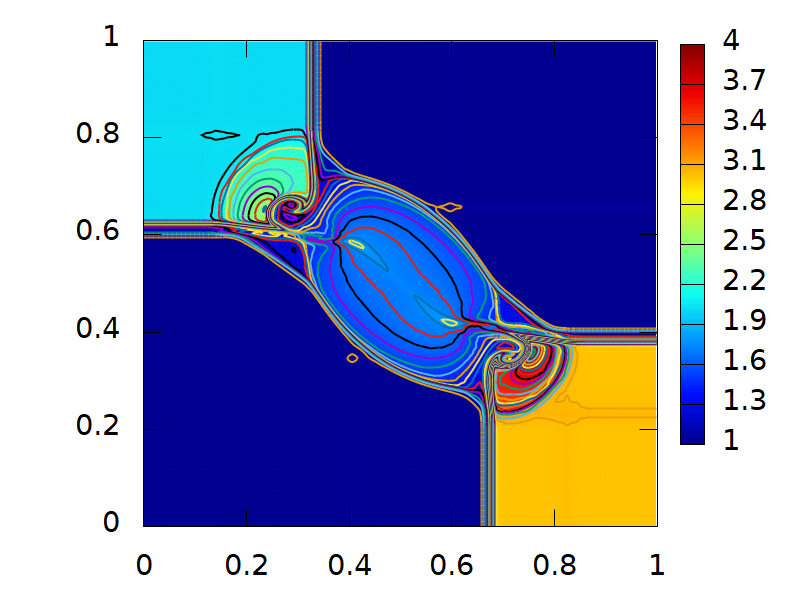}} \\
	{\includegraphics[width=0.28\textwidth]{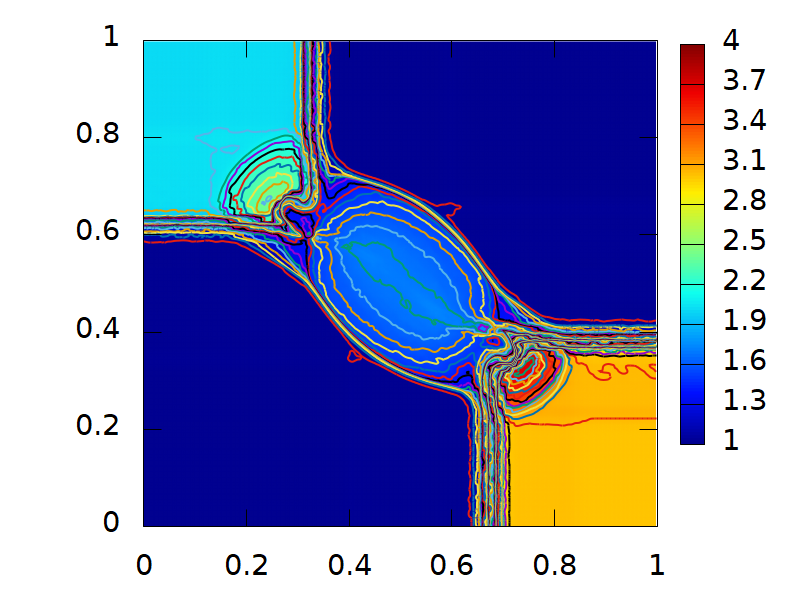}} 
	{\includegraphics[width=0.28\textwidth]{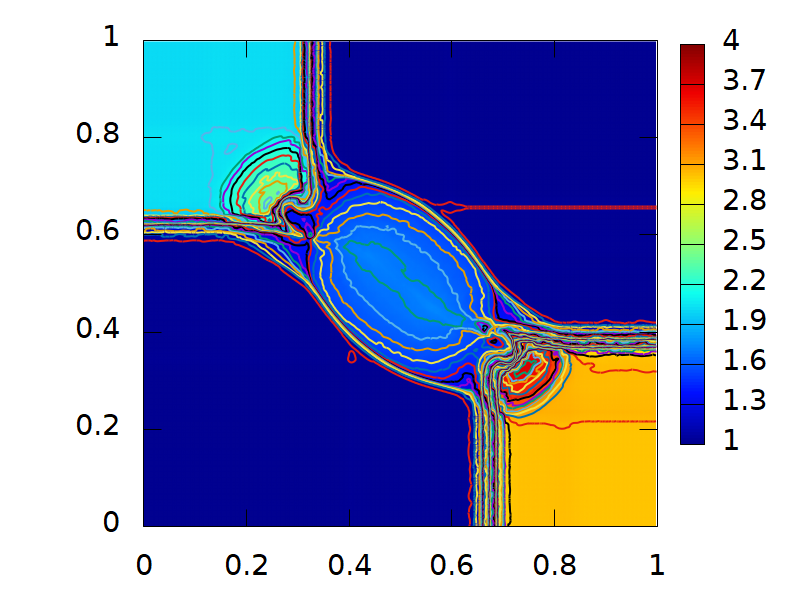}} 
	{\includegraphics[width=0.28\textwidth]{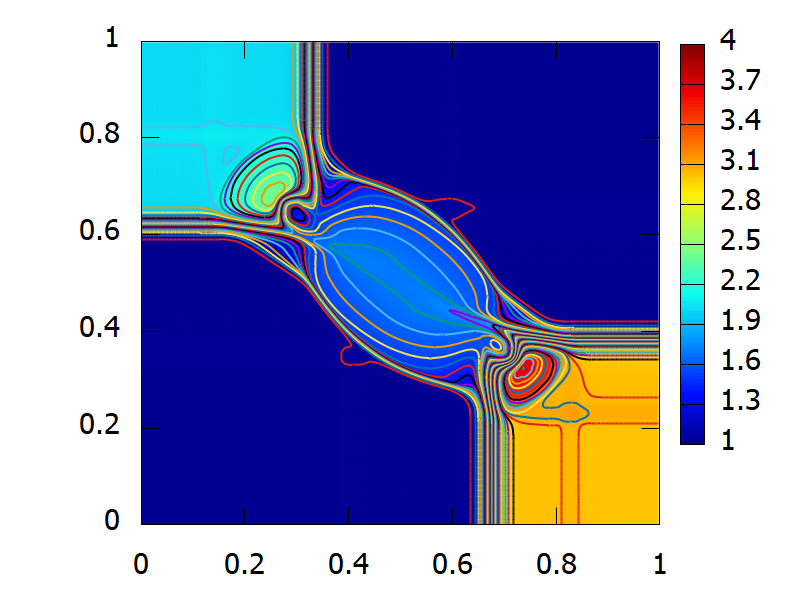}} \\
	{\includegraphics[width=0.28\textwidth]{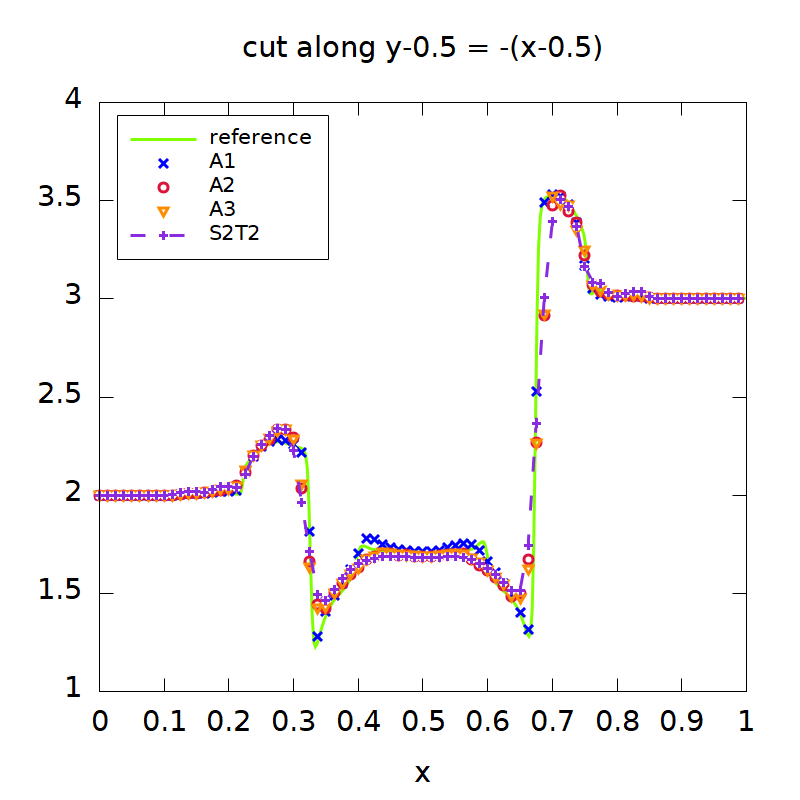}}	
	{\includegraphics[width=0.28\textwidth]{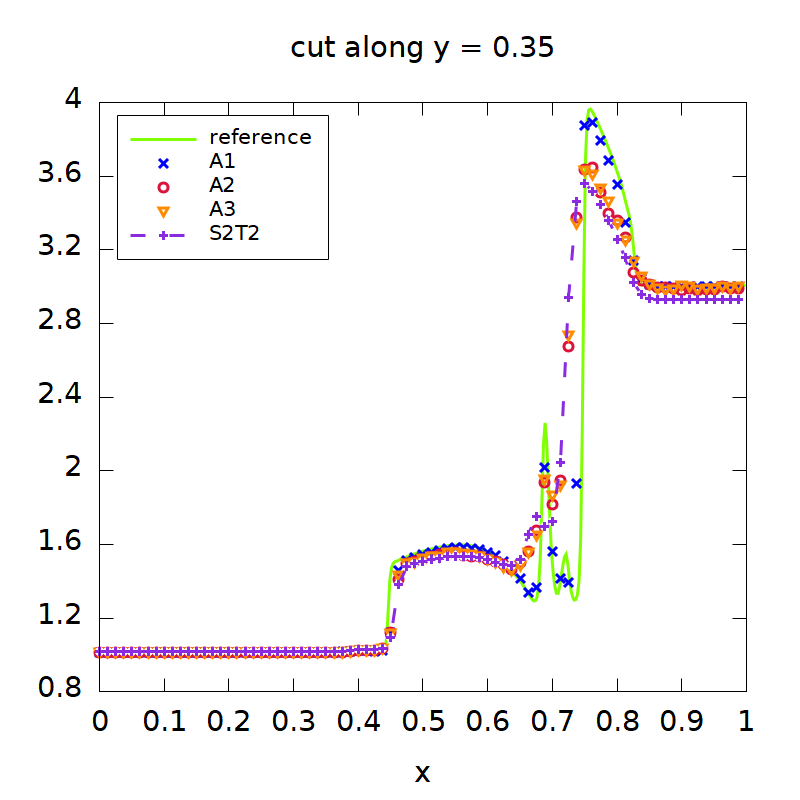}}
	\caption{Example \ref{2Dshock}. Surface plots of the density for 2D Riemann problem \eqref{shock2d2} at $T=0.23$. Mesh grid: $400\times 400$. {Top left: reference solution; top right: the approach A1, splitting with characteristic-wise reconstruction $\nabla_{CW}$; middle left: the approach A2, splitting with only component-wise reconstruction $\nabla_{W}$; middle middle: the approach A3, no splitting; middle right: S2T2. Bottom are 1D cuts of these solutions.}}
	\label{2dshock2}
\end{figure}

\end{exa}

\begin{exa}
	\label{gresho} (Gresho vortex \cite{noelle2014weakly,BoRuSc18,boscarino2019high}.) This is the time-dependent rotational Gresho vortex problem for the full Euler system. Initially a vortex is centered at $(x_0,y_0)=(0.5,0.5)$ with radius $R=0.4$ in the domain $[0,1]^2$. The initial background state is set as:
$\rho_\infty = 1, \,\, {\bf{u}}_\infty =(u_\infty,0), \,\, p_\infty=1, \,\,c_\infty=\sqrt{\gamma p_\infty/\rho_\infty} = \sqrt{\gamma}.
$
The transverse velocity for the vortex is given by
	\begin{align}
	\label{utheta}
	&u_\theta(r) = 
	\left\{
	\begin{array}{lll}
	2r/R,    & \text{ if } & 0\le r < R/2,  \\ [1mm]
	2(1-r/R),& \text{ if } & R/2 \le r < R, \\ [1mm]
	0,       & \text{ if } & r \ge R,
	\end{array}
	\right. 
	\end{align}
	and the corresponding velocity components are
	\[
	u(x,y,0) = u_\infty - \frac{y-y_0}{r}u_\theta(r), \,\, v(x,y,0) = \frac{x-x_0}{r}u_\theta(r).
	\]		
	The centrifugal force is balanced by the pressure gradient, 
%	that is
%	\[
%	\rho_\infty \frac{u_\theta^2}{r} =  \frac{1}{\eps^2}\frac{\partial p_r}{\partial r}, 
%	\]
	so the (scaled) pressure is given by
	\begin{align}
	\label{pr}
	p(r) = p_\infty + \eps^2
	\left\{
	\begin{array}{lll}
	2(r/R)^2 +2 -\log 16,         &\text{ if } & 0 \le  r < R/2, \\ [1mm]
	2(r/R)^2-4(2r/R+\log(r/R))+6, &\text{ if } & R/2\le r <R, \\ [1mm]
	0,                            &\text{ if } & r \ge R,
	\end{array}
	\right. 
	\end{align}
	where $r=\sqrt{(x-0.5)^2+(y-0.5)^2}$. 
	Periodic boundary conditions on both directions are used and the mesh is $N_x\times N_y=100\times 100$. The background velocity is $u_\infty=0.1$ moving in the $x$-direction, and we take $\eps=10^{-1}, 10^{-2}, 10^{-6}$. 
%	For this problem, the vortex transverse velocity for the rotation is $u_\theta$, corresponding to the angular velocity $\omega = u_\theta/r$. 
	The rotation period is $2\pi/\omega=R\pi$ if $r\le R/2$. We take $T=R\pi$ as one rotating period \cite{boscarino2019high}.
	
	We define the ratio between the local and the maximum Mach number as
    \beq
    \label{rMach}
     M_{\text{ratio}} =\sqrt{[(u-u_\infty)^2+v^2]/(\gamma p/\rho)}.
    \eeq			
	In Fig. \ref{2dgresho} (top left), we show the surface plot of $M_{\text{ratio}}$ for the initial condition with $\eps=10^{-2}$, which is very similar for other $\eps$'s. As time evolves, the solution will rotate while moving in the $x$-direction, but its shape will be kept. However, numerically the shape will be damped due to numerical viscosity. This is a standard example to check whether the numerical viscosity greatly depends on the parameter $\eps$ or not. {For better illustration, we show the cutting plots along $x=0.5$ at two different times: $1$ period and $2$ periods, for both schemes ``S4T3'' and ``S2T2'', and compare them to the corresponding initial shapes, for $\eps=10^{-1}, 10^{-2}, 10^{-6}$ respectively. In these plots, the vortex are shifted to the center by $-u_\infty t$ periodically. Two meshes are used: $100\times100$ and $200\times200$.
	We can see that the shapes are preserved relatively well for both "S4T3" and ``S2T2'' schemes, and the results of ``S4T3'' are a little better than those of ``S2T2''. In Fig. \ref{2dgresho}, we also show the time evolution of the averaged kinetic energy, which is defined as
	\beq
	\label{kie}
	\text{kinetic energy } =  \sum_{i,j}\left[(u(x_i,y_j,t)-u_\infty)^2+v(x_i,y_j,t)^2\right]/N_xN_y.
	\eeq
%	on both meshes $100\times100$ and $200\times200$.
	We can see the conservation of kinetic energy is well maintained. For a coarse mesh $100\times100$, the high order scheme ``S4T3'' preserves the kinetic energy clearly better than ``S2T2''. Refining the mesh can greatly improve the conservation of the kinetic energy, and the differences between these two methods are reduced especially for small $\eps$'s.}

%    \begin{figure}[!h]
%	\centering
%	%{\includegraphics[width=0.48\textwidth]{./pic/pic/2D/Gresho/gresho01.png}} 
%	%{\includegraphics[width=0.48\textwidth]{./pic/pic/2D/Gresho/gresho41.png}}
%	{\includegraphics[width=0.48\textwidth]{./pic/pic/2D/Gresho/Gresho02.png}}	
%	{\includegraphics[width=0.48\textwidth]{./pic/pic/2D/Gresho/Gresho02m.png}}
%	%{\includegraphics[width=0.48\textwidth]{./pic/pic/2D/Gresho/gresho42.png}}
%	%{\includegraphics[width=0.48\textwidth]{./pic/pic/2D/Gresho/gresho03.png}}	
%	%{\includegraphics[width=0.48\textwidth]{./pic/pic/2D/Gresho/gresho43.png}}
%	\caption{Example \ref{gresho}, the ratio of the initial local Mach number over the global Mach number $M_{\text{ratio}}$ \eqref{rMach} for the Gresho vortex problem. $\eps=10^{-2}$. Mesh grid: $100\times 100$.}
%	\label{2dgreshoini}
%    \end{figure}

    \begin{figure}[!h]
	\centering
	%{\includegraphics[width=0.48\textwidth]{./pic/pic/2D/Gresho/gresho41.png}} 
	{\includegraphics[width=0.28\textwidth]{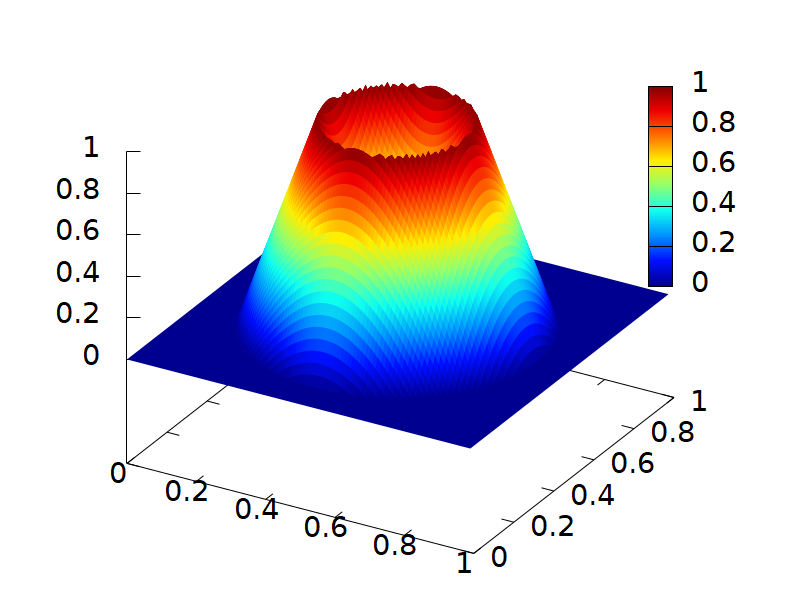}}
	{\includegraphics[width=0.28\textwidth]{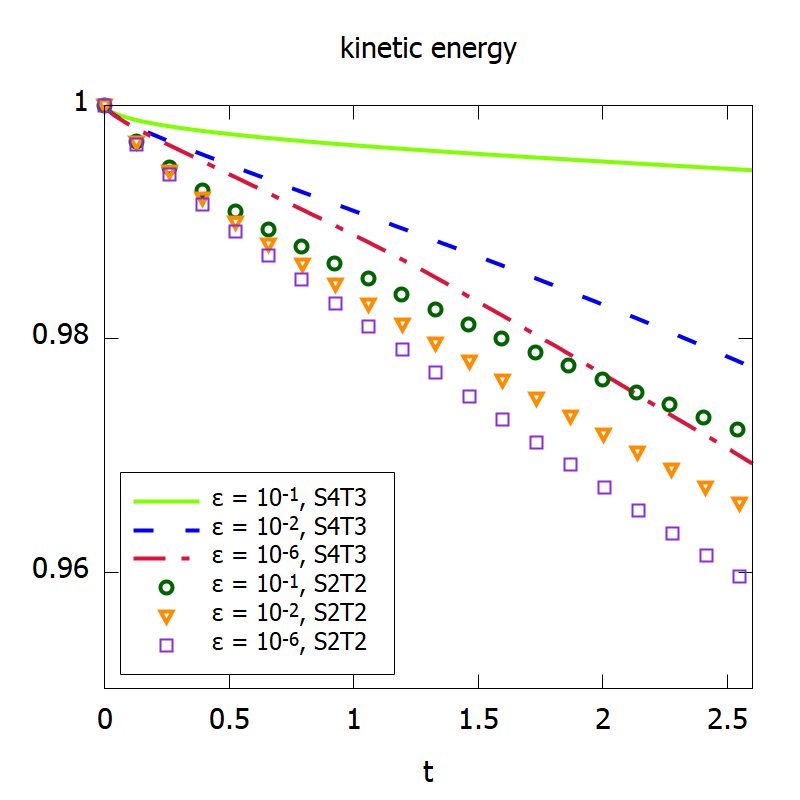}} 
	{\includegraphics[width=0.28\textwidth]{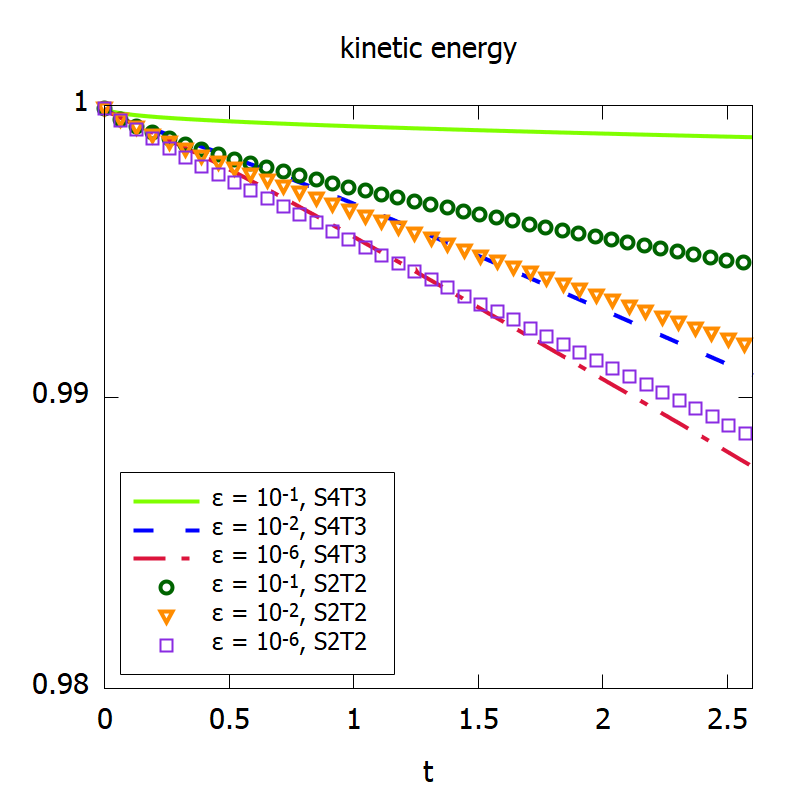}} \\
	{\includegraphics[width=0.28\textwidth]{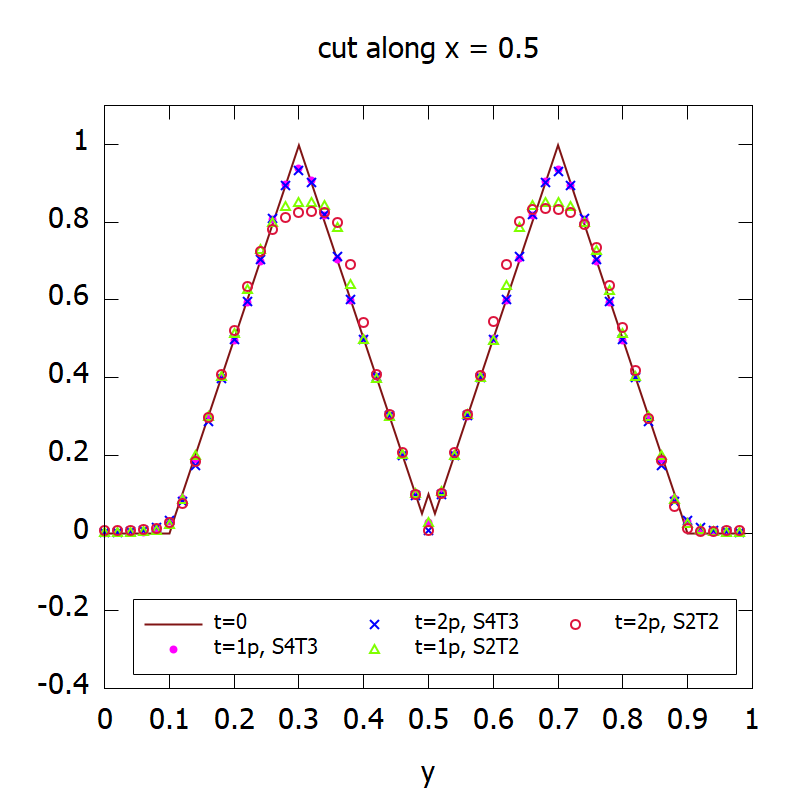}}
	%{\includegraphics[width=0.48\textwidth]{./pic/pic/2D/Gresho/gresho42.png}}	
	{\includegraphics[width=0.28\textwidth]{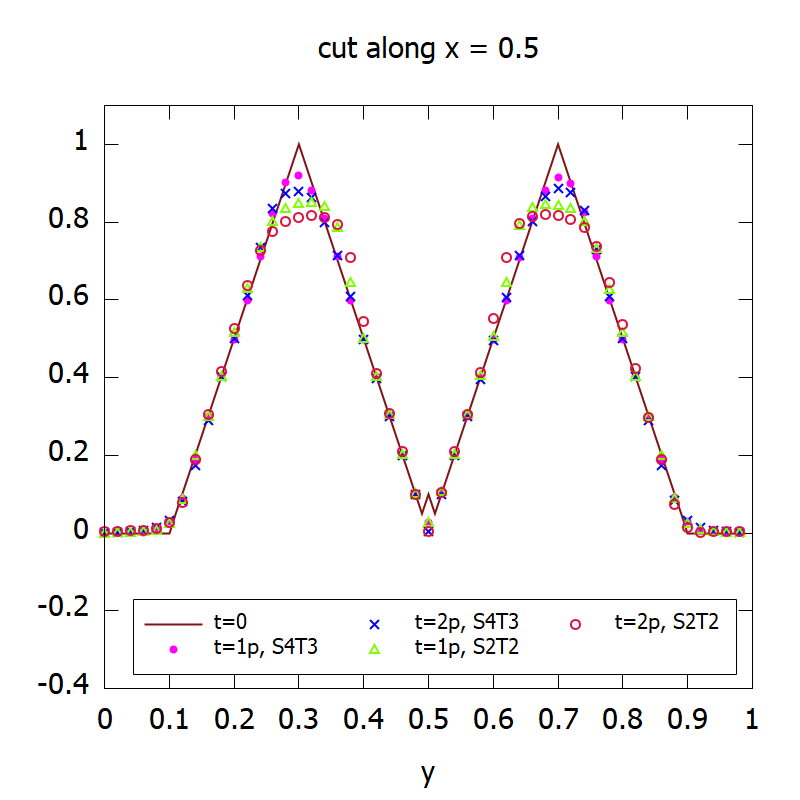}}
	%{\includegraphics[width=0.48\textwidth]{./pic/pic/2D/Gresho/gresho43.png}}	
	{\includegraphics[width=0.28\textwidth]{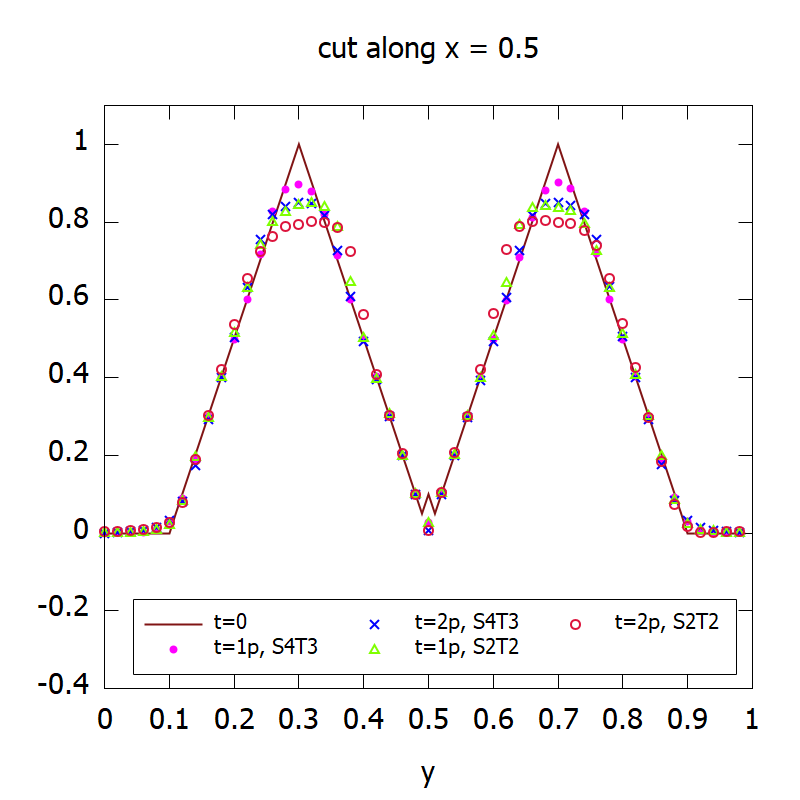}}
	{\includegraphics[width=0.28\textwidth]{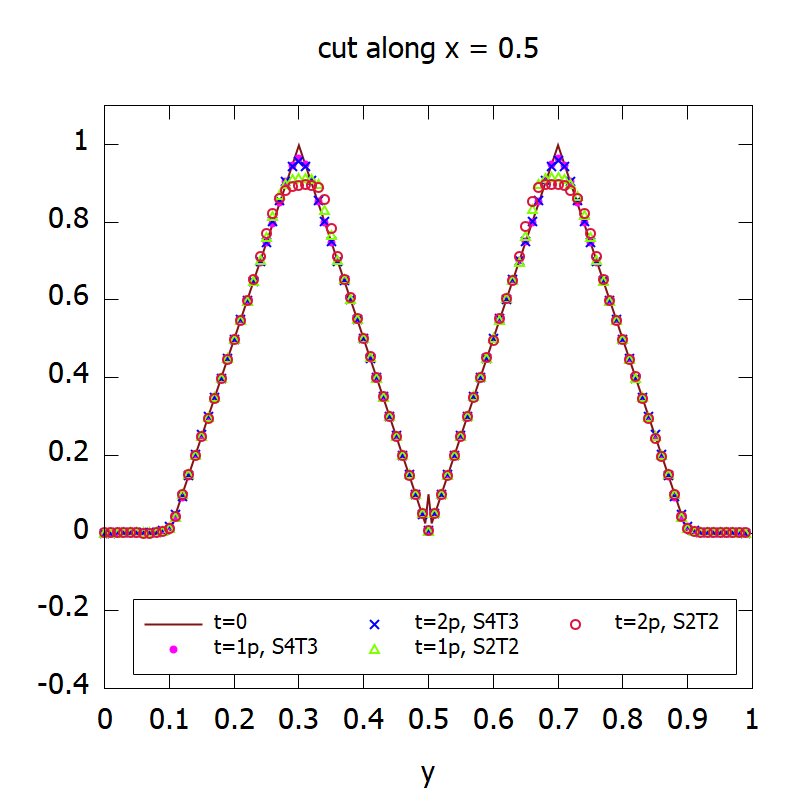}}
%{\includegraphics[width=0.48\textwidth]{./pic/pic/2D/Gresho/gresho42.png}}	
{\includegraphics[width=0.28\textwidth]{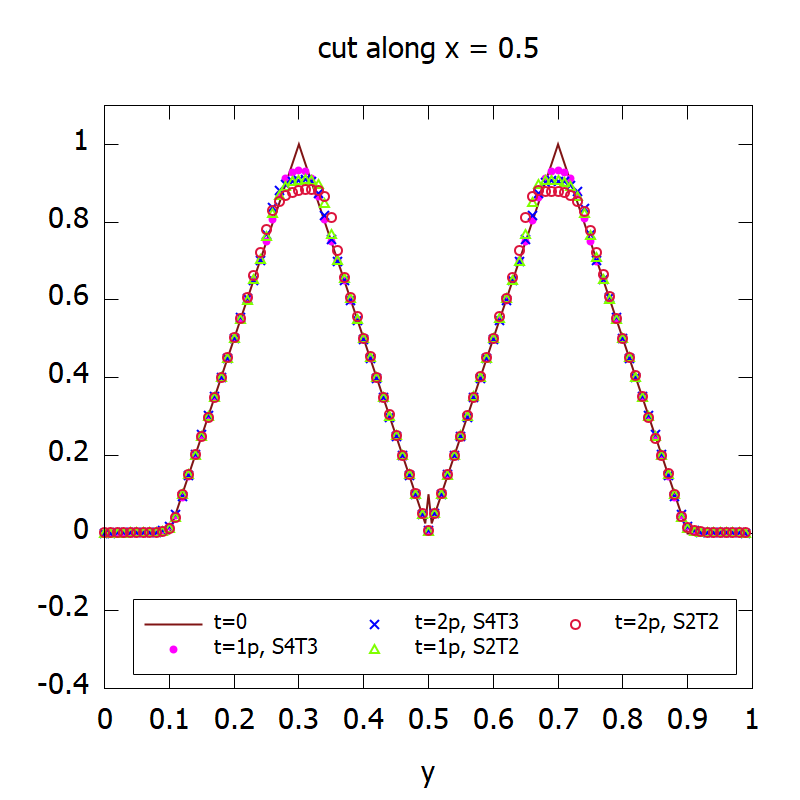}}
%{\includegraphics[width=0.48\textwidth]{./pic/pic/2D/Gresho/gresho43.png}}	
{\includegraphics[width=0.28\textwidth]{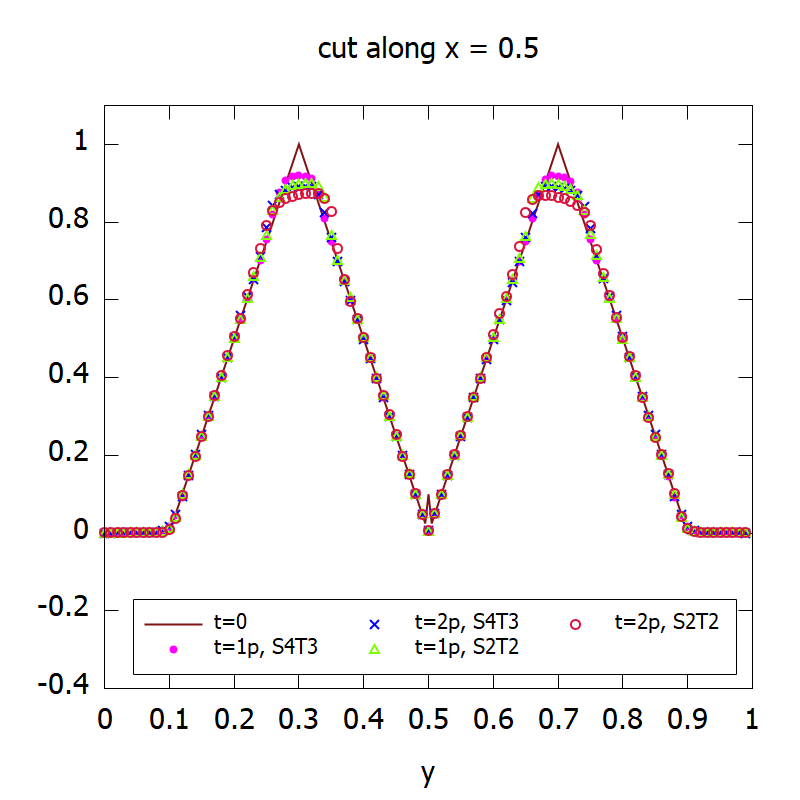}}	
	\caption{Example \ref{gresho}, Gresho vortex problem. Top left: $M_{\text{ratio}}$ in \eqref{rMach} at $T=R\pi$ on the mesh grid $100\times100$.
	{Top middle: time evolution of the kinetic energy \eqref{kie} relative to its initial value on the mesh grid $100\times100$. Top right: time evolution of the kinetic energy \eqref{kie} relative to its initial value on the mesh grid $200\times200$. Middle and bottom:
    cuts along $x=0.5$ at three different times, initial ($t=0$), $1$ period ($t=1 p$), $2$ periods ($t=2 p$). The results are shifted to the center by $-u_\infty t$ periodically. From left to right:  $\eps=10^{-1}$, $\eps=10^{-2}$, $\eps=10^{-6}$. Middle: mesh $100\times 100$; bottom: mesh $200\times200$.}}
	\label{2dgresho}
    \end{figure} 

%    \begin{figure}[!h]
%	\centering
%	{\includegraphics[width=0.24\textwidth]{./pic/pic/2D/Gresho/_kin_eng100.png}} 
%	{\includegraphics[width=0.24\textwidth]{./pic/pic/2D/Gresho/_kin_eng200.png}}
%	\caption{Example \ref{gresho}, Gresho vortex problem. Time evolution of the kinetic energy \eqref{kie} relative to its initial value. Mesh grid: $100\times 100$ (left); $200\times200$ (right).}
%	\label{2dgreshokin}
%    \end{figure} 	

%    \begin{figure}[!h]
%	\centering
%	{\includegraphics[width=0.48\textwidth]{./pic/pic/2D/Gresho/_cutx1.png}} 
%	{\includegraphics[width=0.48\textwidth]{./pic/pic/2D/Gresho/_kin_eng1.png}}\\
%	{\includegraphics[width=0.48\textwidth]{./pic/pic/2D/Gresho/_cutx2.png}} 
%	{\includegraphics[width=0.48\textwidth]{./pic/pic/2D/Gresho/_kin_eng2.png}}\\
%	{\includegraphics[width=0.48\textwidth]{./pic/pic/2D/Gresho/_cutx3.png}} 
%	{\includegraphics[width=0.48\textwidth]{./pic/pic/2D/Gresho/_kin_eng3.png}}
%	\caption{Example \ref{gresho}, Gresho vortex problem. Left: cuts along $x=0.5$, shifted to the center by $-u_\infty t$ periodically, at four different times, initial ($t=0$), $1/4$ period ($t=1/4 p$), $1/2$ period ($t=1/2 p$) and $1$ period ($t=1 p$). Right: time evolution of kinetic energy \eqref{kie}. From top to bottom: $\eps=10^{-1}, 10^{-2}, 10^{-6}$. Mesh grid: $100\times 100$.}
%	\label{2dgresho2}
%    \end{figure}

\end{exa}

\begin{exa} 
\label{2DaccINP_2}
(Incompressible flow.)
Finally we consider two problems in the incompressible flow regime by taking $\eps=10^{-6}$. One is the shear flow problem on $[0,2\pi]^2$ with
\beq \label{shear}
	v(x,y,0) = \delta \cos(x), \quad
	u(x,y,0) = 
	\left\{  
	\begin{array}{ll}
	\tanh((y-\frac{\pi}{2})/\rho), \,\, &\rm{ if } \,\, y \le \pi,  \\[2mm]
	\tanh((\frac{3\pi}{2}-y)/\rho), \,\,& \rm{ if } \,\, y > \pi.
	\end{array}
	\right.
\eeq
The other is the Kelvin-Helmholtz instability problem \cite{crouseilles2010conservative} on $[0,4\pi]\times[0,2\pi]$ with
\beq \label{kh}
u(x,y,0) = \cos(y), \quad v(x,y,0) = 0.03 \sin(0.5 x).
\eeq 
The initial density and pressure for both cases are taken to be $1$.
For the shear flow problem, we run the solution up to $T=6$ on a mesh grid $256\times256$, while $T=40$ for the Kelvin-Helmholtz instability problem on the mesh grid $N_x\times N_y=256\times128$. The vorticity $\omega=v_x-u_y$ for both cases are shown in Fig.~\ref{shear6}, where $v_x$ and $u_y$ are discretized by the 4th order central difference. We observe that it is comparable to the results for the isentropic case as in \cite{boscarino2019high}. 

We also show the time evolution of the divergence error $u_x+v_y$ for the velocity in Fig. \ref{shear6}. The divergence is computed by a 5th order finite difference WENO reconstruction with zero viscosity, which mimics what we have done in the numerical scheme. We also compare it with the linear 4th order central difference discretization. For both cases, the divergence error is increasing with time. When very fine structures are no longer supported by the mesh, the divergence error suffers from a sudden increase. 
We would remark that when the flow is incompressible or weakly incompressible, without discontinuities in the initial condition, high order linear reconstructions would perform better than WENO reconstructions in preserving the divergence and in resolving solution structures. 
%However, our purpose here is designing high order spatial discretizations working in both compressible and incompressible regimes. 

\begin{figure}[!h]
	\centering
	{\includegraphics[width=0.3\textwidth]{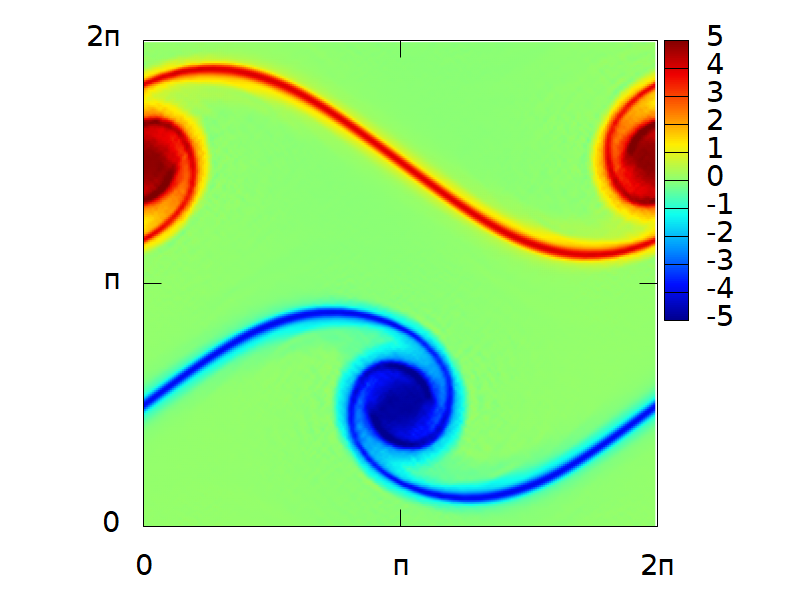}} 
	{\includegraphics[width=0.38\textwidth]{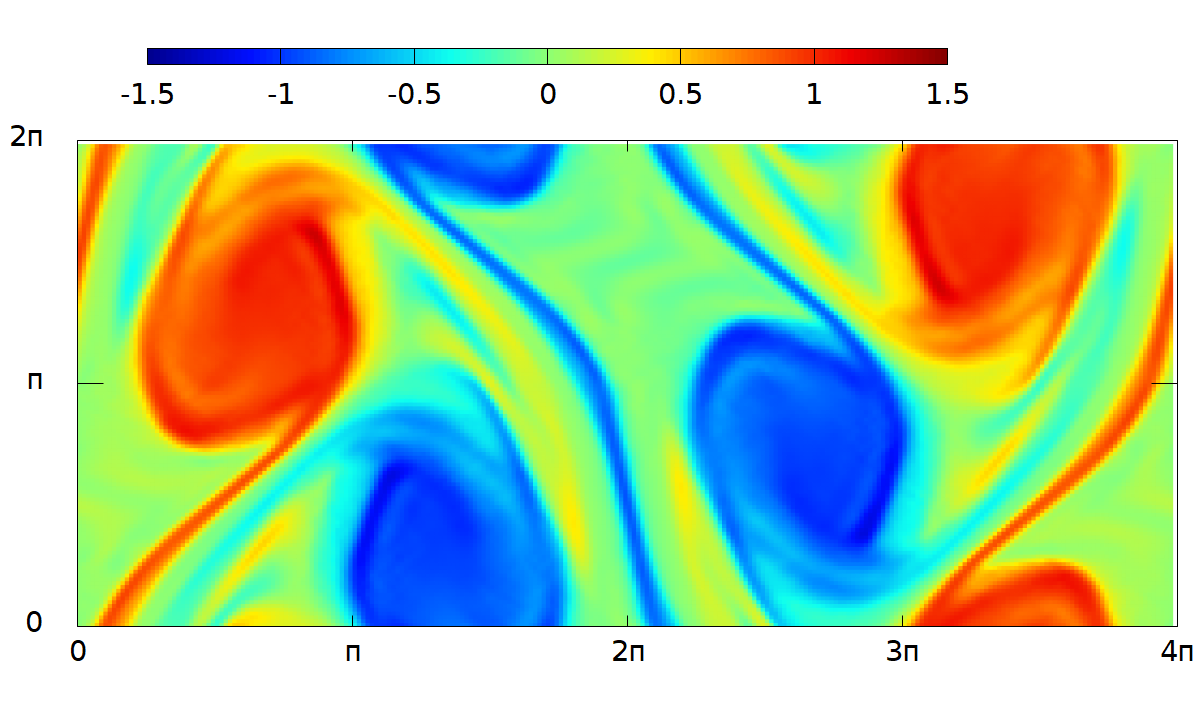}}\\
	{\includegraphics[width=0.30\textwidth, height=0.22\textwidth]{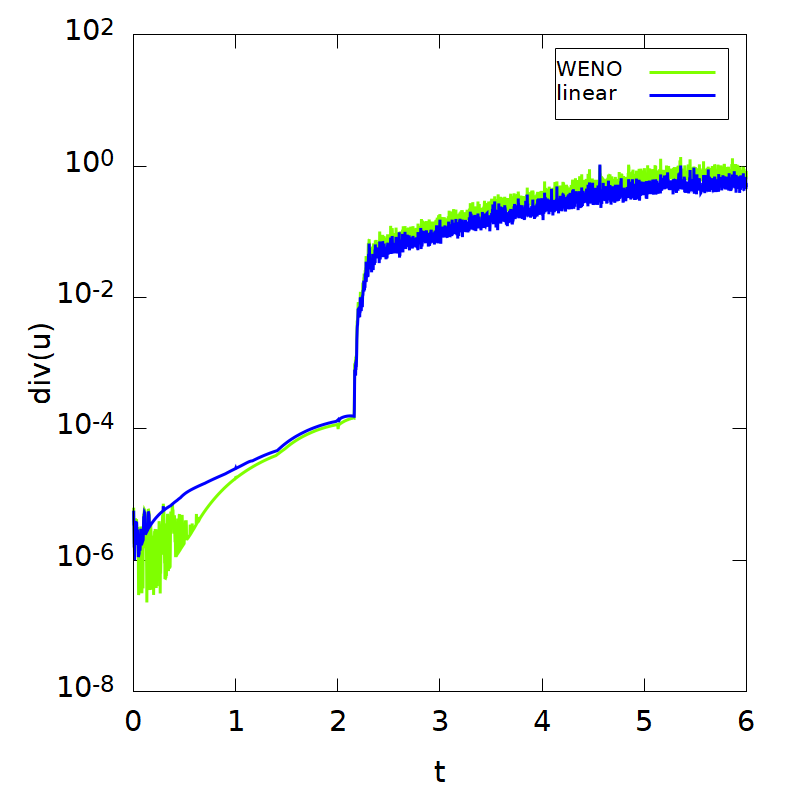}}
	{\includegraphics[width=0.38\textwidth, height=0.22\textwidth]{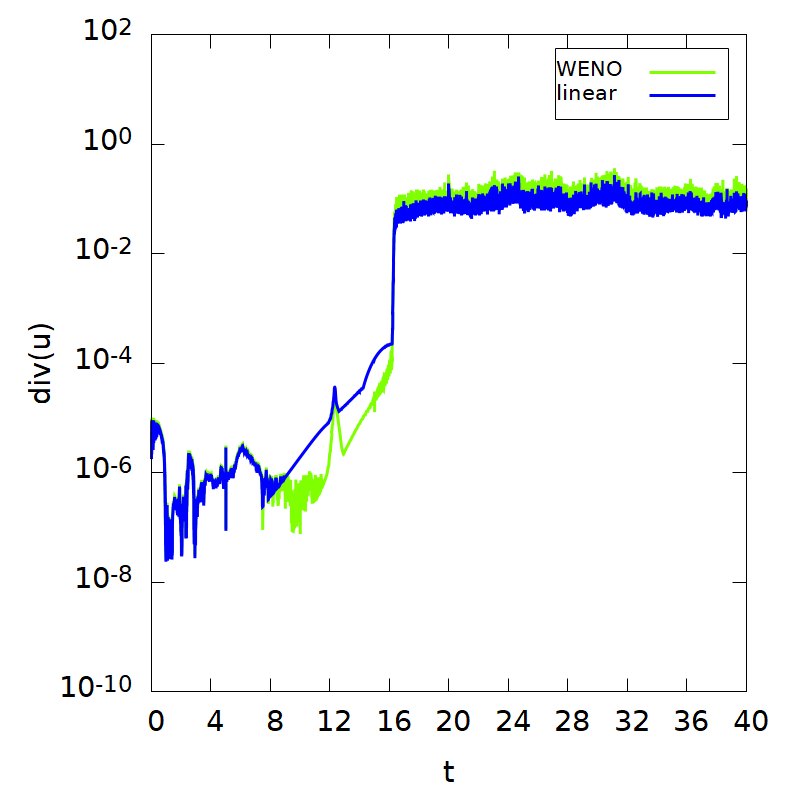}}
	%{\includegraphics[width=0.48\textwidth]{./pic/2D/INP/shear/div_shear6.png}}
	\caption{Example \ref{2DaccINP_2}, the vorticity $\omega=v_x-u_y$ and time history of $L^\infty$ norm for the divergence $u_x+v_y$. Left: the shear flow problem \eqref{shear} at $T=6$ with mesh grid $256\times256$. Right: the Kelvin-Helmholtz instability problem \eqref{kh} at $T=40$ with mesh grid $256\times128$.}
	\label{shear6}
\end{figure} 

%\begin{figure}[!h]
%	\centering
%	{\includegraphics[width=0.24\textwidth]{./pic/2D/INP/shear/div_t.png}}
%	{\includegraphics[width=0.24\textwidth]{./pic/2D/INP/kh/div_t.png}}
%	\caption{Example \ref{2DaccINP_2}. The time history of $L^\infty$ norm for the divergence $u_x+v_y$. ``WENO'' denotes a 5th order finite difference WENO reconstruction with zero viscosity, while ``linear'' denotes a 4th order central difference discretization. Left: shear flow; right: Kelvin-Helmholtz instability. }
%	\label{divt}
%\end{figure} 

\end{exa}
\section{Conclusion}
\label{sec6}

In this paper we present a high-order semi-implicit IMEX RK WENO scheme for the full compressible Euler equations in the case of all-Mach flows. We combine the semi-implicit IMEX RK discretization in time with high order finite difference WENO space discretizations. 
The EOS is treated in a semi-implicit manner, therefore requiring the solution of linearized elliptic systems at each time step. Characteristic reconstructions are proposed in the semi-implicit framework with a fixed splitting for the pressure. The scheme is proven to be asymptotic preserving and asymptotically accurate in the incompressible limit. Numerical tests have demonstrated the efficiency and effectiveness of our proposed approach.

%\section*{Acknowledgement}
%G. Russo and S. Boscarino would like to thank the Italian Ministry of Instruction, University and Research (MIUR) to support this research with funds coming from PRIN Project 2017 (No. 2017KKJP4X entitled ``Innovative numerical methods for evolutionary partial differential equations and applications''). Moreover, G. Russo has been supported by ITN-ETN Horizon 2020 Project ModCompShock, Modeling and Computation on Shocks and Interfaces, Project Reference 642768 and S. Boscarino has been supported by the University of Catania (“Piano della Ricerca 2016/2018, Linea di intervento 2”). S. Boscarino and G. Russo are members of the INdAM Research group GNCS. J.-M. Qiu acknowledges support by NSF grant NSF-DMS-1818924, Air Force Office of Scientific Research FA9550-18-1-0257. T. Xiong acknowledges support by NSFC No. 11971025, NSF of Fujian Province No. 2019J06002, the Strategic Priority Research Program of Chinese Academy of Sciences Grant No. XDA25010401, the Science Challenge Project No. TZ2016002.

%\input{appendix}

\bibliographystyle{siam}
\bibliography{refer}

\newcommand{\noop}[1]{}
\begin{thebibliography}{10}

\bibitem{boscarino2008error}
{\sc S.~Boscarino}, {\em Error analysis of {IMEX} {R}unge-{K}utta methods
  derived from differential-algebraic systems}, SIAM Journal on Numerical
  Analysis, 45 (2008), pp.~1600--1621.

\bibitem{boscarino2015linearly}
{\sc S.~Boscarino, R.~B\"urger, P.~Mulet, G.~Russo, and L.~M. Villada}, {\em
  Linearly implicit {IMEX} {R}unge--{K}utta methods for a class of degenerate
  convection-diffusion problems}, SIAM Journal on Scientific Computing, 37
  (2015), pp.~B305--B331.

\bibitem{boscarino2016linearly}
{\sc S.~Boscarino, R.~B{\"u}rger, P.~Mulet, G.~Russo, and L.~M. Villada}, {\em
  On linearly implicit {IMEX} {R}unge-{K}utta methods for degenerate
  convection-diffusion problems modeling polydisperse sedimentation}, Bulletin
  of the Brazilian Mathematical Society, New Series, 47 (2016), pp.~171--185.

\bibitem{boscarino2016high}
{\sc S.~Boscarino, F.~Filbet, and G.~Russo}, {\em High order semi-implicit
  schemes for time dependent partial differential equations}, Journal of
  Scientific Computing, 68 (2016), pp.~975--1001.

\bibitem{boscarino2014high}
{\sc S.~Boscarino, P.~G. LeFloch, and G.~Russo}, {\em {High-order
  asymptotic-preserving methods for fully nonlinear relaxation problems}}, SIAM
  Journal on Scientific Computing, 36 (2014), pp.~A377--A395.

\bibitem{boscarino2017asymptotic}
{\sc S.~Boscarino and L.~Pareschi}, {\em {On the asymptotic properties of IMEX
  Runge--Kutta schemes for hyperbolic balance laws}}, Journal of Computational
  and Applied Mathematics, 316 (2017), pp.~60--73.

\bibitem{boscarino2013implicit}
{\sc S.~Boscarino, L.~Pareschi, and G.~Russo}, {\em {Implicit-explicit
  Runge--Kutta schemes for hyperbolic systems and kinetic equations in the
  diffusion limit}}, SIAM Journal on Scientific Computing, 35 (2013),
  pp.~A22--A51.

\bibitem{boscarino2019high}
{\sc S.~Boscarino, J.-M. Qiu, G.~Russo, and T.~Xiong}, {\em {A high order
  semi-implicit IMEX WENO scheme for the all-Mach isentropic Euler system}},
  Journal of Computational Physics, 392 (2019), pp.~594--618.

\bibitem{BoRuSc18}
{\sc S.~Boscarino, G.~Russo, and L.~Scandurra}, {\em All mach number second
  order semi-implicit scheme for the {Euler} equations of gas dynamics},
  Journal of Scientific Computing,  (2017), pp.~1--35.

\bibitem{boscheri2021high}
{\sc W.~Boscheri and L.~Pareschi}, {\em High order pressure-based semi-implicit
  {IMEX} schemes for the 3{D} {N}avier-{S}tokes equations at all {M}ach
  numbers}, Journal of Computational Physics, 434 (2021), p.~110206.

\bibitem{butcher2016}
{\sc J.~C. Butcher}, {\em Numerical Methods for Ordinary Differential
  Equations}, Third Edition, John Wiley \& Sons Ltd, 2016.

\bibitem{casulli1984pressure}
{\sc V.~Casulli and D.~Greenspan}, {\em Pressure method for the numerical
  solution of transient, compressible fluid flows}, International Journal for
  Numerical Methods in Fluids, 4 (1984), pp.~1001--1012.

\bibitem{chorin1968numerical}
{\sc A.~J. Chorin}, {\em Numerical solution of the {Navier-Stokes} equations},
  Mathematics of computation, 22 (1968), pp.~745--762.

\bibitem{chorin1993}
{\sc A.~J. Chorin and J.~Marsden}, {\em A Mathematical Introduction to Fluid
  Mechanics}, Springer Verlag, 1993.

\bibitem{cockburn2006advanced}
{\sc B.~Cockburn, C.~Johnson, C.-W. Shu, and E.~Tadmor}, {\em Advanced
  numerical approximation of nonlinear hyperbolic equations: lectures given at
  the 2nd session of the Centro Internazionale Matematico Estivo (CIME) held in
  Cetraro, Italy, June 23-28, 1997}, Springer, 2006.

\bibitem{cordier2012asymptotic}
{\sc F.~Cordier, P.~Degond, and A.~Kumbaro}, {\em An asymptotic-preserving
  all-speed scheme for the {Euler} and {Navier-Stokes} equations}, Journal of
  Computational Physics, 231 (2012), pp.~5685--5704.

\bibitem{crouseilles2010conservative}
{\sc N.~Crouseilles, M.~Mehrenberger, and E.~Sonnendr{\"u}cker}, {\em
  Conservative semi-lagrangian schemes for vlasov equations}, Journal of
  Computational Physics, 229 (2010), pp.~1927--1953.

\bibitem{degond2011all}
{\sc P.~Degond and M.~Tang}, {\em All speed scheme for the low {Mach} number
  limit of the isentropic {Euler} equations}, Communications in Computational
  Physics, 10 (2011), pp.~1--31.

\bibitem{dellacherie2010analysis}
{\sc S.~Dellacherie}, {\em Analysis of {Godunov} type schemes applied to the
  compressible {Euler} system at low {Mach} number}, Journal of Computational
  Physics, 229 (2010), pp.~978--1016.

\bibitem{Denner2020}
{\sc F.~Denner, F.~Evrard, and B.~van Wachem}, {\em Conservative finite-volume
  framework and pressure-based algorithm for flows of incompressible, ideal-gas
  and real-gas fluids at all speeds}, Journal of Computational Physics, 409
  (2020), p.~109348.

\bibitem{Denner2018}
{\sc F.~Denner, C.-N. Xiao, and B.~van Wachem}, {\em Pressure-based algorithm
  for compressible interfacial flows with acoustically-conservative interface
  discretization}, Journal of Computational Physics, 367 (2018), pp.~192--234.

\bibitem{dimarco2017study}
{\sc G.~Dimarco, R.~Loub{\`e}re, and M.-H. Vignal}, {\em Study of a new
  asymptotic preserving scheme for the {Euler} system in the low {Mach} number
  limit}, SIAM Journal on Scientific Computing, 39 (2017), pp.~A2099--A2128.

\bibitem{GoRa2014}
{\sc E.~Godlewski and P.-A. Raviart}, {\em Numerical Approximation of
  Hyperbolic Systems of Conservation Laws}, Springer, 2014.

\bibitem{guermond2006overview}
{\sc J.-L. Guermond, P.~Minev, and J.~Shen}, {\em An overview of projection
  methods for incompressible flows}, Computer Methods in Applied Mechanics and
  Engineering, 195 (2006), pp.~6011--6045.

\bibitem{haack2012all}
{\sc J.~Haack, S.~Jin, and J.-G. Liu}, {\em An all-speed asymptotic-preserving
  method for the isentropic {Euler} and {Navier-Stokes} equations},
  Communications in Computational Physics, 12 (2012), pp.~955--980.

\bibitem{hairer1993solving1}
{\sc E.~Hairer, S.~N{\o}rsett, and G.~Wanner}, {\em Solving ordinary
  differential equations: Nonstiff problems}, vol.~1, Springer Verlag, 1993.

\bibitem{hairer1993solving2}
{\sc E.~Hairer and G.~Wanner}, {\em Solving ordinary differential equations II:
  stiff and differential algebraic problems}, vol.~2, Springer Verlag, 1993.

\bibitem{Harlow1968}
{\sc F.~Harlow and A.~Amdsen}, {\em Numerical calculation of almost
  incompressible flow}, Journal of Computational Physics, 3 (1968), pp.~80--93.

\bibitem{Harlow1971}
\leavevmode\vrule height 2pt depth -1.6pt width 23pt, {\em A numerical fluid
  dynamics calculation method for all flow speeds}, Journal of Computational
  Physics, 8 (1971), pp.~197--213.

\bibitem{jin1999efficient}
{\sc S.~Jin}, {\em Efficient asymptotic-preserving ({AP}) schemes for some
  multiscale kinetic equations}, SIAM Journal on Scientific Computing, 21
  (1999), pp.~441--454.

\bibitem{jin2010asymptotic}
\leavevmode\vrule height 2pt depth -1.6pt width 23pt, {\em Asymptotic
  preserving {(AP)} schemes for multiscale kinetic and hyperbolic equations: a
  review}, Lecture Notes for Summer School on Methods and Models of Kinetic
  Theory (M\&MKT), Porto Ercole (Grosseto, Italy),  (2010), pp.~177--216.

\bibitem{klainerman1981singular}
{\sc S.~Klainerman and A.~Majda}, {\em Singular limits of quasilinear
  hyperbolic systems with large parameters and the incompressible limit of
  compressible fluids}, Communications on pure and applied Mathematics, 34
  (1981), pp.~481--524.

\bibitem{klainerman1982compressible}
\leavevmode\vrule height 2pt depth -1.6pt width 23pt, {\em Compressible and
  incompressible fluids}, Communications on Pure and Applied Mathematics, 35
  (1982), pp.~629--651.

\bibitem{Klein1995}
{\sc R.~Klein}, {\em Semi-implicit extension of a {G}odunov-type scheme based
  on low {M}ach number asymptotics. {I}: One-dimensional flow}, J. Comput.
  Phys., Vol. 121 (1995), pp.~pp. 213--237.

\bibitem{lax1998}
{\sc P.~D. Lax and X.-D. Liu}, {\em Solution of two-dimensional {Riemann}
  problems of gas dynamics by positive schemes}, SIAM Journal on Scientific
  Computing, 19 (1998), pp.~319--340.

\bibitem{leveque2002finite}
{\sc R.~J. LeVeque}, {\em Finite volume methods for hyperbolic problems},
  vol.~31, Cambridge university press, 2002.

\bibitem{Ropke2015Mach}
{\sc F.~Miczek, F.~R\"opke, and P.~Edelmann}, {\em A new numerical solver for
  flows at various {M}ach numbers}, Astronomy \& Astrophysics, Vol. 576 (2015),
  p.~A50.

\bibitem{munz2003extension-low-Mach}
{\sc C.-D. Munz, S.~Roller, R.~Klein, and K.~J. Geratz}, {\em The extension of
  incompressible flow solvers to the weakly compressible regime}, Computers \&
  Fluids, 32 (2003), pp.~173--196.

\bibitem{noelle2014weakly}
{\sc S.~Noelle, G.~Bispen, K.~R. Arun, M.~Luk\'{a}\v{c}ov\'{a}-Medvidov\'{a},
  and C.-D. Munz}, {\em {A weakly asymptotic preserving low Mach number scheme
  for the Euler equations of gas dynamics}}, SIAM Journal on Scientific
  Computing, 36 (2014), pp.~B989--B1024.

\bibitem{pareschi2005implicit}
{\sc L.~Pareschi and G.~Russo}, {\em {Implicit-explicit Runge-Kutta schemes and
  applications to hyperbolic systems with relaxation}}, Journal of Scientific
  computing, 25 (2005), pp.~129--155.

\bibitem{park-munz-2005multiple-all-Mach}
{\sc J.-H. Park and C.-D. Munz}, {\em Multiple pressure variables methods for
  fluid flow at all {M}ach numbers}, International Journal for Numerical
  Methods in Fluids, 49 (2005), pp.~905--931.

\bibitem{shu1998essentially}
{\sc C.-W. Shu}, {\em Essentially non-oscillatory and weighted essentially
  non-oscillatory schemes for hyperbolic conservation laws}, Advanced Numerical
  Approximation of Nonlinear Hyperbolic Equations,  (1998), pp.~325--432.

\bibitem{shu2009high}
\leavevmode\vrule height 2pt depth -1.6pt width 23pt, {\em High order weighted
  essentially nonoscillatory schemes for convection dominated problems}, SIAM
  review, 51 (2009), pp.~82--126.

\bibitem{tang2012}
{\sc M.~Tang}, {\em Second order method for isentropic euler equation in the
  low mach number regime}, Kinetic and Related Models, 5 (2012), pp.~155--184.

\bibitem{tavelli2017pressure}
{\sc M.~Tavelli and M.~Dumbser}, {\em A pressure-based semi-implicit
  space--time discontinuous galerkin method on staggered unstructured meshes
  for the solution of the compressible navier--stokes equations at all mach
  numbers}, Journal of Computational Physics, 341 (2017), pp.~341--376.

\bibitem{Temam1984}
{\sc R.~Temam}, {\em Navier-Stokes Equations: Theory and Numerical Analysis},
  AMS Chelsea Publishing, 1984.

\bibitem{toro2009riemann}
{\sc E.~F. Toro}, {\em Riemann solvers and numerical methods for fluid
  dynamics: a practical introduction}, Springer, 2009.

\bibitem{Turkel1987}
{\sc E.~Turkel}, {\em Preconditioned methods for solving the incompressible and
  low speed compressible equations}, Journal of Computational Physics, 72
  (1987), pp.~277--298.

\bibitem{Turkel1993}
{\sc E.~Turkel, A.~Fiterman, and B.~van Leer}, {\em Preconditioning and the
  limit to the incompressible flow equations}, NASA CR-191500,  (1993).

\bibitem{viozat1997implicit}
{\sc C.~Viozat}, {\em {I}mplicit upwind schemes for low {M}ach number
  compressible flows}, PhD thesis, Inria, 1997.

\bibitem{Jonas2020novel}
{\sc J.~Zeifang, J.~Sch\"{u}tz, K.~Kaiser, A.~Beck,
  M.~Luk\'a\v{c}ov\'a-Medvidov\'a, and S.~Noelle}, {\em {A novel full-Euler low
  Mach number IMEX splitting}}, Communication in Computational Physics, 27
  (2020), pp.~292--320.

\end{thebibliography}
%\bibliography{refer2}

\end{document}